\documentclass[reqno, 11pt]{amsart}
\usepackage{amsmath, amsthm, amscd, amsfonts, amssymb, graphicx, color, mathrsfs, enumitem}
\usepackage{extarrows}
\usepackage[bookmarksnumbered, colorlinks, plainpages]{hyperref}
\hypersetup{colorlinks=true,linkcolor=red, anchorcolor=green,
	citecolor=cyan, urlcolor=red, filecolor=magenta, pdftoolbar=true}

\newtheorem{theorem}{Theorem}[section]
\newtheorem{lemma}[theorem]{Lemma}
\newtheorem{proposition}[theorem]{Proposition}

\theoremstyle{definition}

\newtheorem{example}[theorem]{Example}

\theoremstyle{remark}
\newtheorem{remark}[theorem]{Remark}
\numberwithin{equation}{section}

\begin{document}

	\setcounter{page}{1}
	
	\title[Sum commuting maps with the Aluthge transform]
	{Weighted sum commuting maps with the $(f,g)$-Aluthge transform}

	\author[J.B. Zhou]{Jing-Bin Zhou}
	\address{Jing-Bin Zhou, School of Mathematics,
		Shanghai University of Finance and Economics,
		777 Gouding Road, Shanghai 200433, P. R. China}
	\email{zhoujingbin@stu.sufe.edu.cn}
	
	\author[S. Yang]{Shihai Yang\textsuperscript{*}}
	\address{Shihai Yang, School of Mathematics,
		Shanghai University of Finance and Economics,
		777 Gouding Road, Shanghai 200433, P. R. China}
	\email{yang.shihai@mail.shufe.edu.cn}

	
	
	\begin{abstract}
			Let $\mathbb{M}_n(\mathbb{C})$ denote the
			algebra of all $n\times n$ complex matrices,
			where $n\geq2$, and let $\Delta_{f,g}$
			be the $(f,g)$-Aluthge transform associated
			with nonnegative continuous functions
			$f,g$ on $[0,\infty)$ satisfying
			$f(t)g(t)=t$ for all $t\geq0$.
			In this paper, we characterize all surjective
			maps $\Phi:\mathbb{M}_n(\mathbb{C})
			\to\mathbb{M}_n(\mathbb{C})$ satisfying
			\[
			\Phi(\Delta_{f,g}(\alpha A+\beta B))
			=\Delta_{f,g}(\alpha\Phi(A)+\beta\Phi(B))
			\]
			for all $A,B\in\mathbb{M}_n(\mathbb{C})$,
			where $\alpha,\beta$ are fixed nonzero
			real numbers.
			Without assuming additivity or linearity,
			we prove that when $f$ is nonconstant,
			$\Phi$ must be of the form
			$\Phi(A)=cUAU^*$ or
			$\Phi(A)=cU\overline{A}U^*$
			for all $A$, where $U$ is a fixed unitary
			matrix, $\overline{A}$ denotes the
			entrywise complex conjugate of $A$,
			and $c$ is a nonzero complex scalar. Moreover, we completely determine
			the possible values of $|c|$,
			which form either $\{1\}$,
			$q^{\mathbb{Z}}$ for some $q>1$,
			or $(0,\infty)$. When $f$ is constant, we show that
			$\Phi(A)=S(A)+D$ for all $A$, where
			$S$ is a surjective
			$\mathbb{Q}(\alpha,\beta)$-linear map
			and $D\in\mathbb{M}_n(\mathbb{C})$
			satisfies $(\alpha+\beta-1)D=0$.
	
	\end{abstract}
	\maketitle
	\renewcommand{\thefootnote}{*}
	\footnotetext{Corresponding author: Shihai Yang}
	
	\section{Introduction}
The notation $\mathbb{B}(\mathcal{H})$ is used for the algebra of all bounded linear operators on a complex Hilbert space $\mathcal{H}$. Denote the kernel of an operator $T\in\mathbb{B}(\mathcal{H})$ by $\ker T$. Let $T=U|T|$ be the polar decomposition of $T$, where $|T|=(T^\ast T)^\frac{1}{2}$ is the absolute value of $T\in\mathbb{B}(\mathcal{H})$ and $U$ is the partial isometry satisfying
$\ker U=\ker T$. The operator $\Delta(T)=|T|^\frac{1}{2}U|T|^\frac{1}{2}$ was introduced in \cite{Aluthge-IEOT-1990} to study properties of hyponormal operators and was later called the {\em Aluthge transform} of $T$. Over the years, one of the reasons why
the operator theory of the Aluthge transform
has attracted the attention of many authors
may be its connection with the invariant
subspace problem. Specifically, Jung et al.
proved that $T\in\mathbb{B}(\mathcal{H})$
has a nontrivial invariant subspace if
and only if $\Delta(T)$ does
\cite{Jung-IEOT-2000}.
Subsequently, the same authors also showed
that $T$ and $\Delta(T)$ have the same
spectrum \cite{Jung-IEOT-2001}.
Kim and Ko proved that $T$ has property
$(\beta)$ if and only if $\Delta(T)$
has property $(\beta)$, and that $T$
has the single-valued extension property
(SVEP) if and only if $\Delta(T)$ does
\cite{Kim-GMJ-2005}.
In other words, the Aluthge transform
preserves many properties of the original
operator. Moreover, another noteworthy
contribution is due to Yamazaki, who
established a new spectral radius formula
through iterations of the Aluthge transform
\cite{Yamazaki-PAMS-2002}. Other developments concerning the Aluthge transform
include various extensions of this transform.
In \cite{Cho-Tanahashi-SMJ-2002},
Ch\={o} and Tanahashi generalized this transform
and called the operator
$\Delta_{\lambda}(T)=|T|^\lambda U|T|^{1-\lambda}$
the {\em $\lambda$-Aluthge transform} of $T$,
where $\lambda\in[0,1]$. It should also be mentioned that, for $\lambda=1$, the operator $\Delta_{D}(T)\triangleq\Delta_{1}(T)$ is called the {\em Duggal transform} of $T$ in \cite{Foias-PJM-2003}. While improving estimates for the numerical
radius of operators, Shebrawi and Bakherad
further extended this notion in
\cite{Shebrawi-Filomat-2018} and defined
the operator
$\Delta_{f,g}(T)=f(|T|)Ug(|T|)$
as the {\em $(f,g)$-Aluthge transform}
of $T$, where
\[
C_{\Delta}^{+}([0,+\infty))
=
\left\{
(f,g)\in C^{+}([0,+\infty))
\times C^{+}([0,+\infty)):
f(x)g(x)=x\text{ for all }x\geq0
\right\},
\]
and $C^{+}([0,+\infty))$ denotes
the set of all nonnegative continuous
functions on $[0,+\infty)$. In particular, for $\lambda\in[0,1]$, if $f(x)=x^\lambda$ and $g(x)=x^{1-\lambda}$, then $\Delta_{f,g}(T)$ reduces to $\Delta_{\lambda}(T)$. Recently, the first author extended this transform to the multivariable case and also obtained some new properties of the $(f,g)$-Aluthge transform in \cite{Zhou-glma-2023}. For more results on the Aluthge transform, we refer the reader to \cite{Antezana-AvM-2011, Yamazaki-LAA-2021, Yamazak-Trans-i2025, Garcia-IEOT-2008, Labbane-CAOT-2025}.

Another important direction in functional
analysis and operator theory is the study
of preserver problems. A central objective
in such problems is to determine the
structure of all maps satisfying the
prescribed conditions. One of the earliest studies of preserver
problems related to operator transforms
was carried out by Botelho, Molnár and
Nagy \cite{Botelho-bull-2016}, who
characterized bijective linear maps
between von Neumann factors commuting
with the $\lambda$-Aluthge transform. Subsequently, Chabbabi
\cite{Chabbabi-jmaa-2017} characterized
nonlinear bijective maps on
$\mathbb{B}(\mathcal{H})$ commuting
with the $\lambda$-Aluthge transform
under the usual operator product. In \cite{Chabbabi-jmaa-20172},
Chabbabi and Mbekhta continued
the investigation by replacing
the usual operator product
with the Jordan product.
Motivated by these studies,
Ben Ali Essaleh and Peralta
\cite{Ben-laa-2018}
investigated maps commuting with
the $\lambda$-Aluthge transform
under Jordan-type products
on von Neumann algebras. More recently, Bourhim and Mbekhta
\cite{Bourhim-jmaa-2025}
characterized bijective linear maps
preserving the Aluthge transforms
of unitarily similar operators. Labbane and Aharmim
\cite{Labbane-lma-2024, Labbane-lma-2025}
characterized bijective maps
commuting with the $\lambda$-Aluthge
transform under certain generalized
Jordan-type products.
Recently, they further characterized
maps commuting with the
$(f,g)$-Aluthge transform
under triple and generalized
triple Lie-star Jordan products
\cite{Labbane-rcmp-2026}. To the best of our knowledge,
this is among the few existing works
devoted specifically to preserver
problems for the $(f,g)$-Aluthge transform. For other related preserver problems
involving operator transforms,
we refer the reader to
\cite{Bourhim-lma-2025,
	Bourhim-lma-2026,
	Bourhim-lma-20261}.

Let $\mathbb{N}$, $\mathbb{Z}$, $\mathbb{Q}$,
and $\mathbb{R}$ denote the sets of positive
integers, integers, rational numbers, and
real numbers, respectively, and let
$\mathbb{C}$ denote the complex field. As mentioned above, most existing
studies have focused on characterizing
maps commuting with the Aluthge
transform under various generalized
products of operators, whereas
relatively little attention has
been paid to maps commuting with
the Aluthge transform in connection
with sums of operators. Fix $\alpha,\beta\in\mathbb{R}\setminus\{0\}$.
In this paper, we study surjective maps
$\Phi:M_n(\mathbb{C})\to M_n(\mathbb{C})$
satisfying
\begin{equation}\label{eq:weighted}
	\Phi(\Delta_{f,g}(\alpha A+\beta B))
	=\Delta_{f,g}(\alpha\Phi(A)+\beta\Phi(B))
\end{equation}
for all $A,B\in\mathbb{M}_n(\mathbb{C})$. Given $T\in\mathbb{M}_n(\mathbb{C})$, let $T^*$, $T^{\mathrm{tr}}$, and $\overline{T}$ stand for the adjoint, transpose, and complex conjugate of $T$, respectively. We define
$$G_{f,g}:=\left\{r>0:
\Delta_{f,g}(rT)=r\Delta_{f,g}(T)
\text{ for all }T\in\mathbb{M}_n(\mathbb{C})
\right\}.$$ Let $\mathbb{Q}(\alpha,\beta)$ denote
the field generated by $\alpha$ and $\beta$
over $\mathbb{Q}$, that is,
\[
\mathbb{Q}(\alpha,\beta)
=\left\{
\frac{p(\alpha,\beta)}{q(\alpha,\beta)}
:\ p,q\in\mathbb{Q}[x,y],\
q(\alpha,\beta)\neq0
\right\}.
\]
The following theorem is our main result.

\begin{theorem}\label{thm:main}
	Let $n\geq2$, $(f,g)\in C_{\Delta}^{+}([0,\infty))$,
	and $\alpha,\beta\in\mathbb{R}\setminus\{0\}$.
	A surjective map
	$\Phi:\mathbb{M}_n(\mathbb{C})\to\mathbb{M}_n(\mathbb{C})$
	satisfies \eqref{eq:weighted} if and only if one of the following holds:
	\begin{enumerate}[label=(\roman*)]
		\item $f$ is nonconstant, and there exist
		a unitary matrix $U$ and a nonzero scalar
		$c\in\mathbb{C}$ with $|c|\in G_{f,g}$
		such that either
		\[
		\Phi(A)=cUAU^*
		\quad\text{or}\quad
		\Phi(A)=cU\overline{A}U^*
		\]
		for all $A\in\mathbb{M}_n(\mathbb{C})$.
		Moreover, $G_{f,g}$ is one of
		\[
		\{1\},\quad
		\{\cdots,q^{-2},q^{-1},1,q,q^2,\cdots\}\text{ for some }q>1,
		\quad\text{or}\quad
		(0,\infty).
		\]
		
		\item $f$ is constant, and there exist
		a surjective $\mathbb{Q}(\alpha,\beta)$-linear
		map $S$ on $\mathbb{M}_n(\mathbb{C})$
		and a matrix $D\in\mathbb{M}_n(\mathbb{C})$
		with $(\alpha+\beta-1)D=0$ such that
		\[
		\Phi(A)=S(A)+D
		\qquad(A\in\mathbb{M}_n(\mathbb{C})).
		\]
	\end{enumerate}
\end{theorem}

We remark that when $n=1$, every surjective
map $\Phi$ satisfying \eqref{eq:weighted}
has the same form as described in (ii)
of the theorem, since it is easy to see
that $\Delta_{f,g}(z)=z$ for all
$z\in\mathbb{C}$.
For further details, see the proof of
Theorem \ref{main}.
Moreover, Example \ref{ex} (iv) also
provides, when $n=1$, a map $\Phi$
satisfying \eqref{eq:weighted} that is
neither real-linear nor multiplicative.
Therefore, \textbf{unless otherwise specified,
	throughout the paper we always consider
	$\mathbb{M}_n(\mathbb{C})$ with $n\geq2$.} 
	
	Since $\Phi$ is assumed only to be
surjective and to satisfy \eqref{eq:weighted},
without any additional assumptions of
additivity, linearity, or bijectivity,
and since the functions $f$ and $g$
are not explicitly specified,
the analysis presents considerable
difficulties. Consequently, a substantial
part of the paper is devoted to proving
that such a map $\Phi$ is necessarily
additive and satisfies \eqref{eq:commuting}. The paper is organized as follows.
	In Section 2, we establish some new
	properties of the $(f,g)$-Aluthge transform
	that will be needed in the subsequent sections.
	In Section 3, we characterize surjective maps
	$\Phi:\mathbb{M}_n(\mathbb{C})
	\to\mathbb{M}_n(\mathbb{C})$ under the additional assumptions
	that $\Phi$ is additive and satisfies
	\begin{equation}\label{eq:commuting}
		\Phi\circ\Delta_{f,g}
		=\Delta_{f,g}\circ\Phi.
	\end{equation}
	The main purpose of this section is
	to apply \v{S}emrl's result
	\cite[Corollary 2]{Semrl-laa-1993}
	to obtain certain standard forms of $\Phi$,
	and then use the additional assumptions
	to determine its precise form. In Section 4, by applying Theorem 2.10
	of Costara et al. \cite{Costara-laa-2023},
	we show that every surjective map
	$\Phi$ satisfying \eqref{eq:weighted}
	is additive.
	We then prove that such a map
	also satisfies \eqref{eq:commuting}.
	Consequently, by applying
	Theorem \ref{kkeylem} established
	in Section 3, we obtain our main result,
	Theorem \ref{thm:main}.

	\section{New Properties for the $(f,g)$-Aluthge transform of operators}
	
	We denote by $\mathbb{C}^{n}$ the usual $n$-dimensional complex vector space with the standard inner product $\langle\cdot,\cdot\rangle$, and by $\|\cdot\|$ the norm defined by $\|x\|=\sqrt{\langle x,x\rangle}$.
For $T\in \mathbb{M}_n(\mathbb{C})$, let $\sigma(T)$, $\operatorname{tr}(T)$, $\|T\|$, $\det (T)$, and $\operatorname{ran} T$ denote its spectrum, trace, operator norm, determinant, and range, respectively. The Frobenius
norm of $T$ is given by $\|T\|_F=\sqrt{\operatorname{tr}(T^*T)}$. We denote by $I$ the identity matrix in
$\mathbb{M}_n(\mathbb{C})$, and by
$\operatorname{diag}(\lambda_1,\lambda_{2},\ldots,\lambda_n)$
the diagonal matrix with diagonal entries
$\lambda_1,\lambda_{2},\ldots,\lambda_n$, where $\lambda_{i}\in\mathbb{C}$ for $i=1,2,\cdots,n$. 

As will be needed later, we first establish in this section some further properties of this transform. The results in the proposition below are elementary. Note that for $(f,g)\in C_{\Delta}^{+}([0,+\infty))$, if $f(x)=0$, then $x=f(x)g(x)=0$; hence \(0\) is the only possible zero of \(f\). This fact will be used throughout the paper without further mention. The multivariable counterparts of (i) and (iii) can be found in \cite[Theorem 4.12 and Lemma 4.13]{Zhou-glma-2023}. For completeness and the reader's convenience, we include a brief proof for the finite-dimensional, single-operator setting. 

\begin{proposition}\label{basicprop1}
	Let $(f,g)\in C_{\Delta}^{+}([0,+\infty))$.
	Then, for every $T\in\mathbb{M}_n(\mathbb{C})$,
	the following statements hold:
	\begin{enumerate}[label=(\roman*)]
		\item $\sigma(\Delta_{f,g}(T))=\sigma(T)$;
		
		\item $\operatorname{tr}((\Delta_{f,g}(T))^m)
		=\operatorname{tr}(T^m)$ for every $m\in\mathbb{N}$;
		
		\item $\Delta_{f,g}(UTU^*)
		=U\Delta_{f,g}(T)U^*$
		for every unitary $U\in\mathbb{M}_n(\mathbb{C})$;

	    \item If $f(0)>0$ or $T$ is invertible, then $\Delta_{f,g}(T)
	    =f(|T|)T(f(|T|))^{-1}$;
	   
		\item $\Delta_{f,g}(\omega T)
		=\omega\Delta_{f,g}(T)$
		for every $\omega\in\mathbb{C}$ with $|\omega|=1$;
		
		\item The map $r\mapsto\Delta_{f,g}(rT)$
		is continuous in norm on $\mathbb{R}$;
		
	\item If $g$ is constant, then
	$\Delta_{f,g}(T)=\Delta_D(T)$. Moreover, $\|\Delta_D(T)\|_F\leq\|T\|_F$ and
	$\Delta_D(\lambda T)=\lambda\Delta_D(T)$
	for every $\lambda\in\mathbb{C}$.
	\end{enumerate}
\end{proposition}
\begin{proof}
	Let $T=V|T|$ be the polar decomposition of $T$.
	
(i)	Since $f(x)g(x)=x$, the continuous functional calculus immediately gives $f(|T|)g(|T|)=|T|$. Also, it is well known that $\sigma(AB)=\sigma(BA)$ for all matrices $A,B\in\mathbb{M}_n(\mathbb{C})$. Consequently, $\sigma(T)=\sigma(V|T|)=\sigma(Vg(|T|)f(|T|))=\sigma(f(|T|)Vg(|T|))=\sigma(\Delta_{f,g}(T))$.

(ii) Recall that the trace is cyclic, that is,
\begin{equation}\label{tracecyc}
	\operatorname{tr}(AB)=\operatorname{tr}(BA)
\end{equation}
for all matrices $A,B\in\mathbb{M}_n(\mathbb{C})$. For $m\in\mathbb{N}$, a direct computation gives
$$
\begin{aligned}
	\operatorname{tr}\big((\Delta_{f,g}(T))^m\big)
	&=\operatorname{tr}\big((f(|T|)Vg(|T|))^m\big)\\
	&=\operatorname{tr}\Big(f(|T|)\big(Vg(|T|)f(|T|)\big)^{m-1}Vg(|T|)\Big)\\
	&=\operatorname{tr}\big(f(|T|)(V|T|)^{m-1}Vg(|T|)\big)\\
	&=\operatorname{tr}\big(f(|T|)T^{m-1}Vg(|T|)\big)\\
	&=\operatorname{tr}\big(T^{m-1}Vg(|T|)f(|T|)\big)\\
	&=\operatorname{tr}(T^{m-1}T)=\operatorname{tr}(T^m),
\end{aligned}$$
where the fifth equality follows from \eqref{tracecyc}.

(iii) For unitary matrix $U\in M_n(\mathbb{C})$, it is easy to verify that $UTU^\ast=UVU^\ast|UTU^\ast|$
is the polar decomposition of $UTU^\ast$, where $|UTU^\ast|=U|T|U^\ast$. Since $f$ is continuous on the compact set $\sigma(|T|)$, there exists a sequence of polynomials $\{p_n\}_{n\in\mathbb{N}}$ converging uniformly to $f$ on $\sigma(|T|)$. By the continuous functional calculus, $p_n(U|T|U^*)\longrightarrow f(U|T|U^*)$ in operator norm. Moreover, $p_n(U|T|U^\ast)=Up_n(|T|)U^\ast$, and hence $f(U|T|U^\ast)=Uf(|T|)U^\ast$. Similarly, $g(U|T|U^\ast)=Ug(|T|)U^\ast$. Therefore, $$\begin{aligned}
	\Delta_{f,g}(UTU^\ast)
	&=f(U|T|U^\ast)(UVU^\ast)g(U|T|U^\ast)\\&=Uf(|T|)U^\ast UVU^\ast Ug(|T|)U^\ast\\
	&=Uf(|T|)Vg(|T|)U^\ast\\
	&=U\Delta_{f,g}(T)U^\ast.
\end{aligned}$$

(iv) 
Suppose that $f(0)>0$ or $T$ is invertible.
Since $f(x)>0$ for all $x>0$, we have $0\notin\sigma(f(|T|))
=\{f(x):x\in\sigma(|T|)\}$
in either case. Combining this with
$f(x)g(x)=x$, we obtain $g(|T|)=|T|(f(|T|))^{-1}$. Consequently, $$\Delta_{f,g}(T)=f(|T|)Vg(|T|)=f(|T|)V|T|(f(|T|))^{-1}=f(|T|)T(f(|T|))^{-1}.$$

(v) Let $\omega\in\mathbb{C}$ with $|\omega|=1$. Since $\omega T=(\omega V)|\omega T|=(\omega V)|T|$
is the polar decomposition of $\omega T$, we obtain $\Delta_{f,g}(\omega T)=f(|T|)(\omega V)g(|T|)=\omega f(|T|)Vg(|T|)=\omega\Delta_{f,g}(T).$

(vi) For $r\in\mathbb{R}\setminus\{0\}$, the
polar decomposition of $rT$ is $rT=(\frac{r}{|r|}V)|rT|=(\frac{r}{|r|}V)(|r||T|)$, which gives
\begin{equation}\label{continuous}
\Delta_{f,g}(rT)=\frac{r}{|r|}f(|r||T|)Vg(|r||T|).
\end{equation}
 On the other hand, since $|T|$ is normal, there exists a unitary matrix $U$ such that $
U|T|U^\ast=\operatorname{diag}(\lambda_1,\lambda_2,\dots,\lambda_n)$, where $\lambda_i\in[0,+\infty)$ for $i=1,2,\dots,n$. From the proof of (iii), we have
\begin{equation}\label{continuous1}
\begin{aligned}
	f(|r||T|)
	&=f(|r|U^\ast\operatorname{diag}(\lambda_1,\lambda_2,\dots,\lambda_n)U)\\
	&=U^\ast f(|r|\operatorname{diag}(\lambda_1,\lambda_2,\dots,\lambda_n))U\\
	&=U^\ast\operatorname{diag}(f(|r|\lambda_1),f(|r|\lambda_2),\dots,f(|r|\lambda_n))U.
\end{aligned}
\end{equation}
Therefore, since $f\in C^+([0,+\infty))$, the map $r\longmapsto f(|r||T|)$ is continuous in norm. Similarly, $r\longmapsto g(|r||T|)$ is continuous in norm. Moreover, since $r\mapsto \frac{r}{|r|}$ is continuous on $\mathbb{R}\setminus\{0\}$, and matrix multiplication is continuous in norm, it follows from \eqref{continuous} that $r\longmapsto \Delta_{f,g}(rT)$ is continuous on $\mathbb{R}\setminus\{0\}$. Finally, from \eqref{continuous1} it follows easily that $\lim\limits_{r\to0}f(|r||T|)=f(0)I$. Similarly, $\lim\limits_{r\to0}g(|r||T|)=g(0)I$. Since $f(0)g(0)=0$, we obtain $$\lim\limits_{r\to0}\Delta_{f,g}(rT)=\lim_{r\to0}\frac{r}{|r|}f(|r||T|)Vg(|r||T|)=f(0)g(0)V=0=\Delta_{f,g}(0)=\Delta_{f,g}(0\cdot T).$$ Thus, the map \(r\mapsto\Delta_{f,g}(rT)\) is continuous at \(r=0\). The proof is complete.

(vii) By assumption, there exists \(c>0\) such that \(g(x)=c\) for all \(x\ge0\). Note that if \(c=0\), then $x=f(x)g(x)=0$, which contradicts the arbitrariness of $x$. Since $f(x)g(x)=x$, we have $f(x)=\frac{x}{c}$. Hence a direct computation gives $$\Delta_{f,g}(T)=f(|T|)Vg(|T|)=\frac{1}{c}|T|V(cI)=|T|V=\Delta_D(T).$$
Furthermore, since $V$ is a partial isometry, we have $\|V\|\leq 1$. Moreover, $\|AB\|_F\le\|A\|_F\|B\|$ for all matrices $A,B\in \mathbb{M}_{n}(\mathbb{C})$ \cite[(IV.40) and Proposition IV.2.4]{Bhatia-gtm-1997}. Hence
$$\|\Delta_{D}(T)\|_{F}=\||T|V\|_{F}\le\||T|\|_{F}\|V\|\le\|T\|_{F}.$$	
Finally, let $\lambda\in\mathbb{C}$. If $\lambda=0$, then $\Delta_D(\lambda T)=0=\lambda\Delta_D(T)$. If $\lambda\neq 0$, then $
\lambda T=\big(\frac{\lambda}{|\lambda|}V\big)(|\lambda||T|)$ is the polar decomposition of $\lambda T$.
Therefore, $\Delta_D(\lambda T)=|\lambda||T|(\frac{\lambda}{|\lambda|}V)=\lambda|T|V=\lambda\Delta_D(T)$. The proof is complete.
\end{proof}

The following results describe the $(f,g)$-Aluthge transforms of some special matrices. Recall that $T\in\mathbb{M}_n(\mathbb{C})$ is normal if and only if $\Delta_D(T)=T$. (i) and (ii) below provide
characterizations of normal matrices. In particular, (i) is closely related
to \cite[Theorem 3.3]{Zhou-glma-2023}. On the other hand, rank-one matrices also play an important role in our work. Let $\operatorname{rank}(T)$ denote the rank of $T\in\mathbb{M}_n(\mathbb{C})$. For nonzero $u,v\in\mathbb{C}^n$, the map $u\otimes v$ on $\mathbb{C}^n$ is defined by $x\mapsto\langle x,v\rangle u$, and it satisfies
\begin{equation}\label{rankonepro}
	\|u\otimes v\|_F=\|u\otimes v\|=\|u\|\|v\|,\quad\sigma(u\otimes v)=\{0,\langle u,v\rangle\}.
\end{equation}
We write a column vector $x\in\mathbb{C}^n$ as $x=(x_i)=\begin{pmatrix}
	x_1\\
	x_2\\
	\vdots\\
	x_n
\end{pmatrix}$,
where $x_i\in\mathbb{C}$ for $i=1,\ldots,n$. In particular, $e_j$ ($j\in\{1,2,\cdots,n\}$) denotes the $j$-th standard
basis vector in $\mathbb{C}^n$, whose $j$-th entry
is $1$ and all other entries are $0$. It is well known that $\operatorname{rank}(T)=1$ if and only if $T=u\otimes v$ for some nonzero vectors $u=(u_i)$ and $v=(v_i)$ in $\mathbb{C}^n$. In this case, $T=(t_{ij})=(u_i\overline{v_j})$, where $t_{ij}$ denotes the entry of $T$ in the \(i\)-th row and \(j\)-th column. These properties of $u\otimes v$ are given in \cite[p.~55, 62]{Murphy-AP-1990}. For $S,T\in\mathbb{M}_n(\mathbb{C})$, we say that $S\ge T$ if $\langle (S-T)x,x\rangle\ge0$ for all $x\in\mathbb{C}^n$. The dimension of a linear space $V$ is denoted by $\dim V$. 

\begin{proposition}\label{basicprop2}
	Let $(f,g)\in C_{\Delta}^{+}([0,+\infty))$.
	Then, for $T\in\mathbb{M}_n(\mathbb{C})$,
	the following statements hold:
	\begin{enumerate}[label=(\roman*)]
		\item If $f$ is nonconstant, then
		$\Delta_{f,g}(rT)=rT$ for every $r>0$
		if and only if $T$ is normal.
		If $f$ is constant, then $\Delta_{f,g}(T)=T$;
		
		\item Suppose that $g$ is constant.
		In this case, $\Delta_{f,g}(T)=\Delta_D(T)$. If $\|\Delta_D(T)\|_F=\|T\|_F$
		and there exists a dense subset
		$E\subseteq\mathbb{R}$ such that
		$$\Delta_D(T+r\Delta_D(T))
		=\Delta_D(T)+r\Delta_D(\Delta_D(T))$$
		for every $r\in E$, then $\Delta_D(T)=T$;
		
		\item If $T=u\otimes v$ for some nonzero
		vectors $u,v\in\mathbb{C}^n$, then
		$$\Delta_{f,g}(T)=\frac{f(\|u\|\|v\|)-f(0)}{f(\|u\|\|v\|)}\frac{\langle u,v\rangle}{\|v\|^2}v\otimes v+\frac{f(0)}{f(\|u\|\|v\|)}T;$$
		
		\item If $T$ is a partial isometry, then
		for every $r>0$, $$
		\Delta_{f,g}(rT)
		=f(0)g(r)T
		+(f(r)-f(0))g(r)T^*T^2;$$
		
		\item If $f(0)=0$, then
		$\Delta_{f,g}(T)=0$ if and only if $T^2=0$;
		
		\item Suppose that $f(0)=0$ and
		$g$ is nonconstant. If $\Delta_{f,g}(rT)=r\Delta_{f,g}(T)$
		for every $r>0$, then there exists
		$\lambda\in(0,1)$ such that 
		\begin{equation}\label{aluthgechange}
		\Delta_{f,g}(T)=\Delta_\lambda(T).
		\end{equation} Moreover, $\|\Delta_{f,g}(T)\|_F\leq\|T\|_F$,
		and equality holds if and only if
		$\Delta_{f,g}(T)=T$.
	\end{enumerate}
\end{proposition}
\begin{proof}
	 
(i) 
Let $T=U|T|$ be the polar decomposition of $T$.
If $f$ is constant, then there exists $c>0$
such that $f(x)=c$ for all $x\geq0$. Since $f(x)g(x)=x$, we have $g(x)=\frac{x}{c}$ for all $x\geq0$. It follows immediately that $\Delta_{f,g}(T)=f(|T|)Ug(|T|)=cU\frac{|T|}{c}=T$.

Suppose that $f$ is nonconstant.
If $T$ is normal, then $U|T|=|T|U$, which, by the continuous functional calculus,
implies that $Uf(|T|)=f(|T|)U$. Hence, $$\Delta_{f,g}(T)=f(|T|)Ug(|T|)=Uf(|T|)g(|T|)=T.$$ 

Conversely, suppose that
$\Delta_{f,g}(rT)=rT$ for every $r>0$. By Propositon \ref{basicprop1} (iii) and the spectral theorem,
we may assume that $|T|=\operatorname{diag}
(\lambda_1,\lambda_2,\ldots,\lambda_n)$, where $\lambda_{i}\ge0,i\in\{1,2,\cdots,n\}$. Write $U=(u_{ij})$. It suffices to show that
$\lambda_i=\lambda_j$ whenever $u_{ij}\neq0$
for all $1\leq i,j\leq n$.
Indeed, this implies that
$(\lambda_i-\lambda_j)u_{ij}=0$,
that is, $|T|U=U|T|$, and hence $T$ is normal.
Moreover, if $\lambda_j=0$, then
$e_j\in\ker|T|=\ker U$, and hence
$u_{ij}=0$ for every $i$.
Thus, we only need to consider the case
where $\lambda_j>0$ and $u_{ij}\neq0$. Now, since $rT=U(r|T|)$ is the polar decomposition
of $rT$, the assumption
$\Delta_{f,g}(rT)=rT$
implies that
\begin{equation}\label{normal-equ1}
	f(r\lambda_i)u_{ij}g(r\lambda_j)
	=r\lambda_j u_{ij}
\end{equation}
for all $1\leq i,j\leq n$ and $r>0$.
On the other hand, since $f(x)g(x)=x$
for every $x\geq0$, we have
\[
f(r\lambda_j)g(r\lambda_j)
=r\lambda_j>0.
\]
Combining this with \eqref{normal-equ1},
we obtain
\begin{equation}\label{normal-equ2}
	f(r\lambda_i)=f(r\lambda_j)
\end{equation}
for every $r>0$.
If $\lambda_i=0$, then \eqref{normal-equ2}
implies that
$$f(0)=f(r\lambda_j)$$
for every $r>0$.
Since $\lambda_j>0$, this implies that
$f$ is constant on $[0,+\infty)$,
a contradiction. Hence $\lambda_i>0$.
Suppose that $\lambda_i\neq\lambda_j$.
Set $a=\lambda_i/\lambda_j$.
Then $a>0$ and $a\neq1$.
Replacing $r$ by $x/\lambda_j$
in \eqref{normal-equ2}, we obtain
$f(ax)=f(x)$ for every $x>0$.
By mathematical induction, it follows that
$f(a^k x)=f(x)$ for every $k\in\mathbb Z$
and $x>0$.
Since $a\neq1$, there exists a sequence
of integers $(k_m)$ such that
$a^{k_m}\to0$.
By the continuity of $f$ at zero,
we obtain
$$f(x)=\lim_{m\to\infty}f(a^{k_m}x)=f(0)$$
for every $x>0$. This contradicts the fact that $f$ is nonconstant. Therefore, $\lambda_i=\lambda_j$, as desired.

(ii) By Proposition \ref{basicprop1} (vii), we have $\Delta_{f,g}(T)=\Delta_{D}(T)$. Let $T=U|T|$ be the polar decomposition of $T$. Since $I-UU^\ast$ is an orthogonal projection,
	we have
	\begin{equation}\label{pro}
		(I-UU^\ast)^2=I-UU^\ast.
	\end{equation}
Moreover,  a direct computation from $\|\Delta_{D}(T)\|_{F}=\|T\|_{F}$
gives $\operatorname{tr}(U^\ast|T|^2U)=\operatorname{tr}(T^\ast T)$.
	By the cyclic property of the trace, we obtain
	$\operatorname{tr}(|T|^2UU^\ast)
	=\operatorname{tr}(T^\ast T)$.
	Hence, by \eqref{pro},
	$$0=\operatorname{tr}(|T|^2(I-UU^\ast))=\operatorname{tr}(|T|(I-UU^\ast)|T|)=\|(I-UU^\ast)|T|\|_{F}^2.$$
	Consequently, $(I-UU^\ast)|T|=0$, which yields
	$$
	\operatorname{ran}|T|
	\subseteq\operatorname{ran}UU^\ast
	=\operatorname{ran}U.
	$$
	On the other hand, since
	$\ker U=\ker|T|$, the rank-nullity theorem
	gives
	$\dim\operatorname{ran}U
	=\dim\operatorname{ran}|T|$.
	Combining this with the preceding inclusion,
	we conclude that
	$\operatorname{ran}|T|=\operatorname{ran}U$. Therefore, with respect to the orthogonal
	decomposition
	$\mathbb{C}^n
	=\operatorname{ran}|T|\oplus\ker|T|$,
	we can write $U=
	\begin{pmatrix}
		\widetilde{U}&0\\
		0&0
	\end{pmatrix}$ and $|T|=
	\begin{pmatrix}
		\widetilde{|T|}&0\\
		0&0
	\end{pmatrix}$,
	where $\widetilde{U}$ and
	$\widetilde{|T|}$ are a unitary matrix and an invertible positive
	matrix acting on $\operatorname{ran}|T|$, respectively.
	Thus, without loss of generality, we may
	assume that $T$ is invertible; otherwise,
	we may restrict our argument to
	$\widetilde{T}
	=\widetilde{U}\widetilde{|T|}$. 
	
	Now $U$ is unitary. Set
	$A=U^\ast|T|U$.
	It is clear that it suffices to prove
	$A=|T|$ in order to obtain
	$\Delta_{D}(T)=T$.
	Indeed, since $|T|U=UA$ is a polar
	decomposition, we have
	$\Delta_{D}(\Delta_{D}(T))
	=\Delta_{D}(|T|U)=AU$.
	It follows from the hypothesis that
	\begin{equation}\label{duggal}
		\Delta_{D}(U(|T|+tA))
		=(|T|+tA)U
	\end{equation}
	for every $t\in E$. Furthermore, since $|T|$ is invertible,
	$\det(|T|+tA)$ is a polynomial in $t$
	that is not identically zero.
	Thus, there are only finitely many real
	values of $t$ for which $|T|+tA$
	is singular. Choose any $t\in E\cap(-1,0)$ such that
	$|T|+tA$ is invertible, and let
	\begin{equation}\label{polar1}
		|T|+tA=V\bigl||T|+tA\bigr|
	\end{equation}
	be its polar decomposition.
	Since $|T|+tA$ is self-adjoint, we have
	\begin{equation}\label{polar2}
		V^\ast=V,\qquad V^2=I.
	\end{equation}
Moreover, $U(|T|+tA)
	=UV\bigl||T|+tA\bigr|$
	is a polar decomposition. Combining this with \eqref{duggal}, \eqref{polar1},
	\eqref{polar2}, and the self-adjointness of
	$|T|+tA$, we obtain
		$$\bigl||T|+tA\bigr|UV=(|T|+tA)U=\bigl||T|+tA\bigr|VU.$$
	Since $\bigl||T|+tA\bigr|$ is invertible,
	it follows that $UV=VU$.
	Consequently,
	\begin{equation}\label{polar3}
		\operatorname{tr}(V|T|)
		=\operatorname{tr}(UV|T|U^\ast)
		=\operatorname{tr}(VA).
	\end{equation}
Since $|T|+tA$ is self-adjoint, we have
	$\bigl||T|+tA\bigr|\geq |T|+tA$.
	Therefore,
	$$
	\begin{aligned}
		0
		\leq
		\operatorname{tr}
		\left(\bigl||T|+tA\bigr|-(|T|+tA)\right)&=
		\operatorname{tr}
		\left((V-I)(|T|+tA)\right)\\
		&=
		\operatorname{tr}
		\left((V-I)(A+tA)\right)\\
		&=
		(1+t)\operatorname{tr}((V-I)A)\\
		&=
		-(1+t)\operatorname{tr}
		\left((I-V)^{\frac12}A(I-V)^{\frac12}\right)\le 0.
	\end{aligned}
	$$
	Here, the first equality follows from
	\eqref{polar1} and \eqref{polar2},
	the second equality follows from
	\eqref{polar3} together with
	$\operatorname{tr}(|T|)
	=\operatorname{tr}(A)$,
	and the last equality follows from
	the cyclic property of the trace. Thus, $\operatorname{tr}
	\left(\bigl||T|+tA\bigr|-(|T|+tA)\right)
	=0$. Since
	$\bigl||T|+tA\bigr|-(|T|+tA)\geq0$,
	we conclude that
	$$
	|T|+tA=\bigl||T|+tA\bigr|\geq0.$$
Letting $t\to -1^+$, we obtain $|T|-A\geq0$. On the other hand, $$\operatorname{tr}(|T|-A)=\operatorname{tr}(|T|)-\operatorname{tr}(U^\ast|T|U)=0.$$ Therefore, $|T|=A$, and hence $\Delta_{D}(T)=T$.

(iii) As shown in the proof of
\cite[Lemma 2.13]{Chabbabi-pams-2019},
the polar decomposition of $T=u\otimes v$ is
$T=U|T|$, where
\begin{equation}\label{rankonepolar1}
	U=\frac{1}{\|u\|\|v\|}u\otimes v,
	\quad
	|T|=\frac{\|u\|}{\|v\|}v\otimes v.
\end{equation}
Moreover, by \eqref{rankonepro}, we have
$\sigma(|T|)=\{0,\|u\|\|v\|\}$.
It follows that $f(|T|)=p(|T|)$ and
$g(|T|)=q(|T|)$, where
$p(x)=\frac{f(\|u\|\|v\|)-f(0)}
{\|u\|\|v\|}x+f(0)$
and $q(x)=\frac{g(\|u\|\|v\|)-g(0)}
{\|u\|\|v\|}x+g(0)$. Consequently, a direct computation using
\eqref{rankonepolar1} yields
\begin{align*}
	\Delta_{f,g}(T)
	&=f(|T|)Ug(|T|)\\
	&=p(|T|)\Big(\frac{g(\|u\|\|v\|)-g(0)}
	{\|u\|\|v\|}U|T|+g(0)U\Big)\\
	&=\Big(\frac{f(\|u\|\|v\|)-f(0)}
	{\|u\|\|v\|}|T|+f(0)\Big)g(\|u\|\|v\|)U\\
	&=\Big(\frac{f(\|u\|\|v\|)-f(0)}
	{\|u\|\|v\|}|T|+f(0)\Big)\frac{\|u\|\|v\|}{f(\|u\|\|v\|)}U\\
	&=\frac{f(\|u\|\|v\|)-f(0)}
	{f(\|u\|\|v\|)}
	\frac{\langle u,v\rangle}{\|v\|^2}
	v\otimes v
	+\frac{f(0)}{f(\|u\|\|v\|)}T.
\end{align*}

(iv) Let $T=U|T|$ be the polar decomposition of $T$.
Since $T$ is a partial isometry, it is well known
that $U=T$ and $|T|=T^*T$. Thus,
$rT=T(rT^*T)$ is the polar decomposition of $rT$
for every $r>0$.
Moreover, since $|T|=T^*T$ is an orthogonal
projection, a direct computation yields
\[
f(|rT|)=f(r|T|)
=f(0)I+(f(r)-f(0))T^*T
\]
and
\[
g(|rT|)=g(r|T|)
=g(0)I+(g(r)-g(0))T^*T.
\]
Consequently, using $T(T^*T)=T$ and
$(T^*T)^2=T^*T$, we obtain
\begin{align*}
	\Delta_{f,g}(rT)
	&=f(|rT|)Tg(|rT|)\\
	&=\bigl(f(0)T+(f(r)-f(0))T^*T^2\bigr)
	\bigl(g(0)I+(g(r)-g(0))T^*T\bigr)\\
	&=g(0)\bigl(f(0)T+(f(r)-f(0))T^*T^2\bigr)+(g(r)-g(0))
	\bigl(f(0)T+(f(r)-f(0))T^*T^2\bigr)\\
	&=g(r)\bigl(f(0)T+(f(r)-f(0))T^*T^2\bigr)\\
	&=f(0)g(r)T
	+(f(r)-f(0))g(r)T^*T^2,
\end{align*}
as desired.

(v) follows from \cite[Theorem 3.5]{Zhou-glma-2023}.

(vi) 
	The assertion is trivial if $T=0$.
	Hence, we assume that $T\neq0$.
	Let $T=U|T|$ be the polar decomposition of $T$,
	and write $U=(u_{ij})$.
	By Propositon \ref{basicprop1} (iii) and the spectral theorem,
	we may assume that $|T|=\operatorname{diag}
	(\lambda_1,\lambda_2,\ldots,\lambda_n)$, where $\lambda_{i}\ge0,i\in\{1,2,\cdots,n\}$. Since $rT=U(r|T|)$ is the polar decomposition
	of $rT$, the assumption
	$\Delta_{f,g}(rT)=r\Delta_{f,g}(T)$
	implies that
	\begin{equation}\label{basequ1}
		f(r\lambda_i)u_{ij}g(r\lambda_j)
		=rf(\lambda_i)u_{ij}g(\lambda_j)
	\end{equation}
	for all $1\leq i,j\leq n$. On the other hand, it is easy to see that
	\eqref{aluthgechange} holds if and only if
	there exists $\lambda\in(0,1)$ such that
	\begin{equation}\label{basequ}
		f(\lambda_i)u_{ij}g(\lambda_j)
		=\lambda_i^\lambda u_{ij}
		\lambda_j^{1-\lambda}
	\end{equation}
	for all $1\leq i,j\leq n$. Note that if $\lambda_j=0$, then
	$e_j\in\ker|T|$. Since
	$\ker U=\ker|T|$, we have $Ue_j=0$,
	which implies $u_{ij}=0$.
	Thus, \eqref{basequ} holds.
	If $\lambda_i=0$, then \eqref{basequ}
	also holds because $f(0)=0$. Therefore, if $\lambda_i=\lambda_j$
	for every pair $(i,j)$ satisfying
	$\lambda_i,\lambda_j>0$ and $u_{ij}\neq0$,
	then \eqref{basequ} follows immediately
	from $f(x)g(x)=x$ by taking
	$\lambda=\frac12$.
	
	It remains to consider the case where
	there exist $p,q\in\{1,2,\ldots,n\}$
	such that $\lambda_p,\lambda_q>0$,
	$u_{pq}\neq0$, and $\lambda_p\neq\lambda_q$.
	We first show that
	$$
	\lambda\triangleq
	\frac{\ln f(\lambda_p)-\ln f(\lambda_q)}
	{\ln\lambda_p-\ln\lambda_q}
	\in(0,1),$$
	and then prove that
	$\Delta_{f,g}(T)=\Delta_\lambda(T)$. Indeed, by \eqref{basequ1} and
	$f(x)g(x)=x$, we obtain
	\[
	\frac{f(r\lambda_p)}{f(r\lambda_q)}
	=\frac{f(\lambda_p)}{f(\lambda_q)}
	\]
	for every $r>0$.
	Since $r>0$ is arbitrary, replacing $r$
	by $x/\lambda_q$, where $x>0$,
	and then applying mathematical induction, we obtain
	\begin{equation}\label{prabound}
		f\left(
		\left(\frac{\lambda_p}{\lambda_q}\right)^k x
		\right)
		=
		\left(
		\frac{f(\lambda_p)}{f(\lambda_q)}
		\right)^k f(x)
	\end{equation}
	for every $k\in\mathbb Z$.
	By the definition of $\lambda$, this becomes
	\begin{equation}\label{geze}
		f\left(
		\left(\frac{\lambda_p}{\lambda_q}\right)^k x
		\right)
		=
		\left(\frac{\lambda_p}{\lambda_q}\right)^{k\lambda}
		f(x).
	\end{equation}
	Furthermore, using $f(x)g(x)=x$, we obtain
	\begin{equation}\label{leze}
		g\left(
		\left(\frac{\lambda_p}{\lambda_q}\right)^k x
		\right)
		=
		\left(\frac{\lambda_p}{\lambda_q}\right)^{k(1-\lambda)}
		g(x).
	\end{equation}	Since $\lambda_p\neq\lambda_q$,
	there exists a sequence of integers $(k_m)$
	with $|k_m|\to\infty$ such that
	\[
	\left(\frac{\lambda_p}{\lambda_q}\right)^{k_m}
	\longrightarrow0.
	\]
	By \eqref{geze}, the continuity of $f$
	at zero, and $f(0)=0$, we obtain
	$\lambda>0$.
	Similarly, \eqref{leze} and the continuity
	of $g$ at zero imply that $\lambda\leq1$.
	If $\lambda=1$, then \eqref{leze} yields $g\left(
	\left(\frac{\lambda_p}{\lambda_q}\right)^{k_m}x
	\right)=g(x)$
	for every $x>0$.
	Letting $m\to\infty$, we obtain
	$g(x)=g(0)$ for every $x>0$,
	which contradicts the nonconstancy of $g$.
	Therefore, $\lambda\in(0,1)$.
	
	Next, consider the function $P(x)=\ln f(e^x)-\lambda x$, $x\in\mathbb R$. By \eqref{geze}, we have
	\begin{align*}
		P\left(x+\ln\frac{\lambda_p}{\lambda_q}\right)
		&=\ln f\left(
		e^x\frac{\lambda_p}{\lambda_q}
		\right)
		-\lambda\left(
		x+\ln\frac{\lambda_p}{\lambda_q}
		\right)\\&=\ln\left(
		\left(\frac{\lambda_p}{\lambda_q}\right)^\lambda
		f(e^x)\right)
		-\lambda x
		-\lambda\ln\frac{\lambda_p}{\lambda_q}\\&=\ln f(e^x)-\lambda x=P(x).
	\end{align*}
	Hence $P$ is a continuous periodic function
	on $\mathbb R$, and thus is bounded.
	Consequently,
	\begin{equation}\label{lim}
\lim_{x\to0^+}\frac{\ln f(x)}{\ln x}=\lim_{x\to0^+}\left(\lambda+\frac{P(\ln x)}{\ln x}\right)=\lambda.
	\end{equation}
Now let $1\leq i,j\leq n$ be such that
	$\lambda_i,\lambda_j>0$ and $u_{ij}\neq0$.
	If $\lambda_i=\lambda_j$, then clearly
	\[
	\frac{f(\lambda_i)}{f(\lambda_j)}
	=\left(\frac{\lambda_i}{\lambda_j}\right)^\lambda.
	\]
	Suppose that $\lambda_i\neq\lambda_j$.
	By the same argument as that used to derive
	\eqref{prabound}, with $x=1$, we obtain
	\[
	f\left(
	\left(\frac{\lambda_i}{\lambda_j}\right)^k
	\right)
	=
	\left(
	\frac{f(\lambda_i)}{f(\lambda_j)}
	\right)^k f(1)
	\]
	for every $k\in\mathbb Z$. Choose a sequence of integers $(\ell_m)$
	such that $\left(\frac{\lambda_i}{\lambda_j}\right)^{\ell_m}
	\longrightarrow0$. It follows from \eqref{lim} that
	\begin{align*}
		\lambda
		&=\lim_{m\to\infty}
		\frac{
			\ln f\left(
			\left(\frac{\lambda_i}{\lambda_j}\right)^{\ell_m}
			\right)}
		{\ln\left(
			\left(\frac{\lambda_i}{\lambda_j}\right)^{\ell_m}
			\right)}=\lim_{m\to\infty}
		\frac{
			\ell_m\ln\frac{f(\lambda_i)}{f(\lambda_j)}
			+\ln f(1)}
		{\ell_m\ln\frac{\lambda_i}{\lambda_j}}=\frac{\ln f(\lambda_i)-\ln f(\lambda_j)}
		{\ln\lambda_i-\ln\lambda_j}.
	\end{align*}
	Therefore, $\frac{f(\lambda_i)}{f(\lambda_j)}
	=\left(\frac{\lambda_i}{\lambda_j}\right)^\lambda$. Consequently, for every pair $(i,j)$
	satisfying $\lambda_i,\lambda_j>0$
	and $u_{ij}\neq0$, we have
	\begin{align*}
		f(\lambda_i)u_{ij}g(\lambda_j)
		&=\lambda_j
		\frac{f(\lambda_i)}{f(\lambda_j)}u_{ij}=\lambda_j
		\left(\frac{\lambda_i}{\lambda_j}\right)^\lambda
		u_{ij}=\lambda_i^\lambda
		u_{ij}\lambda_j^{1-\lambda}.
	\end{align*}
	Thus, \eqref{basequ} holds.
	Combining this with the cases
	$\lambda_i=0$ and $\lambda_j=0$,
	we conclude that
	$$
	\Delta_{f,g}(T)=\Delta_\lambda(T),$$
	where
	$$
	\lambda=
	\begin{cases}
		\displaystyle
		\frac{\ln f(\lambda_p)-\ln f(\lambda_q)}
		{\ln\lambda_p-\ln\lambda_q},
		&
		\begin{gathered}
			\text{if there exist $p,q$ such that $\lambda_p,\lambda_q>0, u_{pq}\neq0$, and $\lambda_p\neq\lambda_q$}.
		\end{gathered}
		\\[12pt]
		\displaystyle\frac12,
		&\text{otherwise}.
	\end{cases}$$ 
	
	Finally, recall that in finite dimensions,
	$\Delta_\lambda(T)=T$ if and only if
	$T$ is normal.
	Therefore, the remaining assertions follow
	immediately from \cite[Theorem 3.9]{Antezana-laa-2005} and \eqref{aluthgechange}. 
	
	This completes the proof.
\end{proof}
\begin{remark}
	In general, $\|\Delta_D(T)\|_F=\|T\|_F$ does not imply that $\Delta_D(T)=T$. Indeed, consider $T=\begin{pmatrix}
		0 & 2\\
		1 & 0
	\end{pmatrix}$. Then $T=U|T|$ is the polar decomposition of $T$, where $U=\begin{pmatrix}
		0 & 1\\
		1 & 0
	\end{pmatrix}$ and $|T|=\begin{pmatrix}
		1 & 0\\
		0 & 2
	\end{pmatrix}$. Thus, $\Delta_D(T)=|T|U
	=\begin{pmatrix}
		0 & 1\\
		2 & 0
	\end{pmatrix}\neq T,
	$
	while $\|\Delta_D(T)\|_F=\|T\|_F=\sqrt{5}$.
\end{remark}

To establish the final result in this section, we need the following lemma, which is a strengthening of \v{S}emrl's result \cite[Lemma 2.4]{Semr-laa-l1991}.

\begin{lemma}\label{specequ}
	Let $T,S\in\mathbb{M}_n(\mathbb{C})$. If $\sigma(T+zP)=\sigma(S+zP)$ for every rank-one orthogonal projection $P\in\mathbb{M}_n(\mathbb{C})$ and every $z\in\mathbb{C}$, then $T=S$.
\end{lemma}

\begin{proof}
	Fix a unit vector $u\in\mathbb{C}^n$, and let
	$P_u=u\otimes u$ be the rank-one orthogonal
	projection onto the subspace spanned by $u$.
	Choose $\lambda\in\mathbb{C}$ such that $|\lambda|>\max\{\|T\|,\|S\|\}$. Then both $\lambda I-T$ and $\lambda I-S$ are invertible. For $X\in\{T,S\}$, set
	$R_X=(\lambda I-X)^{-1}$.
	Then, for every $z\in\mathbb{C}$,
	$$
	\lambda I-(X+zP_u)
	=(\lambda I-X)(I-zR_XP_u).$$
	Thus, $\lambda I-(X+zP_u)$ is invertible
	if and only if $I-zR_XP_u$ is invertible.
	By the hypothesis that
	$\sigma(T+zP_u)=\sigma(S+zP_u)$, we obtain
	\begin{equation}\label{rankoneequ}
		I-zR_TP_u\text{ is invertible if and only if }
		I-zR_SP_u\text{ is invertible}.
	\end{equation}
	Recall that, for any $x,y\in\mathbb{C}^n$,
	the operator $I-x\otimes y$ is invertible
	if and only if $1-\langle x,y\rangle\neq 0$.
	Since $zR_XP_u=zR_Xu\otimes u$,
	it follows from \eqref{rankoneequ} that
	$$
	1-z\langle R_Tu,u\rangle\neq 0
	\quad\text{if and only if}\quad
	1-z\langle R_Su,u\rangle\neq 0$$
	for every $z\in\mathbb{C}$.
	Since $z$ is arbitrary, we conclude that$
	\langle R_Tu,u\rangle
	=\langle R_Su,u\rangle$. Since $u$ is an arbitrary unit vector, we have $R_T=R_S$. Consequently,
	$$
	\lambda I-T=R_T^{-1}=R_S^{-1}
	=\lambda I-S,$$
	and hence $T=S$.
\end{proof}

 A characterization of the elements of $G_{f,g}$ is given in (i) below. Two rigidity results arising from identities
involving the $(f,g)$-Aluthge transform
are also established.

\begin{proposition}\label{basicprop3}
	Let $(f,g)\in C_{\Delta}^{+}([0,+\infty))$.
	Then the following statements hold:
	\begin{enumerate}[label=(\roman*)]
		
		\item For $r>0$, we have $r\in G_{f,g}$
		if and only if there exists a constant
		$R>0$ such that
		\[
		f(rx)=Rf(x)
		\]
		for every $x>0$.
	Moreover, $G_{f,g}$ has one of the
	following forms:
	\[
	G_{f,g}=\{1\},\quad
	G_{f,g}=q^{\mathbb{Z}}
	\text{ for some }q>1,\quad
	\text{or}\quad
	G_{f,g}=(0,\infty).
	\]
		
		\item Suppose that $f$ is nonconstant.
		If $T,S\in\mathbb{M}_n(\mathbb{C})$ satisfy
		\begin{equation}\label{defequ}
			\Delta_{f,g}(X+T)=\Delta_{f,g}(X)+S
		\end{equation}
		for every $X\in\mathbb{M}_n(\mathbb{C})$,
		then $T=S=0$.
		
		\item Suppose that $f$ is nonconstant
		and $\gamma:\mathbb{C}\to\mathbb{C}$
		is an additive map. If
	$$\Delta_{f,g}\left(
	\gamma(z)I+\begin{pmatrix}
		0&r\\
		0&0
	\end{pmatrix}\oplus 0
	\right)-\gamma(z)I=
	\Delta_{f,g}\left(
	zI+
	\begin{pmatrix}
		0&r\\
		0&0
	\end{pmatrix}\oplus 0
	\right)-zI$$
		for every $z\in\mathbb{C}$ and $r>0$,
		then $\gamma(z)=z$ for every
		$z\in\mathbb{C}$.
	\end{enumerate}
\end{proposition}
\begin{proof}
(i) Suppose that $r\in G_{f,g}$. Then
\begin{equation}\label{radequ}
	\Delta_{f,g}(rT)=r\Delta_{f,g}(T)
\end{equation}
for every $T\in\mathbb{M}_n(\mathbb{C})$. Taking $T=\begin{pmatrix}
	0&y\\
	x&0
\end{pmatrix}\oplus0$, where $x,y>0$, we have the polar decomposition
$T=U|T|$, where $
U=\begin{pmatrix}
	0&1\\
	1&0
\end{pmatrix}\oplus0$ and $|T|=\begin{pmatrix}
	x&0\\
	0&y
\end{pmatrix}\oplus0$. By a direct computation, \eqref{radequ} yields
$$
f(rx)g(ry)=rf(x)g(y)
$$
for all $x,y>0$.
Since $f(t)g(t)=t$ for every $t\geq0$,
we obtain
$$
\frac{f(rx)}{f(ry)}
=\frac{f(x)}{f(y)}.$$
Taking $y=1$, we have
$
f(rx)=\frac{f(r)}{f(1)}f(x)$
for every $x>0$. Thus, the desired identity holds with
$R=\frac{f(r)}{f(1)}>0$.

Conversely, suppose that there exists
a constant $R>0$ such that
\[
f(rx)=Rf(x)
\]
for every $x>0$.
By the continuity of $f$ at zero,
this identity also holds for $x=0$.
Moreover, since $f(t)g(t)=t$,
we have
$$
g(rx)=\frac{r}{R}g(x)$$
for every $x>0$.
By the continuity of $g$ at zero,
this identity also holds for $x=0$. Now let $T=U|T|$ be the polar decomposition
of an arbitrary matrix
$T\in\mathbb{M}_n(\mathbb{C})$.
By the continuous functional calculus,
we have $f(r|T|)=Rf(|T|)$ and $g(r|T|)=\frac{r}{R}g(|T|)$.
Consequently,
\begin{align*}
	\Delta_{f,g}(rT)=f(r|T|)Ug(r|T|)=Rf(|T|)U\frac{r}{R}g(|T|)=r\Delta_{f,g}(T).
\end{align*}
Hence $r\in G_{f,g}$, as desired.

We can now write $
G_{f,g}
=\left\{
r>0:
f(rx)=\frac{f(r)}{f(1)}f(x)
\text{ for all }x>0
\right\}$. By the continuity of $f$, it is clear that
$G_{f,g}$ is a closed subgroup. Since the logarithmic map is a topological
group isomorphism from $((0,\infty),\cdot)$
onto $(\mathbb{R},+)$,
$H:=\log G_{f,g}$ is a closed subgroup
of $(\mathbb{R},+)$.
By \cite[Proposition 3.1.11]{Aussenhofer-sim-2022},
every proper closed subgroup of
$(\mathbb{R},+)$ is cyclic. Thus
\[
H=\{0\},\quad
H=a\mathbb{Z}\text{ for some }a>0,
\quad\text{or}\quad H=\mathbb{R}.
\]
Applying the exponential map, with
$q=e^a>1$ in the second case, yields
\[
G_{f,g}=\{1\},\quad
G_{f,g}=q^{\mathbb{Z}}
\text{ for some }q>1,\quad
\text{or}\quad
G_{f,g}=(0,\infty).
\]

(ii) Let $\alpha\in\mathbb{C}$, and let $P$ be a rank-one orthogonal projection.
By Proposition \ref{basicprop2} (i), we have $\Delta_{f,g}(\alpha P)=\alpha P$.
Taking $X=\alpha P$ in \eqref{defequ}
and using Proposition~\ref{basicprop1}(i),
we obtain
$\sigma(\alpha P+T)=\sigma(\alpha P+S)$.
It follows from Lemma~\ref{specequ} that $T=S$.
Consequently, $T\in\mathbb{M}_n(\mathbb{C})$
satisfies
\begin{equation}\label{defequ1}
	\Delta_{f,g}(X+T)=\Delta_{f,g}(X)+T
\end{equation}
for every $X\in\mathbb{M}_n(\mathbb{C})$.
We now prove that $T=0$.
To this end, we first show that
$T=zI$ for some $z\in\mathbb{C}$,
and then prove that $z=0$.

Since $f$ is nonconstant, there exists $r>0$
such that $f(r)\neq f(0)$. Hence,
\begin{equation}\label{consne}
	\frac{f(0)}{f(r)}\neq1.
\end{equation}
For any unit vectors $u,v\in\mathbb{C}^n$
satisfying $\langle u,v\rangle=0$,
let $X_{uv}=ru\otimes v$. By Proposition \ref{basicprop2} (iii), we have $\Delta_{f,g}(X_{uv})
=\frac{f(0)}{f(r)}X_{uv}$. Thus, \eqref{defequ1} yields
\[
\Delta_{f,g}(X_{uv}+T)
=\frac{f(0)}{f(r)}X_{uv}+T.
\]
By Proposition \ref{basicprop1} (ii), we obtain
$$\operatorname{tr}((X_{uv}+T)^2)
=\operatorname{tr}\left(
\left(\frac{f(0)}{f(r)}X_{uv}+T\right)^2
\right).$$
Since $X_{uv}^2=0$, expanding both sides gives $\operatorname{tr}\left(
\left(1-\frac{f(0)}{f(r)}\right)
(X_{uv}T+TX_{uv})
\right)=0$.
By \eqref{consne} and the cyclicity of the trace,
we obtain
\[
\langle Tu,v\rangle
=\operatorname{tr}(T(u\otimes v))
=\frac{1}{r}\operatorname{tr}(TX_{uv})=0.
\]
Since $u$ and $v$ are arbitrary orthogonal vectors,
there exists $z\in\mathbb{C}$ such that $T=zI$. We claim that $z=0$. Otherwise, we shall
derive
$\frac{f(0)}{f(r)}=1$,
contrary to \eqref{consne}. Indeed, taking $X=\widetilde{X}\oplus0$ in \eqref{defequ1}, where $\widetilde{X}=
\begin{pmatrix}
	0&r\\
	0&0
\end{pmatrix}$, and using
$\Delta_{f,g}(A\oplus B)
=\Delta_{f,g}(A)\oplus\Delta_{f,g}(B)$
for any matrices $A$ and $B$, we obtain
\[
\Delta_{f,g}(\widetilde{X}+zI)
=\Delta_{f,g}(\widetilde{X})+zI,
\]
where $I$ denotes the $2\times2$ identity matrix. Since $\widetilde{X}
=r\begin{pmatrix}1\\0\end{pmatrix}
\otimes
\begin{pmatrix}0\\1\end{pmatrix}$, Proposition~\ref{basicprop1}(iv) and
Proposition~\ref{basicprop2}(iii) yield
$f(|\widetilde{X}+zI|)(\widetilde{X}+zI)
(f(|\widetilde{X}+zI|))^{-1}=\frac{f(0)}{f(r)}\widetilde{X}+zI$.
Consequently,
\begin{equation}\label{contra}
f(|\widetilde{X}+zI|)\widetilde{X}
=\frac{f(0)}{f(r)}
\widetilde{X}f(|\widetilde{X}+zI|).
\end{equation}
By the above identity and the invertibility of
$\widetilde{X}+zI$, there exist $a,b>0$ such that $$f(|\widetilde{X}+zI|)
=\begin{pmatrix}
	a&0\\
	0&b
\end{pmatrix}.$$ On the other hand, $f(|\widetilde{X}+zI|)$
commutes with
$
(\widetilde{X}+zI)^*(\widetilde{X}+zI)
=
\begin{pmatrix}
	|z|^2&\overline{z}r\\
	zr&|z|^2+r^2
\end{pmatrix}$. Therefore, $\begin{pmatrix}
	0&(a-b)\overline{z}r\\
	(b-a)zr&0
\end{pmatrix}
=0$.
Since $zr\neq0$, we have $a=b$. Therefore, by \eqref{contra}, we have $\frac{f(0)}{f(r)}=1$, as desired.

(iii) Since $f$ is nonconstant, there exists $r>0$
such that $f(r)\neq f(0)$.
Moreover, for any matrices $A$ and $B$, we have $\Delta_{f,g}(A\oplus B)=\Delta_{f,g}(A)\oplus\Delta_{f,g}(B)$. Thus, by the assumption and the additivity of
$\gamma$, we obtain
\begin{equation}\label{assuequ}
	\begin{aligned}
		&\Delta_{f,g}\left(
		\begin{pmatrix}
			q\gamma(z)&r\\
			0&q\gamma(z)
		\end{pmatrix}
		\right)
		-
		\begin{pmatrix}
			q\gamma(z)&0\\
			0&q\gamma(z)
		\end{pmatrix}\\
		&\qquad=
		\Delta_{f,g}\left(
		\begin{pmatrix}
			qz&r\\
			0&qz
		\end{pmatrix}
		\right)
		-
		\begin{pmatrix}
			qz&0\\
			0&qz
		\end{pmatrix}
	\end{aligned}
\end{equation}
for every nonzero $z\in\mathbb{C}$ and every
positive rational number $q$.

On the other hand, since
$\lim_{q\to0^+}
\left|
\begin{pmatrix}
	qz&r\\
	0&qz
\end{pmatrix}
\right|
=
\begin{pmatrix}
	0&0\\
	0&r
\end{pmatrix}$,
the continuous functional calculus yields
$$
\lim_{q\to0^+}
f\left(
\left|
\begin{pmatrix}
	qz&r\\
	0&qz
\end{pmatrix}
\right|
\right)
=
\begin{pmatrix}
	f(0)&0\\
	0&f(r)
\end{pmatrix}.$$
Since $f(r)\neq f(0)$, there exists a
positive rational number $q_0$ such that
\begin{equation}\label{neque}
	(f(|T|))_{11}\neq(f(|T|))_{22},
\end{equation}
where
$T=\begin{pmatrix}
	q_0z&r\\
	0&q_0z
\end{pmatrix}$
and $(f(|T|))_{ii}$ denotes the $(i,i)$-th
entry of $f(|T|)$ for $i=1,2$. 

Set $
S=
\begin{pmatrix}
	q_0\gamma(z)&r\\
	0&q_0\gamma(z)
\end{pmatrix}$. We claim that $f(|S|)$ is invertible.
Indeed, if $f(|S|)$ were singular, then
$S$ would be singular, which implies that
$\gamma(z)=0$.
Since $S$ is then a rank-one nilpotent
matrix with singular values $0$ and $r$,
the singularity of $f(|S|)$ also implies
that $f(0)=0$.
By Proposition \ref{basicprop2} (v), we obtain
$\Delta_{f,g}(S)=0$.
On the other hand, since $T$ is invertible, Proposition \ref{basicprop1} (iv) and
\eqref{assuequ} yield $0=f(|T|)
\begin{pmatrix}
	0&r\\
	0&0
\end{pmatrix}
(f(|T|))^{-1}$, which is impossible because $r>0$.
Therefore, $f(|S|)$ is invertible. By Proposition \ref{basicprop1} (iv) and
\eqref{assuequ}, we further obtain
$$
f(|S|)
\begin{pmatrix}
	0&1\\
	0&0
\end{pmatrix}
(f(|S|))^{-1}
=
f(|T|)
\begin{pmatrix}
	0&1\\
	0&0
\end{pmatrix}
(f(|T|))^{-1}.$$
Thus, $(f(|T|))^{-1}f(|S|)$ commutes with
$\begin{pmatrix}0&1\\0&0\end{pmatrix}$.
Consequently, there exist
$\alpha\in\mathbb{C}\setminus\{0\}$
and $\beta\in\mathbb{C}$ such that $
(f(|T|))^{-1}f(|S|)
=
\begin{pmatrix}
	\alpha&\beta\\
	0&\alpha
\end{pmatrix}$.
Moreover,
$$(f(|T|))^{-1}f(|S|)=(f(|T|))^{-\frac12}
\Big(
(f(|T|))^{-\frac12}
f(|S|)
(f(|T|))^{-\frac12}
\Big)\cdot(f(|T|))^{\frac12}.
$$
Hence $(f(|T|))^{-1}f(|S|)$ is similar
to the self-adjoint matrix
$
(f(|T|))^{-\frac12}
f(|S|)
(f(|T|))^{-\frac12},
$
and is therefore diagonalizable.
It follows that $\beta=0$, and hence $
f(|S|)=\alpha f(|T|)$. Consequently, a direct computation yields
$$	0=[f(|T|),T^*T-S^*S]=
\begin{pmatrix}
	\times&
	q_0r\bigl((f(|T|))_{11}-(f(|T|))_{22}\bigr)
	(\overline{z}-\overline{\gamma(z)})\\
	\times&\times
\end{pmatrix}.
$$
By \eqref{neque}, we obtain
$\gamma(z)=z$ for every nonzero
$z\in\mathbb{C}$.
Moreover, since $\gamma$ is additive,
we have $\gamma(0)=0$.
Therefore, $\gamma(z)=z$ for every
$z\in\mathbb{C}$, as desired.
\end{proof}

	\section{Additive Maps Commuting with the $(f,g)$-Aluthge Transform}

	In this section, we characterize surjective
	additive maps $\Phi:\mathbb{M}_{n}(\mathbb{C})\to\mathbb{M}_{n}(\mathbb{C})$ satisfying
	\begin{equation}\label{commualuthge}
		\Phi\circ\Delta_{f,g}
		=\Delta_{f,g}\circ\Phi.
	\end{equation}
	In the next section, we will show that
	every surjective map satisfying
	\eqref{eq:weighted} is necessarily additive
	and satisfies \eqref{commualuthge},
	thereby establishing our main result,
	Theorem \ref{thm:main}. Our approach is based on the following
	result of \v{S}emrl. Let $\mathfrak{sl}_n(\mathbb{C})$ denote the set
	of all $n\times n$ complex matrices with trace zero.
	
\begin{theorem}[{\cite[Corollary 2]{Semrl-laa-1993}}]\label{semrlprezero}
	Let $n\geq2$ and let
	$\Phi:\mathfrak{sl}_n(\mathbb{C})
	\to\mathfrak{sl}_n(\mathbb{C})$
	be an invertible complex-linear map
	such that $\Phi(A)^2=0$ whenever $A^2=0$.
	Then there exist a nonzero scalar
	$c\in\mathbb{C}$ and an invertible matrix
	$S\in\mathbb{M}_n(\mathbb{C})$
	such that either
	\[
	\Phi(A)=cSAS^{-1}
	\quad\text{or}\quad
	\Phi(A)=cSA^{\mathrm{tr}}S^{-1}
	\]
	for all $A\in\mathfrak{sl}_n(\mathbb{C})$.
\end{theorem}	
	In fact, under condition \eqref{commualuthge},
	we can initially show that the restriction $\Phi|_{\mathfrak{sl}_n(\mathbb{C})}$
	is invertible and real-linear, and that
	$\Phi(A)^2=0$ whenever $A^2=0$. To apply Theorem \ref{semrlprezero},
	we first establish Proposition \ref{realsemrl}
	using Theorem \ref{dimthm}. We then consider the maps
$\Phi(X)$ and $\Phi(\overline{X})$
on $\mathfrak{sl}_n(\mathbb{C})$.
By applying Theorem \ref{semrlprezero},
we obtain Theorem \ref{keyth}.
Notably, Theorem \ref{keyth} can be viewed
as a real-linear analogue of
Theorem \ref{semrlprezero}, as its proof shows.
This in turn leads to the main result
of this section, Theorem \ref{kkeylem} (vi).
	
	We start with the following lemma.
	
	\begin{lemma}\label{ralema}
		Let $\Phi:\mathfrak{sl}_n(\mathbb{C})
		\to\mathfrak{sl}_n(\mathbb{C})$
		be a real-linear map. Suppose that for any nonzero vectors
	$u,v\in\mathbb{C}^n$ with
	$\langle u,v\rangle=0$,
	there exists a scalar
	$\lambda_{u,v}\in\mathbb{C}$ such that $\Phi(u\otimes v)
	=\lambda_{u,v}(u\otimes v)$. Then there exists a scalar
		$\lambda\in\mathbb{C}$ such that $\Phi(T)=\lambda T$ for all $T\in\mathfrak{sl}_n(\mathbb{C})$.
	\end{lemma}
	\begin{proof}
		Let $e_1,e_2,\ldots,e_n$ be the standard
		orthonormal basis of $\mathbb{C}^n$.
		Fix $i\neq j$ and $z\in\mathbb{C}$.
		Since $(e_i-e_j)\perp(e_i+e_j)$ and $(e_i+e_j)\perp(e_i-e_j)$, the hypothesis implies that there exist
		$a_z^{ij},b_z^{ij}\in\mathbb{C}$ such that
		$$
		\Phi(z(e_i-e_j)\otimes(e_i+e_j))
		=a_z^{ij}(e_i-e_j)\otimes(e_i+e_j),\quad\Phi(z(e_i+e_j)\otimes(e_i-e_j))
		=b_z^{ij}(e_i+e_j)\otimes(e_i-e_j).
		$$
		Noting that
		\[
		\begin{aligned}
			&(e_i-e_j)\otimes(e_i+e_j)
			-(e_i+e_j)\otimes(e_i-e_j)\\
			&\qquad=2e_i\otimes e_j-2e_j\otimes e_i
		\end{aligned}
		\]
		and
		\[
		\begin{aligned}
			&(e_i-e_j)\otimes(e_i+e_j)
			+(e_i+e_j)\otimes(e_i-e_j)\\
			&\qquad=2e_i\otimes e_i-2e_j\otimes e_j,
		\end{aligned}
		\]
		by the real-linearity of $\Phi$,
		adding and subtracting the above identities
		yield
		\begin{equation}\label{rankequ1}
			2\Phi(ze_i\otimes e_j)
			-2\Phi(ze_j\otimes e_i)=(a_z^{ij}-b_z^{ij})
			(e_i\otimes e_i-e_j\otimes e_j)+(a_z^{ij}+b_z^{ij})
			(e_i\otimes e_j-e_j\otimes e_i)
		\end{equation}
		and
		\begin{equation}\label{rankequ2}
		2\Phi(z(e_i\otimes e_i-e_j\otimes e_j))=(a_z^{ij}+b_z^{ij})
		(e_i\otimes e_i-e_j\otimes e_j)\quad+(a_z^{ij}-b_z^{ij})
		(e_i\otimes e_j-e_j\otimes e_i).
		\end{equation}
		Moreover, by the hypothesis, there exist
		$c_z^{ij},d_z^{ij}\in\mathbb{C}$ such that
		\begin{equation}\label{scaequ1}
			\Phi(ze_i\otimes e_j)
			=c_z^{ij}(e_i\otimes e_j)
		\end{equation}
		and
		\begin{equation}\label{scaequ2}
			\Phi(ze_j\otimes e_i)
			=d_z^{ij}(e_j\otimes e_i).
		\end{equation}
		Since
		$e_i\otimes e_j$, $e_j\otimes e_i$,
		$e_i\otimes e_i$, and $e_j\otimes e_j$
		are linearly independent,
		comparing coefficients in
		\eqref{rankequ1} gives
		\begin{equation}\label{scaequ3}
			a_z^{ij}=b_z^{ij}=c_z^{ij}=d_z^{ij}.
		\end{equation}
		Consequently, \eqref{rankequ2} yields
		\begin{equation}\label{scaequ4}
			\Phi(z(e_i\otimes e_i-e_j\otimes e_j))
			=a_z^{ij}(e_i\otimes e_i-e_j\otimes e_j).
		\end{equation}Next, observe that
		$(e_i+\mathrm{i}e_j)\perp
		(e_i-\mathrm{i}e_j)$ and
		$$
(e_i+\mathrm{i}e_j)
\otimes(e_i-\mathrm{i}e_j)=e_i\otimes e_i-e_j\otimes e_j+\mathrm{i}(e_i\otimes e_j+e_j\otimes e_i).
		$$
		By \eqref{scaequ1}--\eqref{scaequ4},
		the real-linearity of $\Phi$, and the
		hypothesis, there exists
		$f_z^{ij}\in\mathbb{C}$ such that
		\[
		\begin{aligned}
			&a_z^{ij}(e_i\otimes e_i-e_j\otimes e_j)
			+a_{\mathrm{i}z}^{ij}
			(e_i\otimes e_j+e_j\otimes e_i)\\
			&=\Phi(z(e_i\otimes e_i-e_j\otimes e_j))+\Phi(\mathrm{i}z
			(e_i\otimes e_j+e_j\otimes e_i))\\
			&=\Phi(z(e_i+\mathrm{i}e_j)
			\otimes(e_i-\mathrm{i}e_j))\\
			&=f_z^{ij}(e_i\otimes e_i-e_j\otimes e_j)+\mathrm{i}f_z^{ij}
			(e_i\otimes e_j+e_j\otimes e_i).
		\end{aligned}
		\]
		By the linear independence of
		$e_i\otimes e_j$, $e_j\otimes e_i$,
		$e_i\otimes e_i$, and $e_j\otimes e_j$,
		we obtain $a_z^{ij}=f_z^{ij}$ and $a_{\mathrm{i}z}^{ij}=\mathrm{i}f_z^{ij}$.
		Hence 
\begin{equation}\label{real200}
	a_{\mathrm{i}z}^{ij}=\mathrm{i}a_z^{ij}.
\end{equation}		
		Furthermore, since $\Phi(ze_i\otimes e_j)
		=a_z^{ij}(e_i\otimes e_j)$,
		the map $z\mapsto a_z^{ij}$ from
		$\mathbb{C}$ into itself is real-linear.
		Together with \eqref{real200},
		this implies
		\[
		a_z^{ij}=za_1^{ij}
		\qquad(z\in\mathbb{C}).
		\]
		
		We now claim that $a_1^{ij}=a_1^{kl}$
		for all $i\neq j$ and $k\neq l$
		in $\{1,2,\ldots,n\}$.
		Indeed, by \eqref{scaequ1}--\eqref{scaequ3},
		we have
		\begin{equation}\label{scaequ5}
			a_1^{ij}=a_1^{ji}.
		\end{equation}
		Thus the claim is immediate when $n=2$.
		Suppose that $n\geq3$.
		For any three distinct indices $i,j,k$,
		since $e_i\perp(e_j+e_k)$,
		the hypothesis implies that there exists
		$\mu\in\mathbb{C}$ such that
		\[
		\begin{aligned}
			\mu e_i\otimes e_j
			+\mu e_i\otimes e_k=\mu e_i\otimes(e_j+e_k)&=\Phi(e_i\otimes(e_j+e_k))\\
			&=\Phi(e_i\otimes e_j)
			+\Phi(e_i\otimes e_k)\\
			&=a_1^{ij}(e_i\otimes e_j)
			+a_1^{ik}(e_i\otimes e_k).
		\end{aligned}
		\]
		Since $e_i\otimes e_j$ and
		$e_i\otimes e_k$ are linearly independent,
		we obtain $a_1^{ij}=a_1^{ik}=\mu$. Combining this with \eqref{scaequ5},
		we conclude that $a_1^{ij}$ is independent
		of the choice of $i\neq j$.
		This proves the claim.
		Set $\lambda=a_1^{ij}$. It follows from \eqref{scaequ1}
		and \eqref{scaequ4} that
		$$\Phi(ze_i\otimes e_j)=a_z^{ij}(e_i\otimes e_j)=za_1^{ij}(e_i\otimes e_j)=\lambda z(e_i\otimes e_j)$$
		and $$\Phi(z(e_i\otimes e_i-e_j\otimes e_j))
		=\lambda z(e_i\otimes e_i-e_j\otimes e_j).$$ Finally, since
		$\mathfrak{sl}_n(\mathbb{C})$
		is spanned over $\mathbb{C}$ by $\{e_i\otimes e_j:i\neq j\}$
		and $\{e_i\otimes e_i-e_n\otimes e_n:
		1\leq i\leq n-1\}$, the above identities, together with
		the real-linearity of $\Phi$, imply that $\Phi(T)=\lambda T$ for all $T\in\mathfrak{sl}_n(\mathbb{C})$.
		
		This completes the proof.
	\end{proof}
	
	To establish the next proposition, we need
	the following theorem.
	\begin{theorem}[{\cite[Theorem 4]{Sweet-laa-2009}}]\label{dimthm}
		Let $V$ be a complex-linear subspace of
		$\mathbb{M}_n(\mathbb{C})$ such that
		$T^2=0$ for every $T\in V$. Set $r=\max\{\operatorname{rank} (T):T\in V\}$. Then $\dim_{\mathbb{C}}V\leq r(n-r).
		$
	\end{theorem}
	\begin{proposition}\label{realsemrl}
		Let $\Phi:\mathfrak{sl}_n(\mathbb{C})
		\to\mathfrak{sl}_n(\mathbb{C})$
		be an invertible real-linear map such that
		$\Phi(A)^2=0$ whenever $A^2=0$.
		Then $\Phi$ is either complex-linear
		or conjugate-linear.
	\end{proposition}
	\begin{proof}
		Set $p=\lfloor \frac{n}{2}\rfloor$ and
		$q=\lceil \frac{n}{2}\rceil$.
		Let $V$ be a complex-linear subspace of
		$\mathbb{M}_n(\mathbb{C})$ consisting of
		square-zero matrices.
		For every $T\in V$, we have $T^2=0$,
		which implies $\operatorname{ran}T\subseteq\ker T$.
		By the rank-nullity theorem,
		$\operatorname{rank}T\leq p$.
		Hence, by Theorem \ref{dimthm},
		we obtain $\dim_{\mathbb{C}}V\leq pq$. Let $U$ be a subspace of $\mathbb{C}^n$
		with $\dim_{\mathbb{C}}U=p$ or $q$.
		Define $$S_U=\{A\in\mathbb{M}_n(\mathbb{C}):
		\operatorname{ran}A\subseteq U\subseteq\ker A\}.$$ Clearly, $A\in S_U$ if and only if $A=\begin{pmatrix}
			0&\widetilde{A}\\
			0&0
		\end{pmatrix}$ with respect to the orthogonal decomposition
		$\mathbb{C}^n=U\oplus U^\perp$,
		where $\widetilde{A}:U^\perp\to U$
		is a linear map. Therefore,
		\begin{equation}
			\dim_{\mathbb{C}}S_U=pq.
		\end{equation}
		By the hypothesis, $\Phi(S_U)$ is a
		real-linear subspace consisting of
		square-zero matrices.
		Thus, for any $A_k,A_l\in\Phi(S_U)$,
		we have
		\begin{equation}\label{squarezero}
			A_kA_l+A_lA_k
			=(A_k+A_l)^2-A_k^2-A_l^2=0.
		\end{equation}
		Let $V_U$ denote the complex-linear span
		of $\Phi(S_U)$. We claim that $\Phi(S_U)$ is a
		complex-linear subspace. Indeed, for any $B\in V_U$, we can write
		$B=\sum_k\alpha_k A_k$, where
		$\alpha_k\in\mathbb{C}$ and
		$A_k\in\Phi(S_U)$.
		By \eqref{squarezero}, we obtain $B^2=\sum_k\alpha_k^2 A_k^2+\sum_{k<l}\alpha_k\alpha_l(A_kA_l+A_lA_k)=0$. Hence, Theorem \ref{dimthm} applies to
		$V_U$, yielding
		$\dim_{\mathbb{C}}V_U\leq pq$.
		Since $\Phi$ is invertible, we have
			$$2pq=2\dim_{\mathbb{C}}S_U=\dim_{\mathbb{R}}S_U=\dim_{\mathbb{R}}\Phi(S_U)\leq\dim_{\mathbb{R}}V_U=2\dim_{\mathbb{C}}V_U\leq2pq.$$
Consequently,
		$\dim_{\mathbb{R}}\Phi(S_U)
		=\dim_{\mathbb{R}}V_U$.
		Since $\Phi(S_U)\subseteq V_U$,
		we conclude that $\Phi(S_U)=V_U$.
		Thus $\Phi(S_U)$ is a complex-linear
		subspace.
		
		Now let $u,v\in\mathbb{C}^n$ be
		arbitrary nonzero vectors satisfying
		$\langle u,v\rangle=0$.
		Since $p+q=n$, there exist subspaces
		$V_1,V_2$ of $\mathbb{C}^n$ such that $(\operatorname{span}_{\mathbb{C}}\{v\})^\perp
		=\operatorname{span}_{\mathbb{C}}\{u\}
		\oplus V_1\oplus V_2$,
		where
		$\dim_{\mathbb{C}}V_1=p-1$ and $\dim_{\mathbb{C}}V_2=q-1$.
		Set $U_1=\operatorname{span}_{\mathbb{C}}\{u\}
		\oplus V_1$ and $
		U_2=\operatorname{span}_{\mathbb{C}}\{u\}
		\oplus V_2$.
		For any nonzero
		$T\in S_{U_1}\cap S_{U_2}$,
		we have
		$
		\operatorname{ran}T
		\subseteq U_1\cap U_2
		=\operatorname{span}_{\mathbb{C}}\{u\}$.
		Thus $T$ has rank one, and there exists
		a nonzero vector $w\in\mathbb{C}^n$
		such that $T=u\otimes w$.
		Moreover, $(\operatorname{span}_{\mathbb{C}}\{v\})^\perp=U_1+U_2\subseteq\ker T=(\operatorname{span}_{\mathbb{C}}\{w\})^\perp$.
		It follows that
		$w\in\operatorname{span}_{\mathbb{C}}\{v\}$.
		Hence $S_{U_1}\cap S_{U_2}\subseteq W$, where
		$W=\operatorname{span}_{\mathbb{C}}\{u\otimes v\}$.
		Since the reverse inclusion is immediate,
		we obtain $S_{U_1}\cap S_{U_2}=W$. By the injectivity of $\Phi$, we have $\Phi(W)=\Phi(S_{U_1}\cap S_{U_2})=\Phi(S_{U_1})\cap\Phi(S_{U_2})$.
		Since both $\Phi(S_{U_1})$ and
		$\Phi(S_{U_2})$ are complex-linear
		subspaces, so is $\Phi(W)$.
		Consequently,
		$$\dim_{\mathbb{C}}\Phi(W)=\frac12\dim_{\mathbb{R}}\Phi(W)=\frac12\dim_{\mathbb{R}}W=\dim_{\mathbb{C}}W=1,$$
		where the second equality follows
		from the invertible real-linearity
		of $\Phi$.
		Since $\Phi(u\otimes v)\in\Phi(W)$
		and $\Phi(u\otimes v)\neq0$, we conclude that
		\[
		\Phi(W)=
		\operatorname{span}_{\mathbb{C}}
		\{\Phi(u\otimes v)\}.
		\]
		
		Finally, define a real-linear map
		$\Psi:\mathfrak{sl}_n(\mathbb{C})
		\to\mathfrak{sl}_n(\mathbb{C})$ by $\Psi(T)=\Phi^{-1}(\mathrm{i}\Phi(T))$.
		By the preceding identity, for any
		nonzero $u,v\in\mathbb{C}^n$ satisfying
		$\langle u,v\rangle=0$, there exists
		$\lambda_{u,v}\in\mathbb{C}$ such that
		$$\Psi(u\otimes v)=\Phi^{-1}(\mathrm{i}\Phi(u\otimes v))=\lambda_{u,v}(u\otimes v).$$
		Therefore, by Lemma \ref{ralema},
		there exists $\lambda\in\mathbb{C}$
		such that $\Psi(T)=\lambda T$ for all $T\in\in\mathfrak{sl}_n(\mathbb{C})$. On the other hand, the definition of
		$\Psi$ implies that
		$\Psi^2(T)=-T$ for all $T\in\mathfrak{sl}_n(\mathbb{C})$.
		Thus $\lambda^2=-1$, and hence
		$\lambda=\mathrm{i}$ or $-\mathrm{i}$. If $\lambda=\mathrm{i}$, then
		$\mathrm{i}\Phi(T)=\Phi(\mathrm{i}T)$
		for every $T\in\mathfrak{sl}_n(\mathbb{C})$,
		so $\Phi$ is complex-linear. If $\lambda=-\mathrm{i}$, then $\mathrm{i}\Phi(T)=\Phi(-\mathrm{i}T)$. By the real-linearity of $\Phi$, we obtain
		$\Phi(\mathrm{i}T)=-\mathrm{i}\Phi(T)$,
		which shows that $\Phi$ is conjugate-linear.
		
		This completes the proof.
	\end{proof}
	
	Using Theorem \ref{semrlprezero}, we obtain the following key result. For $1\leq i,j\leq n$, let $E_{ij}$ be the standard
	matrix unit with $1$ in the $(i,j)$-entry and
	$0$ elsewhere.
	
	\begin{theorem}\label{keyth}
		Let $\Phi:\mathfrak{sl}_n(\mathbb{C})
		\to\mathfrak{sl}_n(\mathbb{C})$
		be an invertible real-linear map such that
		$\Phi(T)^2=0$ whenever $T^2=0$,
		and $\Phi(N)$ is normal for every
		normal matrix $N\in\mathfrak{sl}_n(\mathbb{C})$.
		Then there exist a nonzero scalar
		$c\in\mathbb{C}$ and a unitary matrix $U$
		such that one of the following holds:
		\[
		\begin{gathered}
			\Phi(A)=cUAU^*
			\quad\text{or}\quad
			\Phi(A)=cU\overline{A}U^*,\\
			\Phi(A)=cUA^{\mathrm{tr}}U^\ast
			\quad\text{or}\quad
			\Phi(A)=cUA^*U^*
		\end{gathered}
		\]
		for all $A\in\mathfrak{sl}_n(\mathbb{C})$.
	\end{theorem}
	
	\begin{proof}
		By Lemma \ref{realsemrl}, $\Phi$ is either
		complex-linear or conjugate-linear. If $\Phi$ is complex-linear, then
		Theorem \ref{semrlprezero} implies that
		there exist a nonzero scalar
		$c\in\mathbb{C}$ and an invertible matrix
		$S\in\mathbb{M}_n(\mathbb{C})$ such that $\Phi(X)=cSXS^{-1}$ or $\Phi(X)=cSX^{\mathrm{tr}}S^{-1}$ for all $X\in\mathfrak{sl}_n(\mathbb{C})$. If $\Phi$ is conjugate-linear, then the map
		$\Psi:\mathfrak{sl}_n(\mathbb{C})
		\to\mathfrak{sl}_n(\mathbb{C})$
		defined by $\Psi(X)=\Phi(\overline{X})$ is invertible and complex-linear.
		Moreover, $\Psi(X)^2=0$ whenever $X^2=0$.
		Thus, by Theorem \ref{semrlprezero},
		there exist a nonzero scalar
		$c\in\mathbb{C}$ and an invertible matrix
		$S\in\mathbb{M}_n(\mathbb{C})$ such that
		$\Psi(X)=cSXS^{-1}$ or $
		\Psi(X)=cSX^{\mathrm{tr}}S^{-1}$
		for all $X\in\mathfrak{sl}_n(\mathbb{C})$.
		Consequently, $\Phi(X)=cS\overline{X}S^{-1}$ or $\Phi(X)=cSX^*S^{-1}$
		for all $X\in\mathfrak{sl}_n(\mathbb{C})$.
		
			Let $\tau$ be one of the maps
		\[
		X\mapsto X,\quad
		X\mapsto\overline{X},\quad
		X\mapsto X^{\mathrm{tr}},\quad
		X\mapsto X^*
		\]
		on $\mathfrak{sl}_n(\mathbb{C})$.
		Fix such a map $\tau$.
		It remains to show that if
		$\Phi(X)=cS\tau(X)S^{-1}$ for all $X\in\mathfrak{sl}_n(\mathbb{C})$,
		where $c\neq0$ and $S$ is invertible,
		then $S$ is a scalar multiple of
		a unitary matrix. Indeed, let $A\in\mathfrak{sl}_n(\mathbb{C})$
		be an arbitrary self-adjoint matrix.
		Then $\tau^{-1}(A)$ is normal and has
		trace zero. Since $\Phi$ preserves
		normal matrices, it follows that
		$SAS^{-1}$ is normal.
		Moreover, $\sigma(SAS^{-1})
		=\sigma(A)\subseteq\mathbb{R}$,
		which implies that $SAS^{-1}$ is self-adjoint. Hence
		$$
		S^*SA=AS^*S.$$
		Thus $S^*S$ commutes with every
		trace-zero self-adjoint matrix. For any $i\neq j$, both $
		E_{ij}+E_{ji}$ and $\mathrm{i}(E_{ji}-E_{ij})$ are trace-zero self-adjoint matrices.
		Since $E_{ij}
		=\frac12(E_{ij}+E_{ji})
		+\frac{\mathrm{i}}2
		\bigl(\mathrm{i}(E_{ji}-E_{ij})\bigr)$, we obtain $S^*SE_{ij}=E_{ij}S^*S$ for $i\ne j$.
		Consequently, $S^*S$ commutes with
		all off-diagonal matrix units, and hence $S^*S=\lambda I$
		for some $\lambda>0$.
		Therefore, $U=\frac{1}{\sqrt{\lambda}}S$ is unitary, as desired.
	\end{proof}

	Before stating our main result,
	we present one final lemma.
	
	\begin{lemma}\label{traceequ1}
		Let $f:[0,\infty)\to(0,\infty)$ be a
		continuous nonconstant function, and let
		$T\in\mathbb{M}_n(\mathbb{C})$.
		Suppose that
		\begin{equation}\label{qnorequ}
			f(k|T|)T=\frac{f(0)}{f(k)}Tf(k|T|)
		\end{equation}
		for every $k>0$.
		Then $\operatorname{tr}([T^*,T]^3)\geq0$,
		with equality if and only if $T^2=0$.
	\end{lemma}
	
	\begin{proof}
		The assertion is immediate if $T=0$.
		Henceforth, assume that $T\neq0$. Since $f$ is nonconstant, there exists
		$k_0>0$ such that $f(k_0)\neq f(0)$.
		By the spectral mapping theorem, $\sigma(f(k|T|))
		=\{f(k\lambda):\lambda\in\sigma(|T|)\}$.
Since $f$ is strictly positive,
		$f(k|T|)$ is invertible.
		It follows from \eqref{qnorequ} that
		$$\sigma(T)=\sigma\bigl(
			f(k_0|T|)T f(k_0|T|)^{-1}
			\bigr)=\sigma\left(\frac{f(0)}{f(k_0)}T\right)=\frac{f(0)}{f(k_0)}\sigma(T).
		$$
		Since $\frac{f(0)}{f(k_{0})}\neq1$, we obtain
		$\sigma(T)=\{0\}$.
		Thus $T$ is nilpotent.
		Let $m\geq2$ be its nilpotency index. We have
		$
		\{0\}=\ker T^0
		\subsetneq\ker T
		\subsetneq\ker T^2
		\subsetneq\cdots
		\subsetneq\ker T^m=\mathbb{C}^n.$
		Moreover,
		$T(\ker T^i)\subseteq\ker T^{i-1}$
		for every $1\leq i\leq m$. Set $\mathcal{H}_i=\ker T^i\ominus\ker T^{i-1}$
	for $1\leq i\leq m$, that is, $\ker T^i=\ker T^{i-1}\oplus\mathcal{H}_i$. Then, with respect to the orthogonal
		decomposition
		\begin{equation}\label{spacedec}
			\mathbb{C}^n
			=\mathcal{H}_1\oplus\mathcal{H}_2
			\oplus\cdots\oplus\mathcal{H}_m,
		\end{equation}
		the matrix $T$ has the form
		\begin{equation}\label{spaceequ1}
			T=
			\begin{pmatrix}
				0&T_{12}&T_{13}&\cdots&T_{1m}\\
				0&0&T_{23}&\cdots&T_{2m}\\
				\vdots&\vdots&\ddots&\ddots&\vdots\\
				0&0&\cdots&0&T_{(m-1)m}\\
				0&0&\cdots&0&0
			\end{pmatrix}.
		\end{equation}
On the other hand, by induction,
		\eqref{qnorequ} yields
		\begin{equation}\label{consact}
			T^if(k|T|)
			=\left(\frac{f(k)}{f(0)}\right)^i
			f(k|T|)T^i
		\end{equation}
		for every $i\in\mathbb{N}$ and $k>0$. Fix $1\leq i\leq m$ and
		$x\in\ker T^i$.
		By \eqref{consact}, we have $T^if(k|T|)x=\left(\frac{f(k)}{f(0)}\right)^if(k|T|)T^ix=0$. Hence $\ker T^i$ is invariant under
		the self-adjoint matrix $f(k|T|)$. Consequently, $\mathcal{H}_i$
		is also invariant under $f(k|T|)$. Now let $y\in\mathcal{H}_i$.
		By \eqref{consact},
		$$
		\begin{aligned}
			T^{i-1}f(k|T|)y=\left(\frac{f(k)}{f(0)}\right)^{i-1}
			f(k|T|)T^{i-1}y&=\left(\frac{f(k)}{f(0)}\right)^{i-1}
			\lim_{\ell\to\infty}
			p_\ell(|T|)T^{i-1}y\\
			&=\left(\frac{f(k)}{f(0)}\right)^{i-1}
			\lim_{\ell\to\infty}
			p_\ell(0)T^{i-1}y\\
			&=\left(\frac{f(k)}{f(0)}\right)^{i-1}
			f(0)T^{i-1}y,
		\end{aligned}
		$$
		where $\{p_\ell\}$ is a sequence of
		polynomials converging uniformly to
		$x\mapsto f(kx)$ on $[0,\|T\|]$.
		The third equality holds because
		$T^{i-1}y\in\ker T=\ker|T|$. It follows that
		\[
		T^{i-1}\left(
		f(k|T|)y
		-\left(\frac{f(k)}{f(0)}\right)^{i-1}
		f(0)y
		\right)=0.
		\]
		Since
		$\ker T^{i-1}\cap\mathcal{H}_i=\{0\}$
		and $\mathcal{H}_i$ is invariant
		under $f(k|T|)$, we obtain $f(k|T|)|_{\mathcal{H}_i}
		=\left(\frac{f(k)}{f(0)}\right)^{i-1}
		f(0)I_{\mathcal{H}_i}$.
		Therefore, with respect to
		\eqref{spacedec},
		\begin{equation}\label{spaceequ}
			f(k|T|)
		=\operatorname{diag}\Bigg(f(0)I_{\mathcal{H}_1},
		f(k)I_{\mathcal{H}_2},
		\ldots, \left(\frac{f(k)}{f(0)}\right)^{m-1}
		f(0)I_{\mathcal{H}_m}
		\Bigg).
		\end{equation}
		Substituting \eqref{spaceequ1} and
		\eqref{spaceequ} into \eqref{qnorequ},
		we obtain
		
		$$\Big(
		\left(\frac{f(k)}{f(0)}\right)^{i-1}
		-\left(\frac{f(k)}{f(0)}\right)^{j-2}
		\Big)T_{ij}=0$$
for all $1\leq i<j\leq m$ and $k>0$.
		Since $f$ is nonconstant,
		it follows that $T_{ij}=0$
		whenever $j\neq i+1$.
		Thus \eqref{spaceequ1} reduces to
		\begin{equation}\label{spacedec1}
			T=
			\begin{pmatrix}
				0&T_{12}&0&\cdots&0\\
				0&0&T_{23}&\cdots&0\\
				\vdots&\vdots&\ddots&\ddots&\vdots\\
				0&0&\cdots&0&T_{(m-1)m}\\
				0&0&\cdots&0&0
			\end{pmatrix}.
		\end{equation}
By \eqref{spacedec1}, we have $f(k|T|)=\operatorname{diag}\Big(f(0)I_{\mathcal{H}_1},f(k|T_{12}|),\ldots,f(k|T_{(m-1)m}|,f(k|T_{(m-2)(m-1)}|)\Big).$ Comparing this with \eqref{spaceequ},
		we obtain
		\begin{equation}\label{spec3}
			f(k|T_{(i-1)i}|)
			=\left(\frac{f(k)}{f(0)}\right)^{i-1}
			f(0)I_{\mathcal{H}_i}
		\end{equation}
		for every $2\leq i\leq m$. By the spectral mapping theorem,
		for any
		$\alpha,\beta\in\sigma(|T_{(i-1)i}|)$,
		we have $f(k\alpha)=f(k\beta)$ for $k>0$. Suppose that $0\leq\alpha<\beta$,
		and put $r=\alpha/\beta\in[0,1)$.
		Since $k>0$ is arbitrary, we obtain $f(rx)=f(x)$ for $ x>0$.
		Consequently, $f(r^\ell x)=f(r^{\ell-1}x)=\cdots=f(x)$ for every positive integer $\ell$.
		Letting $\ell\to\infty$ and using
		the continuity of $f$ at zero, we get $f(0)=f(x)$ for $x>0$, contrary to the nonconstancy of $f$.
		Thus $\alpha=\beta$, which implies
		that there exists $p_i\geq0$ such that
		\begin{equation}\label{spec1}
			|T_{(i-1)i}|=p_iI_{\mathcal{H}_i}.
		\end{equation}
	We claim that $p_i>0$. Indeed, it suffices to show that
		$\ker T_{(i-1)i}=\{0\}$.
		Let $x\in\mathcal{H}_i$ satisfy
		$T_{(i-1)i}x=0$.
		By \eqref{spacedec1}, we have $Tx=0$.
		Hence $x\in\ker T\subseteq\ker T^{i-1}$. Since
		$\mathcal{H}_i
		=\ker T^i\ominus\ker T^{i-1}$,
		we conclude that $x=0$.
		Therefore $p_i>0$. Set $U_i=p_i^{-1}T_{(i-1)i},
		\qquad 2\leq i\leq m$. By \eqref{spec1}, we have
		$
			U_i^*U_i=p_i^{-2}T_{(i-1)i}^*T_{(i-1)i}=p_i^{-2}|T_{(i-1)i}|^2=I_{\mathcal{H}_i}.
		$
		Thus each $U_i$ is an isometry.
		Consequently, \eqref{spacedec1}
		takes the form
		\begin{equation}\label{spacedec2}
			T=
			\begin{pmatrix}
				0&p_2U_2&0&\cdots&0\\
				0&0&p_3U_3&\cdots&0\\
				\vdots&\vdots&\ddots&\ddots&\vdots\\
				0&0&\cdots&0&p_mU_m\\
				0&0&\cdots&0&0
			\end{pmatrix}.
		\end{equation}
		For every $2\leq i\leq m$, we have
		$\operatorname{tr}(U_iU_i^*)=\operatorname{tr}(U_i^*U_i)=\operatorname{tr}(I_{\mathcal{H}_i})=\dim_{\mathbb{C}}\mathcal{H}_i>0$.
		Hence $P_i=U_iU_i^\ast$
		is a nonzero orthogonal projection. Put $A=T^*T$ and $B=TT^*$.
		By the cyclicity of the trace,
		we have $$\operatorname{tr}([T^*,T]^3)
		=3\operatorname{tr}(AB^2-A^2B).$$
		Moreover, a direct computation using
		\eqref{spacedec2} gives
		\[
		AB^2-A^2B=
		\begin{cases}
			0, & m=2,\\[4pt]
			\operatorname{diag}\Big(
			0,\,
			p_2^2p_3^2(p_3^2-p_2^2)P_3,\,
			\ldots, p_{m-1}^2p_m^2
			(p_m^2-p_{m-1}^2)P_m,\,
			0
			\Big),
			& m\geq3.
		\end{cases}
		\]
Finally, if we can show that
\begin{equation}\label{inseq}
	p_2<p_3<\cdots<p_m,
	\qquad m\in\{3,4,\ldots,n\},
\end{equation}
then we immediately obtain
$\operatorname{tr}(AB^2-A^2B)\geq0$,
and hence the desired conclusion. Indeed, suppose that there exist
$m_1,m_2\in\{2,3,\ldots,m\}$
with $m_1<m_2$ such that
$p_{m_1}\geq p_{m_2}$.
By \eqref{spec3} and \eqref{spec1},
we have
\[
f(kp_{m_1})
=\left(\frac{f(k)}{f(0)}\right)^{m_1-1}f(0)
\]
and
\[
f(kp_{m_2})
=\left(\frac{f(k)}{f(0)}\right)^{m_2-1}f(0).
\]
Hence $
\frac{f(kp_{m_1})}{f(0)}
=\left(\frac{f(kp_{m_2})}{f(0)}\right)^{
	\frac{m_1-1}{m_2-1}}.
$
Set
$\widetilde{r}=\frac{p_{m_2}}{p_{m_1}}\in(0,1]$
and $x=kp_{m_1}$.
Then $\frac{f(x)}{f(0)}
=\left(\frac{f(\widetilde{r}x)}{f(0)}\right)^{
	\frac{m_1-1}{m_2-1}}$. Since $k>0$ is arbitrary, so is $x>0$.
Thus, by induction, we obtain
\begin{equation}\label{induequ}
\begin{aligned}
	\frac{f(\widetilde{r}^{\ell}x)}{f(0)}=\left(
	\frac{f(\widetilde{r}^{\ell-1}x)}{f(0)}
	\right)^{\frac{m_2-1}{m_1-1}}=\left(
	\frac{f(\widetilde{r}^{\ell-2}x)}{f(0)}
	\right)^{\left(\frac{m_2-1}{m_1-1}\right)^2}=\cdots=\left(
	\frac{f(x)}{f(0)}
	\right)^{\left(\frac{m_2-1}{m_1-1}\right)^\ell}.
\end{aligned}
\end{equation}
If $\widetilde{r}\in(0,1)$,
letting $\ell\to\infty$ in the above identity
and using the continuity of $f$, we obtain
\[
1=\frac{f(0)}{f(0)}
=\lim_{\ell\to\infty}
\left(\frac{f(x)}{f(0)}\right)^{
	\left(\frac{m_2-1}{m_1-1}\right)^\ell}.
\]
Hence $f(x)=f(0)$. If $\widetilde{r}=1$,
letting $\ell\to\infty$ in \eqref{induequ},
we obtain
\[
\frac{f(x)}{f(0)}
=\lim_{\ell\to\infty}
\left(\frac{f(x)}{f(0)}\right)^{
	\left(\frac{m_2-1}{m_1-1}\right)^\ell}.
\]
Hence $f(x)=f(0)$. Thus, in either case, $f$ must be constant, a contradiction. Therefore, \eqref{inseq} holds.

This completes the proof.
	\end{proof}

 Set
\[
S_1^{ij}=E_{ij}+E_{ji},\quad
S_2^{ij}=\mathrm{i}(E_{ji}-E_{ij}),\quad
S_3^{ij}=E_{ii}-E_{jj},
\]
for $1\le i<j\le n$. Then, for any distinct $k,l\in\{1,2,3\}$,
we have
\begin{equation}\label{comequzero}
	(S_k^{ij})^2=E_{ii}+E_{jj},\quad
	S_k^{ij}S_l^{ij}=-S_l^{ij}S_k^{ij},\quad
	\operatorname{tr}(S_k^{ij}S_l^{ij})=0.
\end{equation}
For any $z\in\mathbb{C}$, write
$\Re z=\frac{z+\bar{z}}{2}$ and
$\Im z=\frac{z-\bar{z}}{2\mathrm{i}}$,
so that $z=\Re z+\mathrm{i}\Im z$.
Then, for every
$X=(x_{ij})\in\mathfrak{sl}_n(\mathbb{C})$,
we have
\begin{equation}\label{spanequ}
		X=\sum_{1\leq i<j\leq n}
	\Big((\Re z_1^{ij}+\mathrm{i}\Im z_1^{ij})S_1^{ij}+(\Re z_2^{ij}+\mathrm{i}\Im z_2^{ij})S_2^{ij}\Big)+\sum_{i=1}^{n-1}
	(\Re x_{ii}+\mathrm{i}\Im x_{ii})S_3^{in},
\end{equation}
where $z_1^{ij}=\frac{x_{ij}+x_{ji}}{2}$ and $
z_2^{ij}=\frac{x_{ji}-x_{ij}}{2\mathrm{i}}$. 

We are now ready to state the main result
of this section.
	
\begin{theorem}\label{kkeylem}
	Let $(f,g)\in C_{\Delta}^{+}([0,\infty))$
	with $f$ nonconstant, and let
	$\Phi:\mathbb{M}_n(\mathbb{C})\to
	\mathbb{M}_n(\mathbb{C})$
	be a surjective additive map satisfying
	$\Phi\circ\Delta_{f,g}
	=\Delta_{f,g}\circ\Phi$.
	Then the following statements hold:
	\begin{enumerate}[label=(\roman*)]
		\item $\Phi(N)$ is normal for every
		normal matrix $N$.
		
		\item There exists an additive map
		$\gamma:\mathbb{C}\to\mathbb{C}$
		such that
		\[
		\Phi(zI)=\gamma(z)I
		\qquad(z\in\mathbb{C}).
		\]
		
		\item The restriction
		$\Phi|_{\mathfrak{sl}_n(\mathbb{C})}$
		is an invertible real-linear map
		on $\mathfrak{sl}_n(\mathbb{C})$.
		
		\item $\Phi(T)^2=0$ whenever $T\in\mathbb{M}_n(\mathbb{C})$ satisfies $T^2=0$.
		
		\item There exist a unitary matrix $U$
		and a nonzero scalar $c\in\mathbb{C}$
		with $|c|\in G_{f,g}$ such that either
		\[
		\Phi(A)=cUAU^*
		\quad\text{or}\quad
		\Phi(A)=cU\overline{A}U^*
		\]
		for all $A\in\mathfrak{sl}_n(\mathbb{C})$.
		
		\item There exist a unitary matrix $U$
		and a nonzero scalar $c\in\mathbb{C}$
		with $|c|\in G_{f,g}$ such that either
		\[
		\Phi(A)=cUAU^*
		\quad\text{or}\quad
		\Phi(A)=cU\overline{A}U^*
		\]
		for all $A\in\mathbb{M}_n(\mathbb{C})$.
	\end{enumerate}
\end{theorem}
\begin{proof}
(i) Let $q$ be an arbitrary positive rational number
and $N$ an arbitrary normal matrix.
By Proposition \ref{basicprop2} (i), we have
$\Delta_{f,g}(qN)=qN$.
Since $\Phi$ is additive, we have
$\Phi(qN)=q\Phi(N)$.
Thus, from
$\Phi\circ\Delta_{f,g}=\Delta_{f,g}\circ\Phi$,
we obtain
$$\Delta_{f,g}(q\Phi(N))=\Delta_{f,g}(\Phi(qN))=\Phi(\Delta_{f,g}(qN))=\Phi(qN)=q\Phi(N).$$
By Proposition \ref{basicprop1} (vi) and the
above identity, we have
$\Delta_{f,g}(k\Phi(N))=k\Phi(N)$
for every $k>0$.
Hence, by Proposition \ref{basicprop2} (i),
$\Phi(N)$ is normal.

(ii)
Let $B\in\mathbb{M}_n(\mathbb{C})$ be arbitrary.
Since $\Phi$ is surjective, there exists
$A\in\mathbb{M}_n(\mathbb{C})$ such that
$\Phi(A)=B$.
Write $A=\Re(A)+\mathrm{i}\Im(A)$, where $\Re(A)=\frac{A+A^*}{2}$ and $
\Im(A)=\frac{A-A^*}{2\mathrm{i}}$. For any $z\in\mathbb{C}$, both
$zI+\Re(A)$ and $zI-\Re(A)$ are normal.
By (i) and the additivity of $\Phi$,
the matrices
$
\Phi(zI+\Re(A))
=\Phi(zI)+\Phi(\Re(A))$ and $
\Phi(zI-\Re(A))
=\Phi(zI)-\Phi(\Re(A))$
are normal. Consequently,
$$[\Phi(zI)+\Phi(\Re(A)),
(\Phi(zI))^*+(\Phi(\Re(A)))^*]=[\Phi(zI)+\Phi(\Re(A)),
(\Phi(zI)+\Phi(\Re(A)))^*]=0$$
and
$$[\Phi(zI)-\Phi(\Re(A)),
(\Phi(zI))^*-(\Phi(\Re(A)))^*]=[\Phi(zI)-\Phi(\Re(A)),
(\Phi(zI)-\Phi(\Re(A)))^*]=0.$$
Expanding and subtracting the above
two identities yields
\begin{equation}\label{sumequ1}
	[\Phi(zI),(\Phi(\Re(A)))^*]
	+[\Phi(\Re(A)),(\Phi(zI))^*]=0.
\end{equation}
Similarly, we obtain
\begin{equation}\label{sumequ2}
	[\Phi(zI),(\Phi(\mathrm{i}\Im(A)))^*]
	+[\Phi(\mathrm{i}\Im(A)),(\Phi(zI))^*]=0.
\end{equation}
Adding \eqref{sumequ1} and \eqref{sumequ2}
and using the additivity of $\Phi$, we obtain
\begin{equation}\label{sumequ3}
	[\Phi(zI),B^*]+[B,(\Phi(zI))^*]=[\Phi(zI),(\Phi(A))^*]+[\Phi(A),(\Phi(zI))^*]=0.
\end{equation}
Since $B$ is arbitrary, the above identity
also gives
$$-\mathrm{i}[\Phi(zI),B^*]+\mathrm{i}[B,(\Phi(zI))^*]=[\Phi(zI),(\mathrm{i}B)^*]+[\mathrm{i}B,(\Phi(zI))^*]=0.
$$
Hence
\begin{equation}\label{sumequ4}
	-[\Phi(zI),B^*]+[B,(\Phi(zI))^*]=0.
\end{equation}
Subtracting \eqref{sumequ4} from
\eqref{sumequ3}, we obtain
$[\Phi(zI),B^*]=0$
for every $B\in\mathbb{M}_n(\mathbb{C})$.
Thus $\Phi(zI)$ is a scalar multiple of $I$,
that is, there exists
$\gamma(z)\in\mathbb{C}$ such that
$\Phi(zI)=\gamma(z)I$. Finally, for any $z_1,z_2\in\mathbb{C}$,
the additivity of $\Phi$ gives $$\begin{aligned}
	\gamma(z_1+z_2)I=\Phi((z_1+z_2)I)=\Phi(z_1I)+\Phi(z_2I)=\gamma(z_1)I+\gamma(z_2)I.
\end{aligned}$$
Therefore, the map
$\gamma:\mathbb{C}\to\mathbb{C}$
is additive.

(iii) Let $1\leq i<j\leq n$.
If we can show that
\begin{equation}\label{reallinear1}
	\Phi(rS_k^{ij})=r\Phi(S_k^{ij})
\end{equation}
and
\begin{equation}\label{realinear3}
	\Phi(r\mathrm{i}S_k^{ij})
	=r\Phi(\mathrm{i}S_k^{ij})
\end{equation}
for every $r>0$ and $k\in\{1,2,3\}$,
then it is easy to verify from
\eqref{spanequ} that
$\Phi|_{\mathfrak{sl}_n(\mathbb{C})}$
is real-linear. Indeed, if $f(0)=0$, take any $\rho>0$.
If $f(0)>0$, we may choose $\rho>0$
such that $\frac{f(0)}{f(2\rho)}
\in\mathbb{Q}\setminus\{1\}$. This is possible because $f$ is nonconstant.
Thus there exists $r>0$ such that
$f(0)\neq f(r)$.
Consider the continuous function
$x\mapsto f(0)/f(x)$ on $[0,r]$.
By the density of the rational numbers
and the intermediate value theorem,
there exists $0<\rho<r/2$ such that
$f(0)/f(2\rho)\in\mathbb{Q}\setminus\{1\}$. Let $a,b\in\mathbb{R}$ satisfy
$a^2+b^2=\rho^2$, and set $N=aS_k^{ij}+bS_l^{ij}
+\mathrm{i}\rho S_m^{ij}$, where $k,l,m\in\{1,2,3\}$ are pairwise
distinct.
Clearly, $N$ is unitarily similar to
$\widetilde{N}\oplus0$, where
$\widetilde{N}$ is a $2\times2$ matrix.
A direct computation gives $N-N^*=2\mathrm{i}\rho S_m^{ij}\neq0$, so $\widetilde{N}\neq0$.
By \eqref{comequzero}, we have
\[
N^2=(a^2+b^2-\rho^2)
(E_{ii}+E_{jj})=0.
\]
Hence $\widetilde{N}^2=0$, and consequently
$\det\widetilde{N}=0$.
It follows that
$\operatorname{rank}N
=\operatorname{rank}\widetilde{N}=1$.
Thus $N$ is a rank-one matrix, and
there exist nonzero vectors
$u,v\in\mathbb{C}^n$ such that
$N=u\otimes v$. Moreover, by the definition of $N$
and \eqref{comequzero}, a direct
computation yields
\begin{equation}\label{trace1}
	\begin{aligned}
		N^*N&=a^2(S_k^{ij})^2
		+b^2(S_l^{ij})^2
		+\rho^2(S_m^{ij})^2+2\mathrm{i}\rho
		(aS_k^{ij}S_m^{ij}
		+bS_l^{ij}S_m^{ij})\\&=2\rho^2(E_{ii}+E_{jj})+2\mathrm{i}\rho
		(aS_k^{ij}S_m^{ij}
		+bS_l^{ij}S_m^{ij}).
	\end{aligned}
\end{equation}
Since $N=u\otimes v$, we also have
\begin{equation}\label{trace2}
	N^*N=(v\otimes u)(u\otimes v)
	=\|u\|^2v\otimes v.
\end{equation}
By \eqref{comequzero}, taking traces
on both sides of \eqref{trace1}
and \eqref{trace2} gives
$\operatorname{tr}(N^*N)=4\rho^2=\|u\|^2\operatorname{tr}(v\otimes v)=\|u\|^2\|v\|^2$.
Furthermore, $N^2=0$ implies that $\langle x,v\rangle
\langle u,v\rangle u=0$
for every $x\in\mathbb{C}^n$.
Hence $\langle u,v\rangle=0$.
By Proposition \ref{basicprop2} (iii),
we obtain
\begin{equation}\label{powerzero}
	\Delta_{f,g}(N)
	=\frac{f(0)}{f(\|u\|\|v\|)}N
	=\frac{f(0)}{f(2\rho)}N.
\end{equation}
If $f(0)=0$, then the above identity gives $\Delta_{f,g}(\Phi(N))=\Phi(\Delta_{f,g}(N))=\Phi(0)=0$. Hence, by Proposition \ref{basicprop2} (v),
we have $\Phi(N)^2=0$. If $f(0)>0$, then by \eqref{powerzero}
and Proposition \ref{basicprop1} (i), we have
$$
\sigma(\Phi(N))=\sigma(\Delta_{f,g}(\Phi(N)))=\sigma\left(
\Phi\left(\frac{f(0)}{f(2\rho)}N\right)
\right)=\sigma\left(
\frac{f(0)}{f(2\rho)}\Phi(N)
\right)=\frac{f(0)}{f(2\rho)}
\sigma(\Phi(N)),
$$
where the third equality follows from
the additivity of $\Phi$ and the fact that
$f(0)/f(2\rho)\in\mathbb{Q}$.
Thus $\sigma(\Phi(N))=\{0\}$,
so $\Phi(N)$ is nilpotent. Therefore, in either case, we have
$\Phi(N)^n=0$.
By the additivity of $\Phi$, it follows that
\[
\Big(
\Phi(N-\mathrm{i}\rho S_m^{ij})
+\Phi(\mathrm{i}\rho S_m^{ij})
\Big)^n
=\Phi(N)^n=0.
\]
Expanding the above identity gives
\[
\begin{aligned}
	\bigl(\Phi(N-\mathrm{i}\rho S_m^{ij})\bigr)^n=-\sum_{\substack{
			T_1,\ldots,T_n\in
			\{\Phi(N-\mathrm{i}\rho S_m^{ij}),
			\Phi(\mathrm{i}\rho S_m^{ij})\}\\
			(T_1,\ldots,T_n)\neq
			(\Phi(N-\mathrm{i}\rho S_m^{ij}),\ldots,
			\Phi(N-\mathrm{i}\rho S_m^{ij}))
	}}
	\prod_{k=1}^n T_k.
\end{aligned}
\]
Consequently,
\[
\begin{aligned}
	\left\|
	\bigl(\Phi(N-\mathrm{i}\rho S_m^{ij})\bigr)^n
	\right\|&\leq
	\sum_{r=1}^n\binom nr
	\|\Phi(N-\mathrm{i}\rho S_m^{ij})\|^{n-r}
	\|\Phi(\mathrm{i}\rho S_m^{ij})\|^r\\
	&=
	\Big(
	\|\Phi(N-\mathrm{i}\rho S_m^{ij})\|
	+\|\Phi(\mathrm{i}\rho S_m^{ij})\|
	\Big)^n-
	\|\Phi(N-\mathrm{i}\rho S_m^{ij})\|^n,
\end{aligned}
\]
where
$\binom nr=\frac{n!}{r!(n-r)!}$. Clearly, $N-\mathrm{i}\rho S_m^{ij}$
is normal.
By (i), its image
$\Phi(N-\mathrm{i}\rho S_m^{ij})$
is also normal.
Hence
$
\left\|
\bigl(\Phi(N-\mathrm{i}\rho S_m^{ij})\bigr)^n
\right\|
=
\|\Phi(N-\mathrm{i}\rho S_m^{ij})\|^n$.
It follows from the above inequality that
\begin{equation}\label{reallinear2}
\|\Phi(aS_k^{ij}+bS_l^{ij})\|=\|\Phi(N-\mathrm{i}\rho S_m^{ij})\|\leq
\frac{\|\Phi(\mathrm{i}\rho S_m^{ij})\|}
{2^{1/n}-1}.
\end{equation}
Similarly, we have $\|\Phi(aS_k^{ij}-bS_l^{ij})\|
\leq
\frac{\|\Phi(\mathrm{i}\rho S_m^{ij})\|}
{2^{1/n}-1}$. Thus, by the additivity of $\Phi$,
\[
\begin{aligned}
	\|2\Phi(aS_k^{ij})\|=\|\Phi(2aS_k^{ij})\|&=\|\Phi(aS_k^{ij}+bS_l^{ij}
	+aS_k^{ij}-bS_l^{ij})\|\\&\leq
	\|\Phi(aS_k^{ij}+bS_l^{ij})\|+
	\|\Phi(aS_k^{ij}-bS_l^{ij})\|\\&\leq
	2\frac{\|\Phi(\mathrm{i}\rho S_m^{ij})\|}
	{2^{1/n}-1}.
\end{aligned}
\]
Therefore, for every real number $a$
satisfying $|a|\leq\rho$, we have
\[
\|\Phi(aS_k^{ij})\|
\leq
\frac{\|\Phi(\mathrm{i}\rho S_m^{ij})\|}
{2^{1/n}-1}.
\]
Now consider the additive map
$\widehat{\Phi}:\mathbb{R}\to
\mathbb{M}_n(\mathbb{C})$
defined by
$\widehat{\Phi}(a)=\Phi(aS_k^{ij})$.
For any $\varepsilon>0$, choose
a positive integer $M$ such that $\frac{\|\Phi(\mathrm{i}\rho S_m^{ij})\|}
{M(2^{1/n}-1)}<\varepsilon$.
Then, whenever $|a|\leq\rho/M$,
we have $
\begin{aligned}
	\|\widehat{\Phi}(a)\|=\frac{1}{M}\|\widehat{\Phi}(Ma)\|\leq
	\frac{\|\Phi(\mathrm{i}\rho S_m^{ij})\|}
	{M(2^{1/n}-1)}
	<\varepsilon.
\end{aligned}$
Thus $\widehat{\Phi}$ is continuous at zero,
and hence continuous on $\mathbb{R}$. For any $t\in\mathbb{R}$,
choose a sequence of rational numbers
$\{t_l\}$ such that $t_l\to t$.
Then
\[
\begin{aligned}
	\Phi(tS_k^{ij})=\widehat{\Phi}(t)=\lim_{l\to\infty}
	\widehat{\Phi}(t_l)=\lim_{l\to\infty}
	\Phi(t_lS_k^{ij})=\lim_{l\to\infty}
	t_l\Phi(S_k^{ij})=t\Phi(S_k^{ij}).
\end{aligned}
\]
Thus \eqref{reallinear1} holds. Replacing $N$ by $\mathrm{i}N$
in the above argument, we similarly
obtain
\[
\begin{aligned}
	&\|\Phi(\mathrm{i}
	(aS_k^{ij}+bS_l^{ij}))\|=\|\Phi(\mathrm{i}N+\rho S_m^{ij})\|\leq
	\frac{\|\Phi(\rho S_m^{ij})\|}
	{2^{1/n}-1}.
\end{aligned}
\]
Continuing with the same argument
yields \eqref{realinear3}.
Therefore,
$\Phi|_{\mathfrak{sl}_n(\mathbb{C})}$
is real-linear.

Now take $a=\rho$ and $b=0$
in the definition of $N$.
Since $\Phi(N)$ is nilpotent
and $\Phi$ is real-linear, we have
\begin{equation}\label{realtrace1}
	\begin{aligned}
		0=\operatorname{tr}(\Phi(N))=\operatorname{tr}
		(\Phi(\rho S_k^{ij}+\mathrm{i}\rho S_m^{ij}))=\rho\operatorname{tr}(\Phi(S_k^{ij}))
		+\rho\operatorname{tr}
		(\Phi(\mathrm{i}S_m^{ij})).
	\end{aligned}
\end{equation}
Moreover, replacing $N$ by $N^*$,
it is easy to verify that
\eqref{powerzero} still holds.
Hence $\Phi(N^*)^2=0$, and consequently
\begin{equation}\label{realtrace2}
	\begin{aligned}
		0=\operatorname{tr}(\Phi(N^*))=\operatorname{tr}
		(\Phi(\rho S_k^{ij}
		-\mathrm{i}\rho S_m^{ij}))=\rho\operatorname{tr}(\Phi(S_k^{ij}))
		-\rho\operatorname{tr}
		(\Phi(\mathrm{i}S_m^{ij})).
	\end{aligned}
\end{equation}
By \eqref{realtrace1}, \eqref{realtrace2},
and $\rho>0$, we obtain $\operatorname{tr}(\Phi(S_k^{ij}))
=\operatorname{tr}
(\Phi(\mathrm{i}S_m^{ij}))=0$
for all $k,m\in\{1,2,3\}$.
Therefore, by \eqref{spanequ}
and the linearity of the trace,
we have $\operatorname{tr}(\Phi(X))=0$ for every
$X\in\mathfrak{sl}_n(\mathbb{C})$.
Thus $\Phi|_{\mathfrak{sl}_n(\mathbb{C})}:
\mathfrak{sl}_n(\mathbb{C})
\to\mathfrak{sl}_n(\mathbb{C})$ is well-defined.

Finally, we show that
$\Phi|_{\mathfrak{sl}_n(\mathbb{C})}$
is invertible.
Since it is a real-linear map on
the finite-dimensional real vector space
$\mathfrak{sl}_n(\mathbb{C})$,
it suffices to prove that it is surjective. Let $Y\in\mathfrak{sl}_n(\mathbb{C})$
be arbitrary.
Since $\Phi$ is surjective,
there exists
$X\in\mathbb{M}_n(\mathbb{C})$
such that $\Phi(X)=Y$.
Set $\widetilde{X}
=X-\frac{\operatorname{tr}(X)}n I
\in\mathfrak{sl}_n(\mathbb{C})$. By the additivity of $\Phi$
and (ii), we have
\[
\begin{aligned}
	\Phi(\widetilde{X})=\Phi\left(
	X-\frac{\operatorname{tr}(X)}n I
	\right)=\Phi(X)+
	\gamma\left(
	-\frac{\operatorname{tr}(X)}n
	\right)I=Y+
	\gamma\left(
	-\frac{\operatorname{tr}(X)}n
	\right)I.
\end{aligned}
\]
Taking traces on both sides
and using
$\Phi(\mathfrak{sl}_n(\mathbb{C}))
\subseteq\mathfrak{sl}_n(\mathbb{C})$,
we obtain
$
0=0+n\gamma\left(
-\frac{\operatorname{tr}(X)}n
\right)$. Thus $\gamma\left(
-\frac{\operatorname{tr}(X)}n
\right)=0$. Consequently,
$\Phi(\widetilde{X})=Y$,
as desired.

(iv) The assertion is immediate if $T=0$. Hence assume that $T\neq0$ and $T^2=0$. If $f(0)=0$, then by Proposition \ref{basicprop2} (v), we have $\Delta_{f,g}(T)=0$.
It follows that $\Delta_{f,g}(\Phi(T))=\Phi(\Delta_{f,g}(T))=\Phi(0)=0$. Hence, again by Proposition \ref{basicprop2} (v), we obtain $\Phi(T)^2=0$.

Suppose now that $f(0)>0$.
We first claim that
$\Phi(W)^2=0$
for every partial isometry $W$
satisfying $W^2=0$. Indeed, by Proposition \ref{basicprop2} (iv),
we have
\begin{equation}\label{pariso1}
	\Delta_{f,g}(kW)
	=f(0)g(k)W
	=f(0)\frac{k}{f(k)}W
\end{equation}
for every $k>0$.
Since
$W\in\mathfrak{sl}_n(\mathbb{C})$
and
$\Phi|_{\mathfrak{sl}_n(\mathbb{C})}$
is real-linear, the above identity gives
\[
\begin{aligned}
	\Delta_{f,g}(k\Phi(W))=\Delta_{f,g}(\Phi(kW))=\Phi(\Delta_{f,g}(kW))=k\frac{f(0)}{f(k)}\Phi(W).
\end{aligned}
\]
On the other hand,
since $|k\Phi(W)|=k|\Phi(W)|$, Proposition \ref{basicprop1} (iv) yields $$\Delta_{f,g}(k\Phi(W))=f(k|\Phi(W)|)\,k\Phi(W)(f(k|\Phi(W)|))^{-1}.$$
Combining the above two identities,
we obtain $f(k|\Phi(W)|)\Phi(W)=\frac{f(0)}{f(k)}
	\Phi(W)f(k|\Phi(W)|)$.
Therefore, by Lemma \ref{traceequ1},
we have
\begin{equation}\label{traceequ2}
	\operatorname{tr}
	([(\Phi(W))^*,\Phi(W)]^3)\geq0.
\end{equation}
Since $(W^*)^2=0$ and $W^*$
is also a partial isometry,
the same argument gives
\begin{equation}\label{traceequ3}
	\operatorname{tr}
	([(\Phi(W^*))^*,\Phi(W^*)]^3)\geq0.
\end{equation}
Moreover, both $W+W^*$ and $W-W^*$
are normal.
By (i) and the additivity of $\Phi$,
we have
\[
\begin{aligned}
	0=[(\Phi(W)+\Phi(W^*))^*,
	\Phi(W)+\Phi(W^*)]+[(\Phi(W)-\Phi(W^*))^*,
	\Phi(W)-\Phi(W^*)].
\end{aligned}
\]
Expanding the above identity yields $[(\Phi(W))^*,\Phi(W)]=-[(\Phi(W^*))^*,\Phi(W^*)]$. Thus, by \eqref{traceequ2},
\eqref{traceequ3}, and the above identity,
we obtain
\[
\begin{aligned}
	0\leq
	\operatorname{tr}
	([(\Phi(W))^*,\Phi(W)]^3)=
	\operatorname{tr}
	(-[(\Phi(W^*))^*,\Phi(W^*)]^3)\leq0.
\end{aligned}
\]
Hence $\operatorname{tr}
([(\Phi(W))^*,\Phi(W)]^3)=0$. By Lemma \ref{traceequ1},
we obtain the desired conclusion
$\Phi(W)^2=0$. Now consider the polar decomposition
$T=U|T|$ of a matrix $T$
satisfying $T^2=0$.
Clearly,
\begin{equation}\label{polarpro1}
	\operatorname{ran}U
	=\operatorname{ran}T
	\subseteq\ker T=\ker|T|
\end{equation}
and
\begin{equation}\label{polarpro2}
	U^*U=P_{\operatorname{ran}|T|}.
\end{equation}
Write the spectral decomposition
$|T|=\sum_{i=1}^m s_iP_i$,
where
$0<s_1<s_2<\cdots<s_m$
and, for each $i\in\{1,2,\ldots,m\}$,
$P_i$ is the orthogonal projection onto
$\ker(|T|-s_iI)$. Since
$\ker(|T|-s_iI)
\subseteq(\ker|T|)^\perp$
for every $i\in\{1,2,\ldots,m\}$, it follows from \eqref{polarpro1} that
$P_iU=0$. Hence
$
(UP_i)^2=0$
and
$
\begin{aligned}
	(U(P_i+P_j))^2=U(P_i+P_j)U(P_i+P_j)=0.
\end{aligned}
$
By \eqref{polarpro2}, for any distinct
$i,j\in\{1,2,\ldots,m\}$,
we have $(UP_i)^*UP_i=P_iU^*UP_i=P_iP_{\operatorname{ran}|T|}P_i=P_i$
and $	(U(P_i+P_j))^*U(P_i+P_j)=(P_i+P_j)
P_{\operatorname{ran}|T|}
(P_i+P_j)=P_i+P_j$. Since $P_i$ and $P_i+P_j$
are orthogonal projections,
both $UP_i$ and $U(P_i+P_j)$
are partial isometries.
By the assertion proved above,
we have
\[
\Phi(UP_i)^2=0,
\quad
\Phi(U(P_i+P_j))^2=0
\]
for any distinct
$i,j\in\{1,2,\ldots,m\}$. By the additivity of $\Phi$,
expanding the second identity gives
$\Phi(UP_i)\Phi(UP_j)+\Phi(UP_j)\Phi(UP_i)=0$.
Therefore, by (iii),
\[
\begin{aligned}
	\Phi(T)^2
	&=\Phi\left(
	\sum_{i=1}^m s_iUP_i
	\right)^2=\left(
	\sum_{i=1}^m s_i\Phi(UP_i)
	\right)^2\\
	&=\sum_{i=1}^m
	s_i^2\Phi(UP_i)^2+
	\sum_{1\leq i<j\leq m}
	s_is_j
	\bigl(
	\Phi(UP_i)\Phi(UP_j)+
	\Phi(UP_j)\Phi(UP_i)
	\bigr)\\
	&=0.
\end{aligned}
\]
This completes the proof of (iv).

(v) By (i), (iii), (iv), and
Theorem \ref{keyth},
there exist a nonzero scalar
$c\in\mathbb{C}$ and a unitary matrix $U$
such that one of the following holds:
\[
\begin{gathered}
	\Phi(X)=cUXU^*
	\quad\text{or}\quad
	\Phi(X)=cU\overline{X}U^*,\\
	\Phi(X)=cUX^{\mathrm{tr}}U^*
	\quad\text{or}\quad
	\Phi(X)=cUX^*U^*
\end{gathered}
\]
for all $X\in\mathfrak{sl}_n(\mathbb{C})$.

When $n=2$, for every
$X\in\mathfrak{sl}_2(\mathbb{C})$,
we have
$
X^{\mathrm{tr}}=-JXJ^*,
$
where
$
J=
\begin{pmatrix}
	0&-1\\
	1&0
\end{pmatrix}
$
is unitary.
Thus, if
$\Phi(X)=cUX^{\mathrm{tr}}U^*$,
then
$\Phi(X)=-cVXV^*$.
If
$\Phi(X)=cUX^*U^*$,
then
\[
\begin{aligned}
	\Phi(X)=cU(\overline{X})^{\mathrm{tr}}U^*=-cV\overline{X}V^*,
\end{aligned}
\]
where $V=UJ$ is unitary.
Therefore, (v) holds when $n=2$.

Suppose now that $n\geq3$.
We claim that neither $\Phi(X)=cUX^{\mathrm{tr}}U^*$
nor $
\Phi(X)=cUX^*U^*$ can hold for all
$X\in\mathfrak{sl}_n(\mathbb{C})$. Indeed, suppose that
$\Phi(X)=cUX^{\mathrm{tr}}U^*$.
By the identity
$\Phi\circ\Delta_{f,g}
=\Delta_{f,g}\circ\Phi$ and Proposition \ref{basicprop1} (iii),
we have
$\begin{aligned}
	cU(\Delta_{f,g}(X))^{\mathrm{tr}}U^*=\Delta_{f,g}(cUX^{\mathrm{tr}}U^*)=U\Delta_{f,g}(cX^{\mathrm{tr}})U^*.
\end{aligned}$
Consequently,
\begin{equation}\label{noclude}
	c(\Delta_{f,g}(X))^{\mathrm{tr}}
	=\Delta_{f,g}(cX^{\mathrm{tr}}).
\end{equation}
Let $s,t>0$ be arbitrary,
and consider
$
A=
\begin{pmatrix}
	0&s&0\\
	0&0&t\\
	0&0&0
\end{pmatrix}
\oplus0$.
A direct computation gives
\[
\Delta_{f,g}(A)=
\begin{pmatrix}
	0&f(0)g(s)&0\\
	0&0&f(s)g(t)\\
	0&0&0
\end{pmatrix}
\oplus0
\]
and
\[
\Delta_{f,g}(cA^{\mathrm{tr}})
=\frac{c}{|c|}
\begin{pmatrix}
	0&0&0\\
	f(|c|t)g(|c|s)&0&0\\
	0&f(0)g(|c|t)&0
\end{pmatrix}
\oplus0.
\]
Thus \eqref{noclude} implies
$
f(0)g(|c|t)=|c|f(s)g(t)$. If $f(0)=0$, this contradicts
the positivity of the right-hand side.
If $f(0)>0$, fixing $t$ gives
$
f(s)=
\frac{f(0)g(|c|t)}{|c|g(t)}
$
for every $s>0$.
Hence $f$ is constant, a contradiction.
Therefore,
$\Phi(X)=cUX^{\mathrm{tr}}U^*$
cannot hold for all
$X\in\mathfrak{sl}_n(\mathbb{C})$.

Moreover, since
$A^{\mathrm{tr}}=A^*$
and
$
(\Delta_{f,g}(A))^{\mathrm{tr}}
=(\Delta_{f,g}(A))^*$, the preceding computation immediately
shows that
$\Phi(X)=cUX^*U^*$
cannot hold for all
$X\in\mathfrak{sl}_n(\mathbb{C})$. Finally, since
$\Phi(X)=cUXU^*$
or
$\Phi(X)=cU\overline{X}U^*$,
let $x,b>0$ be arbitrary and consider $B=
\begin{pmatrix}
	0&b\\
	x&0
\end{pmatrix}
\oplus0
\in\mathfrak{sl}_n(\mathbb{C})$. A direct computation, similar to that
for the matrix $A$ above, gives
$
\Delta_{f,g}(|c|B)
=|c|\Delta_{f,g}(B)$,
which implies
$$\begin{pmatrix}
	0&f(|c|x)g(|c|b)\\
	f(|c|b)g(|c|x)&0
\end{pmatrix}
\oplus0=
|c|
\begin{pmatrix}
	0&f(x)g(b)\\
	f(b)g(x)&0
\end{pmatrix}
\oplus0.$$
Hence $
f(|c|x)g(|c|b)
=|c|f(x)g(b)$. Since
$(f,g)\in C_{\Delta}^{+}([0,\infty))$,
we obtain
$
f(|c|x)
=\frac{f(|c|b)}{f(b)}f(x)$. Therefore, by Proposition \ref{basicprop3} (i),
we conclude that
$|c|\in G_{f,g}$.
This completes the proof of (v).

(vi) By (v), there exist a scalar
$c\in\mathbb{C}$ with
$|c|\in G_{f,g}$ and a unitary matrix $U$
such that
\[
\Phi(X)=cUXU^*
\quad\text{or}\quad
\Phi(X)=cU\overline{X}U^*
\]
for all $X\in\mathfrak{sl}_n(\mathbb{C})$.

Suppose that
$\Phi(X)=cUXU^*$
(resp. $\Phi(X)=cU\overline{X}U^*$)
for every
$X\in\mathfrak{sl}_n(\mathbb{C})$.
Define the map
$\Psi:\mathbb{M}_n(\mathbb{C})
\to\mathbb{M}_n(\mathbb{C})$
by
$
\Psi(Y)=cUYU^*
$
(resp.
$\Psi(Y)=cU\overline{Y}U^*$).
Its inverse is given by
$
\Psi^{-1}(Y)=\frac1c U^*YU
$
(resp.
$\Psi^{-1}(Y)=
\frac1{\overline{c}}\,
\overline{U^*YU}$).
Clearly,
\begin{equation}\label{conide}
	(\Psi^{-1}\circ\Phi)(X)=X
\end{equation}
for every
$X\in\mathfrak{sl}_n(\mathbb{C})$. Moreover, since $|c|\in G_{f,g}$,
we have
$
\Psi\circ\Delta_{f,g}
=\Delta_{f,g}\circ\Psi$, and consequently
\begin{equation*}
	\Psi^{-1}\circ\Delta_{f,g}
	=\Delta_{f,g}\circ\Psi^{-1}.
\end{equation*}
Since $\Phi$ is a surjective additive map,
the composition
$\Psi^{-1}\circ\Phi$
is also a surjective additive map
on $\mathbb{M}_n(\mathbb{C})$.
Furthermore, the above identity gives
\begin{equation}\label{commuequ}
	\begin{aligned}
		(\Psi^{-1}\circ\Phi)\circ\Delta_{f,g}=\Psi^{-1}\circ\Delta_{f,g}\circ\Phi=\Delta_{f,g}\circ(\Psi^{-1}\circ\Phi).
	\end{aligned}
\end{equation}
Hence, by (ii), there exists
an additive map
$\beta:\mathbb{C}\to\mathbb{C}$
such that
\begin{equation}\label{conide1}
	(\Psi^{-1}\circ\Phi)(zI)=\beta(z)I
\end{equation}
for every $z\in\mathbb{C}$. Let $k>0$ be arbitrary, and set
$
T=
\begin{pmatrix}
	0&k\\
	0&0
\end{pmatrix}
\oplus0
\in\mathfrak{sl}_n(\mathbb{C})$. For any $z\in\mathbb{C}$,
by Proposition \ref{basicprop1} (ii)
and the linearity of the trace,
we have $
\Delta_{f,g}(zI+T)-zI
\in\mathfrak{sl}_n(\mathbb{C})$. Therefore, by \eqref{conide},
\eqref{commuequ}, \eqref{conide1},
and the additivity of
$\Psi^{-1}\circ\Phi$, we obtain
\[
\begin{aligned}
	\Delta_{f,g}(zI+T)-zI&=(\Psi^{-1}\circ\Phi)
	(\Delta_{f,g}(zI+T)-zI)\\
	&=(\Psi^{-1}\circ\Phi)
	(\Delta_{f,g}(zI+T))
	-\beta(z)I\\
	&=\Delta_{f,g}\circ\Psi^{-1}\circ\Phi
	(zI+T)-\beta(z)I\\
	&=\Delta_{f,g}\bigl(
	(\Psi^{-1}\circ\Phi)(zI)
	+(\Psi^{-1}\circ\Phi)(T)
	\bigr)-\beta(z)I\\
	&=\Delta_{f,g}(\beta(z)I+T)
	-\beta(z)I.
\end{aligned}
\]
Hence, by Proposition \ref{basicprop3} (iii),
we obtain
\begin{equation}\label{betafun}
	\beta(z)=z.
\end{equation}
Finally, let
$Y\in\mathbb{M}_n(\mathbb{C})$
be arbitrary and write
$
Y=
Y-\frac{\operatorname{tr}(Y)}n I
+\frac{\operatorname{tr}(Y)}n I$. Since $
Y-\frac{\operatorname{tr}(Y)}n I
\in\mathfrak{sl}_n(\mathbb{C})$, by the additivity of
$\Psi^{-1}\circ\Phi$ and $\Phi$,
together with \eqref{conide},
\eqref{conide1}, and \eqref{betafun},
we have
\[
\begin{aligned}
	(\Psi^{-1}\circ\Phi)(Y)&=(\Psi^{-1}\circ\Phi)
	\left(
	Y-\frac{\operatorname{tr}(Y)}n I
	+\frac{\operatorname{tr}(Y)}n I
	\right)\\
	&=\Psi^{-1}\left(
	\Phi\left(
	Y-\frac{\operatorname{tr}(Y)}n I
	\right)
	+\Phi\left(
	\frac{\operatorname{tr}(Y)}n I
	\right)
	\right)\\
	&=\Psi^{-1}\left(
	\Phi\left(
	Y-\frac{\operatorname{tr}(Y)}n I
	\right)
	\right)+\Psi^{-1}\left(
	\Phi\left(
	\frac{\operatorname{tr}(Y)}n I
	\right)
	\right)\\&=Y-\frac{\operatorname{tr}(Y)}n I
	+\beta\left(
	\frac{\operatorname{tr}(Y)}n
	\right)I=Y.
\end{aligned}
\]
Therefore,
$\Phi(Y)=\Psi(Y)$
for every
$Y\in\mathbb{M}_n(\mathbb{C})$.
This completes the proof.
\end{proof}

\section{Proof of the Main Theorem}
	
	To derive our main theorem, Theorem \ref{thm:main},
	from Theorem \ref{kkeylem} (vi), we first
	show in this section that every map $\Phi$
	satisfying \eqref{eq:weighted} is additive. The following recent result of Costara et al.
	plays a key role in the proof. For $T=(t_{ij})\in\mathbb{M}_n(\mathbb{C})$,
	write $T^\tau=(\tau(t_{ij}))\in
	\mathbb{M}_n(\mathbb{C})$,
	where $\tau:\mathbb{C}\to\mathbb{C}$ is a map. 

\begin{theorem}[{\cite[Theorem 2.10]{Costara-laa-2023}}]\label{thm:costara}
	Let $n\geq2$. Suppose that
	$\varphi_1,\varphi_2:
	\mathbb{M}_n(\mathbb{C})
	\to\mathbb{M}_n(\mathbb{C})$
	are surjective maps satisfying
	\[
	\det(A+B)=0
	\iff
	\det(\varphi_1(A)+\varphi_2(B))=0
	\]
	for all $A,B\in\mathbb{M}_n(\mathbb{C})$.
	Then there exist
	$R\in\mathbb{M}_n(\mathbb{C})$
	and a map
	$J:\mathbb{M}_n(\mathbb{C})
	\to\mathbb{M}_n(\mathbb{C})$
	such that
	\[
	\varphi_1(A)=J(A)+R,\quad
	\varphi_2(A)=J(A)-R
	\]
	for all $A\in\mathbb{M}_n(\mathbb{C})$,
	where $J$ has one of the forms
	$
	J(A)=PA^\tau Q$ or $
	J(A)=P(A^\tau)^{\mathrm{tr}}Q
	$
	for all $A\in\mathbb{M}_n(\mathbb{C})$,
	with $P,Q$ invertible and $\tau$
	a field automorphism of $\mathbb{C}$.
\end{theorem}

The following proposition allows us to
establish Theorem \ref{addthm}, and hence the additivity of $\Phi$.

\begin{proposition}\label{grouppro}
	Let
	$\{R_T\}_{T\in\mathbb{M}_n(\mathbb{C})}$
	and
	$\{S_T\}_{T\in\mathbb{M}_n(\mathbb{C})}$
	be families of self-maps on
	$\mathbb{M}_n(\mathbb{C})$
	satisfying the following conditions:
	\begin{enumerate}[label=(\roman*)]
		\item $R_0=S_0=\operatorname{id}$ (the identity map);
		
		\item
		$R_{T_1+T_2}=R_{T_1}R_{T_2}$
		and
		$S_{T_1+T_2}=S_{T_1}S_{T_2}$
		for all
		$T_1,T_2\in\mathbb{M}_n(\mathbb{C})$;
		
		\item
		$\sigma(R_T(X)+S_{-T}(Y))=\sigma(X+Y)$
		for all
		$T,X,Y\in\mathbb{M}_n(\mathbb{C})$;
		
		\item
		$\{R_T(0):T\in\mathbb{M}_n(\mathbb{C})\}
		=\mathbb{M}_n(\mathbb{C})$.
	\end{enumerate}
	Then, for every
	$T\in\mathbb{M}_n(\mathbb{C})$,
	there exists a unique
	$\widetilde{T}\in\mathbb{M}_n(\mathbb{C})$
	such that
	$
	R_T(X)=S_T(X)=X+\widetilde{T}
	$
	for all $X\in\mathbb{M}_n(\mathbb{C})$.
	Moreover, $
	\widetilde{T_1+T_2}
	=\widetilde{T_1}+\widetilde{T_2}$.
\end{proposition}

\begin{proof}
	Fix $T\in\mathbb{M}_n(\mathbb{C})$.
	It follows immediately from (i) and (ii)
	that, for every
	$T\in\mathbb{M}_n(\mathbb{C})$,
	the maps $R_T$ and $S_T$ are invertible
	maps on $\mathbb{M}_n(\mathbb{C})$,
	with inverses $R_{-T}$ and $S_{-T}$,
	respectively. Moreover, by (iii), for all
	$X,Y\in\mathbb{M}_n(\mathbb{C})$,
	the matrix $X+Y$ is singular if and only if
	$R_T(X)+S_{-T}(Y)$ is singular.
	Hence, by Theorem \ref{thm:costara}, there exist
	$\widetilde{T}\in\mathbb{M}_n(\mathbb{C})$
	and a map
	$J_T:\mathbb{M}_n(\mathbb{C})
	\to\mathbb{M}_n(\mathbb{C})$
	such that
	\begin{equation}\label{basform}
		R_T(X)=J_T(X)+\widetilde{T},
		\quad
		S_{-T}(X)=J_T(X)-\widetilde{T}
	\end{equation}
	for all $X\in\mathbb{M}_n(\mathbb{C})$,
	where $
	J_T(X)=A_TX^\tau B_T$ or $J_T(X)=A_T(X^\tau)^{\mathrm{tr}}B_T$, $A_T,B_T$ are invertible matrices,
	and $\tau$ is an automorphism of $\mathbb{C}$. By \eqref{basform}, taking $Y=0$ in (iii),
	we obtain
	\begin{equation}\label{speequ}
		\sigma(J_T(X))=\sigma(X)
	\end{equation}
	for every $X\in\mathbb{M}_n(\mathbb{C})$.
	Taking $X=I$ in \eqref{speequ}, we have $\sigma(A_TB_T)=\sigma(I)=\{1\}$. For any $z\in\mathbb{C}$,
	taking $X=zI$ in \eqref{speequ} gives
	$$	\{\tau(z)\lambda:\lambda\in\sigma(A_TB_T)\}=\sigma(\tau(z)A_TB_T)=\sigma(zI)=\{z\}.$$
	Hence $\tau(z)=z$, that is, $\tau$
	is the identity automorphism of
	$\mathbb{C}$. We now consider the following two cases.
	
	\noindent\textbf{Case 1.}
	Suppose that $J_T(X)=A_TXB_T$.
	By \eqref{speequ} and Jacobson's lemma,
	we have
	\[
	\sigma(B_TA_TX)
	=\sigma(A_TXB_T)
	=\sigma(X)
	\]
	for every $X\in\mathbb{M}_n(\mathbb{C})$.
	Hence, by \cite[Theorem 2.6]{Bresar-jmaa-2012},
	we obtain
	$
	B_TA_T=I$.
	
	\noindent\textbf{Case 2.}
	Suppose that
	$J_T(X)=A_TX^{\mathrm{T}}B_T$.
	Since
	$\sigma(X^{\mathrm{T}})=\sigma(X)$
	for every
	$X\in\mathbb{M}_n(\mathbb{C})$,
	an argument similar to that in Case 1
	also gives
	$
	B_TA_T=I$.
	
	Suppose that
	$
	J_T(X)=A_TX^{\mathrm{T}}B_T
	=A_TX^{\mathrm{T}}A_T^{-1}$. Then, for all
	$X,Y\in\mathbb{M}_n(\mathbb{C})$,
	we have
	\begin{equation}\label{comequ}
		J_T(XY)=J_T(Y)J_T(X).
	\end{equation}
	On the other hand, by (ii), we have $R_T(X)=R_{\frac{T}{2}}
	R_{\frac{T}{2}}(X)$
	and
	$
	R_T(0)=R_{\frac{T}{2}}
	R_{\frac{T}{2}}(0)$. A direct computation then gives $J_T(X)=J_{\frac{T}{2}}
	J_{\frac{T}{2}}(X)$ for every
	$X\in\mathbb{M}_n(\mathbb{C})$.
	Consequently, for all
	$X,Y\in\mathbb{M}_n(\mathbb{C})$,
	we have $\begin{aligned}
		J_T(XY)=J_{\frac{T}{2}}
		\bigl(J_{\frac{T}{2}}(XY)\bigr)=J_T(X)J_T(Y).
	\end{aligned}$
	Since $J_T$ is surjective,
	this identity and \eqref{comequ}
	imply that any two matrices commute,
	a contradiction. Therefore, in what follows, we only
	consider the case
	$
	J_T(X)=A_TXA_T^{-1}
	$
	for all
	$T,X\in\mathbb{M}_n(\mathbb{C})$. By (iv), there exists
	$S\in\mathbb{M}_n(\mathbb{C})$
	such that
	\begin{equation}\label{domain}
		\widetilde{S}=R_S(0)=2\widetilde{T}.
	\end{equation}
	We show that
	\begin{equation}\label{map1}
		J_S(X)=2J_T(X)-X
	\end{equation}
	for every
	$X\in\mathbb{M}_n(\mathbb{C})$.
	
	By (iv), it suffices to show that
	$
	J_S(\widetilde{X})
	=2J_T(\widetilde{X})-\widetilde{X}.
	$
	Indeed, by (ii), for all
	$T_1,T_2\in\mathbb{M}_n(\mathbb{C})$,
	we have
	$
	R_{T_1+T_2}(0)
	=R_{T_1}R_{T_2}(0)
	=R_{T_2}R_{T_1}(0)$. Hence
	\[
	J_{T_1}(\widetilde{T_2})
	+\widetilde{T_1}
	=
	J_{T_2}(\widetilde{T_1})
	+\widetilde{T_2}.
	\]
	By \eqref{domain}, applying the above
	identity twice yields
	\[
	\begin{aligned}
		J_S(\widetilde{X})-\widetilde{X}=J_X(\widetilde{S})-\widetilde{S}=2\bigl(
		J_X(\widetilde{T})-\widetilde{T}
		\bigr)=2\bigl(
		J_T(\widetilde{X})-\widetilde{X}
		\bigr).
	\end{aligned}
	\]
	This gives the desired conclusion.
	
	Now define the map
	$D=J_T-\operatorname{id}$
	on $\mathbb{M}_n(\mathbb{C})$.
	Let
	$X,Y\in\mathbb{M}_n(\mathbb{C})$
	be arbitrary.
	Since
	$J_T=\operatorname{id}+D$
	and
	$J_T(XY)=J_T(X)J_T(Y)$,
	we have
	$$
	XY+D(XY)
	=(X+D(X))(Y+D(Y)).$$
Expanding the above identity gives
	\begin{equation}\label{map2}
		D(XY)=D(X)Y+XD(Y)+D(X)D(Y).
	\end{equation}
	By \eqref{map1}, we have
	$J_S=\operatorname{id}+2D$.
	Since
	$J_S(XY)=J_S(X)J_S(Y)$,
	we obtain
	\[
	XY+2D(XY)
	=(X+2D(X))(Y+2D(Y)).
	\]
	Expanding this identity yields
	\begin{equation}\label{map3}
		2D(XY)=2D(X)Y+2XD(Y)+4D(X)D(Y).
	\end{equation}
	Therefore, by \eqref{map2}
	and \eqref{map3}, we have
	\[
	D(X)D(Y)=0,
	\quad
	D(XY)=D(X)Y+XD(Y).
	\]
	Since $X,Y$ are arbitrary,
	we also have
	\[
	D(X)D(YX)=0,
	\quad
	D(YX)=D(Y)X+YD(X).
	\]
	Consequently,
	\begin{equation}\label{map4}
		\begin{aligned}
			0=D(X)D(YX)=D(X)D(Y)X+D(X)YD(X)=D(X)YD(X).
		\end{aligned}
	\end{equation}
	
	We claim that $D(X)=0$
	for every
	$X\in\mathbb{M}_n(\mathbb{C})$.
	Otherwise, there exists
	$X_0\in\mathbb{M}_n(\mathbb{C})$
	such that the $(i,j)$-entry
	$a_{ij}$ of $D(X_0)$ is nonzero.
	Taking $Y=E_{ji}$, we see that
	the $(i,j)$-entry of
	$D(X_0)E_{ji}D(X_0)$
	is $a_{ij}^2\neq0$,
	contradicting \eqref{map4}.
	Thus $D=0$, and hence
	\begin{equation}\label{ide}
		J_T(X)=X
	\end{equation}
	for every
	$X\in\mathbb{M}_n(\mathbb{C})$. Finally, by (ii), for all
	$T,T_1,T_2\in\mathbb{M}_n(\mathbb{C})$,
	we have
	\[
	\begin{aligned}
		\widetilde{T_1+T_2}=R_{T_1+T_2}(0)=R_{T_1}R_{T_2}(0)=\widetilde{T_1}+\widetilde{T_2}
	\end{aligned}
	\]
	and
	$
	R_0(0)=R_{T-T}(0)=R_TR_{-T}(0).
	$
	Combining this with (i)
	and \eqref{basform}, we obtain
	$
	0=R_T(\widetilde{-T})
	=\widetilde{-T}+\widetilde{T}.
	$
	Thus
	$
	\widetilde{-T}=-\widetilde{T}$. Therefore, by \eqref{ide}
	and \eqref{basform}, we obtain
	the desired conclusion.
\end{proof}

\begin{theorem}\label{addthm}
	Let $(f,g)\in C_{\Delta}^{+}([0,\infty))$
	with $f$ nonconstant, and let
	$\alpha,\beta$ be given nonzero real numbers.
	Suppose that
	$\Phi:\mathbb{M}_n(\mathbb{C})
	\to\mathbb{M}_n(\mathbb{C})$
	is a surjective map satisfying
	\begin{equation}\label{keyequ11}
		\Phi(\Delta_{f,g}(\alpha A+\beta B))
		=\Delta_{f,g}(\alpha\Phi(A)+\beta\Phi(B))
	\end{equation}
	for all $A,B\in\mathbb{M}_n(\mathbb{C})$.
	Then:
	\begin{enumerate}[label=(\roman*)]
		\item $\Phi(0)=0$;
		
		\item $\Phi(A+B)=\Phi(A)+\Phi(B)$
		for all $A,B\in\mathbb{M}_n(\mathbb{C})$;
		
		\item There exists a surjective additive map
		$\Psi:\mathbb{M}_n(\mathbb{C})
		\to\mathbb{M}_n(\mathbb{C})$
		such that
		\[
		\Psi(\alpha T)=\alpha\Phi(T),\quad
		\Psi(\beta T)=\beta\Phi(T),\quad
		\Delta_{f,g}(\Psi(T))
		=\Phi(\Delta_{f,g}(T))
		\]
		for all $T\in\mathbb{M}_n(\mathbb{C})$.
	\end{enumerate}
\end{theorem}

\begin{proof}
(i) Since
	$\Phi:\mathbb{M}_n(\mathbb{C})
	\to\mathbb{M}_n(\mathbb{C})$
	is surjective, we can define surjective maps
	$F,G:\mathbb{M}_n(\mathbb{C})
	\to\mathbb{M}_n(\mathbb{C})$
	by $
	F(X)=\alpha\Phi\left(\frac{X}{\alpha}\right)$, $G(X)=\beta\Phi\left(\frac{X}{\beta}\right)$. Then, by \eqref{keyequ11} and Proposition \ref{basicprop1} (i), we have
	\begin{equation}\label{equivequ}
		\sigma(F(X)+G(Y))
		=\sigma(\Phi(\Delta_{f,g}(X+Y)))
	\end{equation}
	for all $X,Y\in\mathbb{M}_n(\mathbb{C})$. Suppose that
	$A,B\in\mathbb{M}_n(\mathbb{C})$
	satisfy $F(A)=F(B)$.
	By the above identity, for all
	$T,X\in\mathbb{M}_n(\mathbb{C})$, we have
	\[
	\begin{aligned}
		\sigma(F(A+T)+G(X))
		&=\sigma(\Phi(\Delta_{f,g}(A+T+X)))\\
		&=\sigma(F(A)+G(T+X))\\
		&=\sigma(F(B)+G(T+X))\\
		&=\sigma(\Phi(\Delta_{f,g}(B+T+X)))\\
		&=\sigma(F(B+T)+G(X)).
	\end{aligned}
	\]
	Since $G$ is surjective and
	$X\in\mathbb{M}_n(\mathbb{C})$ is arbitrary,
	Lemma \ref{specequ} implies that
	$F(A+T)=F(B+T)$.
	Similarly, we obtain
	$G(A+T)=G(B+T)$. Since $F$ and $G$ are surjective,
	for every $T\in\mathbb{M}_n(\mathbb{C})$,
	we can define self-maps $R_T$ and $S_T$
	on $\mathbb{M}_n(\mathbb{C})$ by
	\begin{equation}\label{promap}
		R_T(F(X))=F(X+T),\quad
		S_T(G(X))=G(X+T).
	\end{equation}
	By \eqref{equivequ}, for all
	$A,B\in\mathbb{M}_n(\mathbb{C})$, we have
	\[
	\begin{aligned}
		\sigma(R_T(F(A))+S_{-T}(G(B)))&=\sigma(F(A+T)+G(B-T))\\
		&=\sigma(\Phi(\Delta_{f,g}(A+T+B-T)))\\
		&=\sigma(\Phi(\Delta_{f,g}(A+B)))\\
		&=\sigma(F(A)+G(B)).
	\end{aligned}
	\]
	Since $F$ is surjective, there exists
	$X_0\in\mathbb{M}_n(\mathbb{C})$
	such that $F(X_0)=0$.
	Hence
	\begin{equation}\label{sum4}
		R_T(0)=R_T(F(X_0))=F(X_0+T).
	\end{equation}
	Since $F$ is surjective, we further have $\{R_T(0):T\in\mathbb{M}_n(\mathbb{C})\}
	=\mathbb{M}_n(\mathbb{C})$. Moreover, it is easy to verify that
	$R_T$ and $S_T$ satisfy conditions
	(i) and (ii) of Proposition \ref{grouppro}.
	Thus, there exists a unique
	$\widetilde{T}\in\mathbb{M}_n(\mathbb{C})$
	such that $R_T(X)=S_T(X)=X+\widetilde{T}$ for every $X\in\mathbb{M}_n(\mathbb{C})$. Combining this with \eqref{promap},
	we obtain
	\begin{equation}\label{sum5}
		F(X)+\widetilde{T}=F(X+T),\quad
		G(X)+\widetilde{T}=G(X+T).
	\end{equation}
	Taking $X=0$, we have
	\begin{equation}\label{staform}
		F(T)=\widetilde{T}+F(0),\quad
		G(T)=\widetilde{T}+G(0)
	\end{equation}
	for every $T\in\mathbb{M}_n(\mathbb{C})$.
	Therefore, by \eqref{sum5}, we obtain
	\begin{equation}\label{addprop}
		F(X)+F(T)-F(0)=F(X+T).
	\end{equation}
	Taking $X=-T$ in the above identity,
	we have
	\[
	2F(0)=F(-T)+F(T).
	\]
	Replacing $T$ by $\alpha T$ and using
	the definition of $F$, we obtain
	\begin{equation}\label{sum3}
		2\Phi(0)=\Phi(-T)+\Phi(T).
	\end{equation}
Taking $A=T/\alpha$ and $B=0$
	in \eqref{keyequ11}, we obtain
	from \eqref{staform} that
	\begin{equation}\label{sum1}
		\begin{aligned}
			\Phi(\Delta_{f,g}(T))=\Delta_{f,g}(F(T)+G(0))=\Delta_{f,g}(
			\widetilde{T}+F(0)+G(0)).
		\end{aligned}
	\end{equation}
	Taking $A=-T/\alpha$ and $B=0$
	in \eqref{keyequ11}, by Proposition \ref{basicprop1} (v), we have
	\begin{equation}\label{sum2}
		\begin{aligned}
			\Phi(\Delta_{f,g}(-T))=\Delta_{f,g}(
			\widetilde{-T}+F(0)+G(0))&=\Delta_{f,g}(
			-\widetilde{T}+F(0)+G(0))\\&=-\Delta_{f,g}(
			\widetilde{T}-F(0)-G(0)).
		\end{aligned}
	\end{equation}
	Adding \eqref{sum1} and \eqref{sum2},
	and replacing $T$ by $\Delta_{f,g}(T)$
	in \eqref{sum3}, we obtain
	\[
	\begin{aligned}
		2\Phi(0)=\Delta_{f,g}(
		\widetilde{T}+F(0)+G(0))-\Delta_{f,g}(
		\widetilde{T}-F(0)-G(0)).
	\end{aligned}
	\]
	Set
	$W=\widetilde{T}-F(0)-G(0)$.
	Then
	\[
	\begin{aligned}
		2\Phi(0)=\Delta_{f,g}(
		W+2(F(0)+G(0)))-\Delta_{f,g}(W).
	\end{aligned}
	\]
	Note that, by \eqref{sum4},
	$\widetilde{T}=F(X_0+T)$.
	Since $F$ is surjective,
	$W$ ranges over all of
	$\mathbb{M}_n(\mathbb{C})$.
	Therefore, by the above identity and Proposition \ref{basicprop3} (ii), we obtain
	$\Phi(0)=0$.

(ii) The assertion follows immediately from
	\eqref{addprop}, the definition of $F$,
	and (i).
	
(iii) Let $\Psi$ be the additive map
	on $\mathbb{M}_n(\mathbb{C})$
	defined by
	$\Psi(T)=\widetilde{T}$.
	Then, by \eqref{sum1} and (i), we have
	$
	\Delta_{f,g}(\Psi(T))
	=\Phi(\Delta_{f,g}(T))$. Since $\widetilde{T}=R_T(0)$,
	it follows from \eqref{sum4} that
	$\Psi$ is surjective. Moreover, by \eqref{staform} and (i),
	we have
	$
	\alpha\Phi(T)
	=F(\alpha T)
	=\widetilde{\alpha T}
	=\Psi(\alpha T)$
	and
	$
	\beta\Phi(T)
	=G(\beta T)
	=\widetilde{\beta T}
	=\Psi(\beta T)$.
	
	This completes the proof.
\end{proof}

Next, if we can show that the map $\Psi$
in Theorem \ref{addthm} (iii) coincides
with $\Phi$, then $\Phi$ satisfying
\eqref{eq:weighted} is a surjective
additive map. We can then apply
Theorem \ref{kkeylem} to establish
our main theorem, Theorem \ref{thm:main}. We need the following lemmas.

\begin{lemma}\label{tracepair}
	Let $\Psi,\Phi:\mathbb{M}_n(\mathbb{C})
	\to\mathbb{M}_n(\mathbb{C})$
	be additive maps. Suppose that
	\begin{equation}\label{tracepair3}
		\operatorname{tr}((\Psi(T))^2)
		=\operatorname{tr}((\Phi(T))^2)
	\end{equation}
	for every normal or square-zero matrix $T$.
	Then
	\begin{equation}\label{tp7}
		\operatorname{tr}(\Psi(X)\Psi(Y))
		=\operatorname{tr}(\Phi(X)\Phi(Y))
	\end{equation}
	for all $X,Y\in\mathbb{M}_n(\mathbb{C})$.
\end{lemma}

\begin{proof}
	Since $\Phi$ and $\Psi$ are additive,
	we can define a biadditive map
	\[
	\begin{array}{rcl}
		\Gamma:\mathbb{M}_n(\mathbb{C})
		\times\mathbb{M}_n(\mathbb{C})
		&\to&\mathbb{C}\\
		(X,Y)&\mapsto&
		\operatorname{tr}(\Psi(X)\Psi(Y))
		-\operatorname{tr}(\Phi(X)\Phi(Y)).
	\end{array}
	\]
	If $T,S\in\mathbb{M}_n(\mathbb{C})$
	are normal matrices such that
	$T+S$ is also normal, then
	\[
	\begin{aligned}
		2\Gamma(T,S)=\Gamma(T,S)+\Gamma(S,T)=\Gamma(T+S,T+S)
		-\Gamma(T,T)-\Gamma(S,S).
	\end{aligned}
	\]
	Combining this with \eqref{tracepair3},
	we obtain
	\begin{equation}\label{tp5}
		\Gamma(T,S)=0.
	\end{equation}
Therefore, for any self-adjoint matrices
	$H_1,H_2$ and any $t\in\mathbb{R}$,
	we have
	\begin{equation}\label{tracetp}
		\begin{aligned}
			0
			&=\Gamma(H_1+H_2,
			\mathrm{i}(tH_1+tH_2))\\
			&=\Gamma(H_1,\mathrm{i}tH_1)+\Gamma(H_1,\mathrm{i}tH_2)+\Gamma(H_2,\mathrm{i}tH_1)+\Gamma(H_2,\mathrm{i}tH_2)\\
			&=\Gamma(H_1,\mathrm{i}tH_2)+\Gamma(H_2,\mathrm{i}tH_1),
		\end{aligned}
	\end{equation}
	where, by \eqref{tp5}, we have
$\Gamma(H_1+H_2, \mathrm{i}(tH_1+tH_2))=\Gamma(H_1,\mathrm{i}tH_1)=\Gamma(H_2,\mathrm{i}tH_2)=0$. Similarly, we obtain $$0=\Gamma(tH_1+H_2,\mathrm{i}(tH_1+H_2))=\Gamma(tH_1,\mathrm{i}H_2)+\Gamma(H_2,\mathrm{i}tH_1).$$
	Thus, by the above identity and
	\eqref{tracetp}, we have
	\begin{equation}\label{comreal}
		\Gamma(tH_1,\mathrm{i}H_2)
		=\Gamma(H_1,\mathrm{i}tH_2)
	\end{equation}
	for every $t\in\mathbb{R}$. Now let
	$X,Y\in\mathbb{M}_n(\mathbb{C})$
	be arbitrary. Write
	$X=X_1+\mathrm{i}X_2$ and $Y=Y_1+\mathrm{i}Y_2$, where $X_i,Y_i$ are self-adjoint
	matrices for $i\in\{1,2\}$. A direct computation gives
	\begin{equation}\label{tracepair4}
		\begin{aligned}
			\Gamma(X,Y)=\Gamma(X_1,Y_1)
			+\Gamma(X_1,\mathrm{i}Y_2)+\Gamma(\mathrm{i}X_2,Y_1)
			+\Gamma(\mathrm{i}X_2,\mathrm{i}Y_2).
		\end{aligned}
	\end{equation}
	Moreover, by \eqref{tp5}, we have
	\begin{equation}\label{tp6}
		\Gamma(X_1,Y_1)=0,\quad
		\Gamma(\mathrm{i}X_2,\mathrm{i}Y_2)=0.
	\end{equation}
If we can show that
	\begin{equation}\label{keyine}
		\Gamma(H_1,\mathrm{i}H_2)=0
	\end{equation}
	for any self-adjoint matrices $H_1,H_2$,
	then, since
	$\Gamma(X,Y)=\Gamma(Y,X)$,
	it follows from \eqref{tracepair4}
	and \eqref{tp6} that
	$\Gamma(X,Y)=0$, which is precisely
	\eqref{tp7}. Furthermore, since $H_1$ is self-adjoint,
	there exists a unitary matrix $U$
	such that
	$
	U^*H_1U
	=\operatorname{diag}
	(\lambda_1,\lambda_2,\ldots,\lambda_n)$, where $\lambda_i\in\mathbb{R}$
	for $i\in\{1,2,\ldots,n\}$.
	Write $H_2=\widetilde{H}
	+\sum_{1\leq i<j\leq n}K_{ij}$, where $\widetilde{H}
	=U\operatorname{diag}
	(z_{11},z_{22},\ldots,z_{nn})U^*$
	and
	$
	K_{ij}
	=U(z_{ij}E_{ij}
	+\overline{z_{ij}}E_{ji})U^*$. Since
	$H_1\widetilde{H}
	=\widetilde{H}H_1$,
	the matrix
	$H_1+\mathrm{i}\widetilde{H}$
	is normal. Hence, by \eqref{tp5},
	we have
	\[
	\Gamma(H_1,\mathrm{i}\widetilde{H})=0.
	\]
	Therefore, to obtain \eqref{keyine},
	by the above identity and the
	biadditivity of $\Gamma$, it suffices
	to show that
	\begin{equation}\label{keyineq}
		\Gamma(H_1,\mathrm{i}K_{ij})=0
	\end{equation}
	for every $1\leq i<j\leq n$. Indeed, if $z_{ij}=0$, then
	$K_{ij}=0$, and the above identity
	clearly holds. Henceforth, assume
	that $z_{ij}\neq0$. Take the permutation matrix
	$
	P_{ij}
	=(e_i,e_j,e_{k_1},e_{k_2},
	\ldots,e_{k_{n-2}})$,
	where
	$
	k_l\in\{1,2,\ldots,n\}
	\setminus\{i,j\}$, $l\in\{1,2,\ldots,n-2\}$, satisfy
	$k_1<k_2<\cdots<k_{n-2}$.
	
	Set
	$
	V=UP_{ij}$ and $Q_{ij}=U(E_{ii}-E_{jj})U^*$, where $V$ is unitary, and let
	$0_{n-2}$ denote the zero matrix
	of order $n-2$.
	Then it is easy to verify that
	$
	V^*H_1V
	=\operatorname{diag}
	(\lambda_i,\lambda_j,
	\lambda_{k_1},\ldots,\lambda_{k_{n-2}})$,
	
	\[
	V^*Q_{ij}V
	=
	\begin{pmatrix}
		1&0&0\\
		0&-1&0\\
		0&0&0_{n-2}
	\end{pmatrix},\quad 
	V^*K_{ij}V
	=
	\begin{pmatrix}
		0&z_{ij}&0\\
		\overline{z_{ij}}&0&0\\
		0&0&0_{n-2}
	\end{pmatrix}.
	\]
Set
	$
	R=|z_{ij}|Q_{ij}
	+\mathrm{i}K_{ij}.$
	A direct computation gives
	$(V^*RV)^2=0$, and hence $R^2=0$. Thus, for any $t\in\mathbb{R}$,
	since $tQ_{ij}$ and $tK_{ij}$
	are normal matrices and $(tR)^2=0$,
	a direct computation using
	\eqref{tracepair3} gives
	\[
	\begin{aligned}
		0=\Gamma(tR,tR)&=\Gamma(t|z_{ij}|Q_{ij},
		t|z_{ij}|Q_{ij})+
		2\Gamma(t|z_{ij}|Q_{ij},
		\mathrm{i}tK_{ij})+
		\Gamma(\mathrm{i}tK_{ij},
		\mathrm{i}tK_{ij})\\&=2\Gamma(t|z_{ij}|Q_{ij},
		\mathrm{i}tK_{ij})=2\Gamma(t^2|z_{ij}|Q_{ij},
		\mathrm{i}K_{ij}),
	\end{aligned}
	\]
	where the last equality follows
	from \eqref{comreal}. Since $|z_{ij}|>0$, the above identity
	implies that
	\[
	\Gamma(sQ_{ij},\mathrm{i}K_{ij})=0
	\]
	for every $s\geq0$.
	Moreover, since $\Gamma$ is biadditive,
	we obtain
	\begin{equation}\label{comreal1}
		\Gamma(rQ_{ij},\mathrm{i}K_{ij})=0
	\end{equation}
	for every $r\in\mathbb{R}$. Now take
	$
	r=\frac{\lambda_i-\lambda_j}{2}$. A direct computation gives
	\[
	\begin{aligned}
		V^*(H_1-rQ_{ij})V=\operatorname{diag}\left(
		\frac{\lambda_i+\lambda_j}{2},
		\frac{\lambda_i+\lambda_j}{2},
		\lambda_{k_1},\lambda_{k_2},
		\ldots,\lambda_{k_{n-2}}
		\right).
	\end{aligned}
	\]
	Consequently,
	$
	\begin{aligned}
		V^*(H_1-rQ_{ij})VV^*K_{ij}V=V^*K_{ij}VV^*(H_1-rQ_{ij})V,
	\end{aligned}$ that is,
	\[
	(H_1-rQ_{ij})K_{ij}
	=K_{ij}(H_1-rQ_{ij}).
	\]Since $H_1-rQ_{ij}$ and
	$K_{ij}$ are normal matrices,
	the above identity implies that
	$H_1-rQ_{ij}+\mathrm{i}K_{ij}$
	is normal.
	Hence, by \eqref{tp5}, we have $\Gamma(H_1-rQ_{ij},
	\mathrm{i}K_{ij})=0$. Combining this with \eqref{comreal1}
	and the biadditivity of $\Gamma$,
	we obtain
	\[
	\begin{aligned}
		\Gamma(H_1,\mathrm{i}K_{ij})
		&=\Gamma(H_1-rQ_{ij},
		\mathrm{i}K_{ij})+
		\Gamma(rQ_{ij},\mathrm{i}K_{ij})=0.
	\end{aligned}
	\]
	Thus \eqref{keyineq} holds.
	This completes the proof.
\end{proof}

\begin{lemma}\label{mapcom}
	Let $(f,g)\in C_{\Delta}^{+}([0,\infty))$
	with $f(0)=0$, and fix a nonzero real
	number $\alpha$.
	Suppose that
	$\Psi,\Phi:\mathbb{M}_n(\mathbb{C})
	\to\mathbb{M}_n(\mathbb{C})$
	are surjective additive maps satisfying
	$
	\Delta_{f,g}(\Psi(T))
	=\Phi(\Delta_{f,g}(T))
	$
	and
	$
	\Psi(\alpha T)=\alpha\Phi(T)
	$
	for every $T\in\mathbb{M}_n(\mathbb{C})$.
	Then there exists a complex-linear bijection
	$\widetilde{\Phi}:\mathbb{M}_n(\mathbb{C})
	\to\mathbb{M}_n(\mathbb{C})$
	such that
	$\Phi=\widetilde{\Phi}\circ\Psi$
	and $|\det\widetilde{\Phi}|=1$.
	Moreover, for every normal matrix
	$N\in\mathbb{M}_n(\mathbb{C})$,
	every $t>0$, and every $k\in\mathbb{N}\cup\{0\}$,
	we have
	\begin{equation}\label{fgtran}
		\Delta_{f,g}
		(t\widetilde{\Phi}^{k}(\Psi(N)))
		=t\widetilde{\Phi}^{k+1}(\Psi(N)).
	\end{equation}
\end{lemma}

\begin{proof}
	Let $N\in\mathbb{M}_n(\mathbb{C})$
	be an arbitrary normal matrix.
	By Proposition \ref{basicprop1} (ii), Proposition \ref{basicprop2} (i)
	and the hypotheses, we have
	\begin{equation}\label{trace11}
		\begin{aligned}
			\operatorname{tr}((\Psi(N))^2)=\operatorname{tr}
			((\Delta_{f,g}(\Psi(N)))^2)=\operatorname{tr}
			((\Phi(\Delta_{f,g}(N)))^2)=\operatorname{tr}((\Phi(N))^2).
		\end{aligned}
	\end{equation}
Let $S\in\mathbb{M}_n(\mathbb{C})$
	be arbitrary with $S^2=0$.
	By the hypotheses and
	Proposition \ref{basicprop2} (v), we have
	$
	\Delta_{f,g}(\Psi(S))
	=\Phi(\Delta_{f,g}(S))=0$, and hence $\Psi(S)^2=0$.
	Since $(\alpha S)^2=0$,
	replacing $S$ by $\alpha S$
	in the preceding argument immediately
	gives $\Psi(\alpha S)^2=0$.
	Combining this with the hypothesis
	$\Psi(\alpha S)=\alpha\Phi(S)$,
	we obtain $\Phi(S)^2=0$.
	Thus
	\begin{equation}\label{trace12}
		\operatorname{tr}((\Psi(S))^2)
		=\operatorname{tr}((\Phi(S))^2)=0.
	\end{equation}
Therefore, by Lemma \ref{tracepair},
	\eqref{trace11}, and \eqref{trace12},
	we have
	\begin{equation}\label{tracepair1}
		\operatorname{tr}(\Psi(X)\Psi(Y))
		=\operatorname{tr}(\Phi(X)\Phi(Y))
	\end{equation}
	for all $X,Y\in\mathbb{M}_n(\mathbb{C})$. Take $Z\in\ker\Psi$.
	By the above identity, we have
	\[
	\operatorname{tr}(\Phi(Z)\Phi(Y))=0
	\]
	for every $Y\in\mathbb{M}_n(\mathbb{C})$.
	Since $\Phi$ is surjective,
	there exists $\widetilde{Z}$ such that
	$\Phi(\widetilde{Z})=(\Phi(Z))^*$.
	Hence
	\[
	\begin{aligned}
		\|\Phi(Z)\|_2=\sqrt{
			\operatorname{tr}
			(\Phi(Z)(\Phi(Z))^*)}=\sqrt{
			\operatorname{tr}
			(\Phi(Z)\Phi(\widetilde{Z}))}=0.
	\end{aligned}
	\]
	Thus $\Phi(Z)=0$, that is,
	$\ker\Psi\subseteq\ker\Phi$.
	Since $\Psi$ and $\Phi$ are
	surjective additive maps,
	we can therefore define a
	surjective additive map
	$\widetilde{\Phi}$ on
	$\mathbb{M}_n(\mathbb{C})$ by
	\begin{equation}\label{defmap}
		\widetilde{\Phi}(\Psi(T))=\Phi(T)
	\end{equation}
	for every $T\in\mathbb{M}_n(\mathbb{C})$. Since $\Psi$ is surjective,
	for any $X,Y\in\mathbb{M}_n(\mathbb{C})$,
	there exist $\widetilde{X},\widetilde{Y}$
	such that
	$\Psi(\widetilde{X})=X$ and $\Psi(\widetilde{Y})=Y$. By \eqref{tracepair1} and
	\eqref{defmap}, we obtain
	\begin{equation}\label{tracepair2}
		\operatorname{tr}(XY)
		=\operatorname{tr}
		(\widetilde{\Phi}(X)\widetilde{\Phi}(Y)).
	\end{equation}
Let $\lambda\in\mathbb{C}$ be arbitrary.
	By the above identity, we have
	\[
	\begin{aligned}
		\operatorname{tr}\Big(
		(\widetilde{\Phi}(\lambda X)
		-\lambda\widetilde{\Phi}(X))
		\widetilde{\Phi}(Y)\Big)=\operatorname{tr}
		(\widetilde{\Phi}(\lambda X)
		\widetilde{\Phi}(Y))-\lambda
		\operatorname{tr}
		(\widetilde{\Phi}(X)
		\widetilde{\Phi}(Y))=\operatorname{tr}(\lambda XY)
		-\lambda\operatorname{tr}(XY)=0.
	\end{aligned}
	\]
	Since $Y$ is arbitrary,
	the above identity implies that $\widetilde{\Phi}(\lambda X)
	=\lambda\widetilde{\Phi}(X)$. Moreover, since $\Phi$ is
	surjective and additive,
	$\widetilde{\Phi}$ is a
	complex-linear bijection. Furthermore, consider the
	nondegenerate bilinear form
	$
	B(X,Y)=\operatorname{tr}(XY)
	$
	on $\mathbb{M}_n(\mathbb{C})$.
	Choose a basis $\mathscr{E}$ of
	$\mathbb{M}_n(\mathbb{C})$.
	Then $B(X,Y)=[X]^{\mathrm{T}}G[Y]$, where $[X]$ and $[Y]$ are the
	coordinate column vectors of
	$X$ and $Y$ with respect to
	$\mathscr{E}$, $[X]^{\mathrm{T}}$ denotes the transpose
	of $[X]$, and $G$ is the
	invertible matrix representing
	$B$ with respect to $\mathscr{E}$. Let $A$ be the matrix representing
	$\widetilde{\Phi}$ with respect
	to $\mathscr{E}$.
	Then
	$
	\begin{aligned}
		B(\widetilde{\Phi}(X),
		\widetilde{\Phi}(Y))
		=[X]^{\mathrm{T}}A^{\mathrm{tr}}GA[Y].
	\end{aligned}
$	
	Hence, by \eqref{tracepair2},
	we have $A^{\mathrm{tr}}GA=G$. Taking determinants on both sides
	gives $|\det A|=1$.
	Therefore,
	$|\det\widetilde{\Phi}|=1$.
	
	Finally, let
	$N\in\mathbb{M}_n(\mathbb{C})$
	be a normal matrix.
	By the additivity of $\Psi$,
	Proposition \ref{basicprop2} (i),
	and the hypotheses, we have
	$
	\begin{aligned}
		\Delta_{f,g}(r\Psi(N))=\Delta_{f,g}(\Psi(rN))=\Phi(\Delta_{f,g}(rN))=\Phi(rN)=r\Phi(N)
	\end{aligned}
	$
	for every $r\in\mathbb{Q}$.
	Thus, by Proposition \ref{basicprop1} (vi),
	we obtain
	\begin{equation}\label{raequ}
		\Delta_{f,g}(r\Psi(N))
		=r\Phi(N)
	\end{equation}
	for every $r\in\mathbb{R}$. Taking $r=t>0$ in \eqref{raequ}
	and using \eqref{defmap}, we obtain
	\[
	\Delta_{f,g}(t\Psi(N))
	=t\Phi(N)
	=t\widetilde{\Phi}(\Psi(N)),
	\]
	which proves \eqref{fgtran} for $k=0$. Moreover, by the hypotheses,
	\eqref{defmap}, and the linearity
	of $\widetilde{\Phi}$, it is easy
	to show by induction that
	\[
	\Psi(\alpha^kN)
	=\alpha^k\widetilde{\Phi}^k(\Psi(N))
	\]
	for every $k\in\mathbb{N}$. Take $t>0$.
	By the above identity,
	\eqref{raequ}, and \eqref{defmap},
	we have
	\[
	\begin{aligned}
		\Delta_{f,g}
		(t\widetilde{\Phi}^k(\Psi(N)))=\Delta_{f,g}
		\left(\frac{t}{\alpha^k}
		\Psi(\alpha^kN)\right)&=\frac{t}{\alpha^k}\Phi(\alpha^kN)=\frac{t}{\alpha^k}
		\widetilde{\Phi}(\Psi(\alpha^kN))\\&=\frac{t}{\alpha^k}
		\widetilde{\Phi}
		\left(\alpha^k
		\widetilde{\Phi}^k(\Psi(N))\right)=t\widetilde{\Phi}^{k+1}(\Psi(N)),
	\end{aligned}
	\]
	where the last equality follows from
	the linearity of $\widetilde{\Phi}$.
	This completes the proof.
\end{proof}

\begin{lemma}\label{inequeq}
	Let $\widetilde{\Phi}:\mathbb{M}_n(\mathbb{C})
	\to\mathbb{M}_n(\mathbb{C})$
	be a complex-linear bijection satisfying
	$|\det\widetilde{\Phi}|=1$.
	Let $\Psi$ be a surjective additive map
	on $\mathbb{M}_n(\mathbb{C})$.
	Suppose that
	$
	\|\widetilde{\Phi}^{k+1}(\Psi(N))\|_F
	\leq
	\|\widetilde{\Phi}^{k}(\Psi(N))\|_F
	$
	for every normal matrix
	$N\in\mathbb{M}_n(\mathbb{C})$
	and every $k\in\mathbb{N}\cup\{0\}$.
	Then
	$
	\|\widetilde{\Phi}^{k+1}(\Psi(N))\|_F
	=
	\|\widetilde{\Phi}^{k}(\Psi(N))\|_F
	$.
\end{lemma}

\begin{proof}
	Let $X\in\mathbb{M}_n(\mathbb{C})$
	be arbitrary.
	Since $\Psi$ is surjective,
	there exist normal matrices
	$Y_1,Y_2\in\mathbb{M}_n(\mathbb{C})$
	such that $\Psi(Y_1+Y_2)=X$.
	Then, by the hypothesis, we immediately
	obtain
	\begin{align*}
		\|\widetilde{\Phi}^k(X)\|_F=\|\widetilde{\Phi}^k
		(\Psi(Y_1+Y_2))\|_F&\leq
		\|\widetilde{\Phi}^k(\Psi(Y_1))\|_F
		+\|\widetilde{\Phi}^k(\Psi(Y_2))\|_F\\
		&\leq
		\|\Psi(Y_1)\|_F+\|\Psi(Y_2)\|_F
	\end{align*}
	for every $k\in\mathbb{N}\cup\{0\}$.
	Hence, by the uniform boundedness principle,
	we have
	$$
	\sup_{k\in\mathbb{N}\cup\{0\}}
	\|\widetilde{\Phi}^k\|<\infty,$$
	where
	$
	\|\widetilde{\Phi}^k\|
	=\sup\limits_{\|X\|_F=1}
	\|\widetilde{\Phi}^k(X)\|_F$. Moreover, the set
	\[
	\mathscr{S}
	=\left\{
	\Phi:\mathbb{M}_n(\mathbb{C})
	\to\mathbb{M}_n(\mathbb{C})
	\;\middle|\;
	\Phi\text{ is complex-linear and }
	\|\Phi\|\leq
	\sup_{k\in\mathbb{N}\cup\{0\}}
	\|\widetilde{\Phi}^k\|
	\right\}
	\]
	is compact.
	Thus there exist a strictly increasing
	sequence $\{k_j\}$ in $\mathbb{N}$
	and $\Phi_0\in\mathscr{S}$ such that
	$
	\lim_{j\to\infty}
	\widetilde{\Phi}^{k_j}=\Phi_0$. Set $m_j=k_{j+1}-k_j$.
	Since
	$
	|\det\widetilde{\Phi}^{k_j}|
	=|\det\widetilde{\Phi}|^{k_j}=1$,
	by the continuity of the determinant
	function, we have
	$|\det\Phi_0|=1$.
	Hence $\Phi_0$ is invertible, and
	\[
	\begin{aligned}
		\lim_{j\to\infty}
		\widetilde{\Phi}^{m_j}=
		\lim_{j\to\infty}
		\widetilde{\Phi}^{k_{j+1}}
		\lim_{j\to\infty}
		\widetilde{\Phi}^{-k_j}=\Phi_0\Phi_0^{-1}
		=\operatorname{id}.
	\end{aligned}
	\]
	Therefore, by the above identity,
	for every $k\in\mathbb{N}\cup\{0\}$
	and every normal matrix
	$N\in\mathbb{M}_n(\mathbb{C})$,
	we have
	\begin{equation}\label{limequ}
		\lim_{j\to\infty}
		\widetilde{\Phi}^{k+m_j}(\Psi(N))
		=
		\widetilde{\Phi}^{k}(\Psi(N)).
	\end{equation}
Moreover, since $m_j\geq1$,
	by the hypothesis, we have
	\[
	\begin{aligned}
		\|\widetilde{\Phi}^{k+m_j}(\Psi(N))\|_F\leq
		\|\widetilde{\Phi}^{k+1}(\Psi(N))\|_F\leq
		\|\widetilde{\Phi}^{k}(\Psi(N))\|_F.
	\end{aligned}
	\]
	Letting $j\to\infty$ and using
	\eqref{limequ}, we obtain
	$
	\|\widetilde{\Phi}^{k+1}(\Psi(N))\|_F
	=
	\|\widetilde{\Phi}^{k}(\Psi(N))\|_F$.
	
	This completes the proof.
\end{proof}

We now prove that the maps $\Psi$ and $\Phi$
in Theorem \ref{addthm} (iii) coincide.

\begin{proposition}\label{mapequ}
	Let $(f,g)\in C_{\Delta}^{+}([0,\infty))$
	with $f$ nonconstant, and fix nonzero
	real numbers $\alpha,\beta$.
	Suppose that
	$\Psi,\Phi:\mathbb{M}_n(\mathbb{C})
	\to\mathbb{M}_n(\mathbb{C})$
	are surjective additive maps satisfying
	\begin{equation}\label{equ1}
		\begin{cases}
			\Delta_{f,g}(\Psi(T))
			=\Phi(\Delta_{f,g}(T)),\\
			\Psi(\alpha T)=\alpha\Phi(T),\\
			\Psi(\beta T)=\beta\Phi(T)
		\end{cases}
	\end{equation}
	for every $T\in\mathbb{M}_n(\mathbb{C})$.
	Then $\Psi=\Phi$.
\end{proposition}

\begin{proof}
	Since, for every
	$X\in\mathbb{M}_n(\mathbb{C})$, we have
	$
	X=\frac{X+X^*}{2}
	+\mathrm{i}\frac{X-X^*}{2\mathrm{i}}$,
	where
	$\frac{X+X^*}{2}$ and
	$\mathrm{i}\frac{X-X^*}{2\mathrm{i}}$
	are normal matrices, by the additivity
	of $\Psi$ and $\Phi$, to prove that
	$\Psi(X)=\Phi(X)$ for every $X$,
	it suffices to show that
	$\Phi(N)=\Psi(N)$ for every normal matrix $N$.
	
	\noindent\textbf{Case 1:} $f(0)>0$. Let $N$ be an arbitrary normal matrix,
	and set $Y=\Psi(N)$.
	For every positive rational number $q$,
	by the additivity of $\Psi$ and $\Phi$
	and \eqref{equ1}, we have
	\begin{equation}\label{equ11}
		\begin{aligned}
			\Delta_{f,g}(qY)=\Delta_{f,g}(\Psi(qN))=\Phi(\Delta_{f,g}(qN))=\Phi(qN)=q\Phi(N).
		\end{aligned}
	\end{equation}
	Moreover, since $f(0)>0$,
	by Proposition \ref{basicprop1} (iv), we have
	$$
	\frac{\Delta_{f,g}(qY)}{q}
	=f(q|Y|)Y(f(q|Y|))^{-1}.
	$$
	By the continuous functional calculus,
	we have
	$
	\lim_{q\to0^+}f(q|Y|)=f(0)I,
	$
	and consequently
	$
	\lim_{q\to0^+}
	(f(q|Y|))^{-1}=(f(0))^{-1}I$. Hence, by the above identity, we obtain
	$
	\lim_{q\to0^+}
	\frac{\Delta_{f,g}(qY)}{q}=Y$. Combining this with \eqref{equ11},
	we obtain $\Phi(N)=\Psi(N)$.

	\noindent\textbf{Case 2:}
	$f(0)=0$ and $g$ is nonconstant. Let $N$ be an arbitrary normal matrix.
	By Lemma \ref{mapcom}, taking $t=1$
	in \eqref{fgtran}, we have
	\begin{align}\label{tequ1}
		\widetilde{\Phi}^{k+1}(\Psi(N))
		=\Delta_{f,g}
		(\widetilde{\Phi}^k(\Psi(N))).
	\end{align}
	Moreover, for every $t>0$,
	combining the above identity with
	\eqref{fgtran}, we obtain
	$$
	\begin{aligned}
		\Delta_{f,g}
		(t\widetilde{\Phi}^k(\Psi(N)))=t\widetilde{\Phi}^{k+1}(\Psi(N))=t\Delta_{f,g}
		(\widetilde{\Phi}^k(\Psi(N))).
	\end{aligned}
	$$
	Thus, by the above identity
	and Proposition \ref{basicprop2} (vi), we have
	\[
	\begin{aligned}
		\|\widetilde{\Phi}^{k+1}(\Psi(N))\|_F
		&=\|\Delta_{f,g}
		(\widetilde{\Phi}^k(\Psi(N)))\|_F\leq
		\|\widetilde{\Phi}^k(\Psi(N))\|_F.
	\end{aligned}
	\]
	Hence, by Lemma \ref{inequeq}, we obtain
	\[
	\|\widetilde{\Phi}^{k+1}(\Psi(N))\|_F
	=
	\|\widetilde{\Phi}^k(\Psi(N))\|_F.
	\]
	Taking $k=0$ in the above identity
	and using \eqref{tequ1}, we have
	$$
	\|\Delta_{f,g}(\Psi(N))\|_F=\|\widetilde{\Phi}(\Psi(N))\|_F=\|\Psi(N)\|_F.
	$$
	Then, by Proposition \ref{basicprop2} (vi), we obtain
	$
	\Delta_{f,g}(\Psi(N))=\Psi(N)$. Moreover, by \eqref{tequ1}
	and $\Phi=\widetilde{\Phi}\circ\Psi$,
	we have
	$
	\Delta_{f,g}(\Psi(N))
	=\widetilde{\Phi}(\Psi(N))
	=\Phi(N)$. Therefore, $\Phi(N)=\Psi(N)$.

	\noindent\textbf{Case 3:}
	$f(0)=0$ and $g$ is constant. Let $N$ be an arbitrary normal matrix.
	By Lemma \ref{mapcom}, taking $t=1$
	in \eqref{fgtran}, we have
	\begin{align}\label{tequ2}
		\widetilde{\Phi}^{k+1}(\Psi(N))
		=\Delta_D(\widetilde{\Phi}^k(\Psi(N))).
	\end{align}
	Moreover, by Proposition \ref{basicprop1} (vii),
	we have
	$
	\|\Delta_D
	(\widetilde{\Phi}^k(\Psi(N)))\|_F
	\leq
	\|\widetilde{\Phi}^k(\Psi(N))\|_F$. Hence
	$$
	\|\widetilde{\Phi}^{k+1}(\Psi(N))\|_F
	\leq
	\|\widetilde{\Phi}^k(\Psi(N))\|_F.$$ Therefore, taking $k=0$, by \eqref{tequ2} and Lemma \ref{inequeq}, we obtain
	\[
	\begin{aligned}
		\|\Delta_D(\Psi(N))\|_F=\|\widetilde{\Phi}(\Psi(N))\|_F=\|\Psi(N)\|_F.
	\end{aligned}
	\]
Moreover, by \eqref{equ1},
	\eqref{tequ2}, and Proposition \ref{basicprop1} (vii), we have
	\[
	\begin{aligned}
		\Phi(\alpha N)=\Phi(\Delta_D(\alpha N))=\Delta_D(\Psi(\alpha N))=\Delta_D(\alpha\Phi(N))=\alpha\Delta_D(\Phi(N))=\alpha\Delta_D(\Delta_D(\Psi(N))),
	\end{aligned}
	\]
	where the first equality holds because
	$\Delta_D(T)=T$ for every normal matrix $T$. Set $Y=\Psi(N)$. For any $q\in\mathbb{Q}$, we have
	\[
	\begin{aligned}
		\Delta_D(Y+q\alpha\Delta_D(Y))&=\Delta_D(\Psi(N)+q\alpha\Delta_D(\Psi(N)))=\Delta_D(\Psi(N)+q\alpha\Phi(N))\\&=\Delta_D(\Psi(N+q\alpha N))=\Phi(\Delta_D(N+q\alpha N))\\
		&=\Phi(N+q\alpha N)=\Phi(N)+q\Phi(\alpha N)\\
		&=\Delta_D(Y)
		+q\alpha\Delta_D(\Delta_D(Y)).
	\end{aligned}
	\]
	Hence, by Proposition \ref{basicprop2} (ii),
	we obtain $\Delta_D(Y)=Y$, that is,
	$
	\Delta_D(\Psi(N))=\Psi(N)$. Combining this with \eqref{tequ2}
	for $k=0$, we have $\Phi(N)=\widetilde{\Phi}(\Psi(N))=\Delta_D(\Psi(N))$. Thus $\Phi(N)=\Psi(N)$.
	
	This completes the proof.
\end{proof}
	
We now prove our main result.

\begin{theorem}\label{main}
	Let $n\geq2$, $(f,g)\in C_{\Delta}^{+}([0,\infty))$,
	and $\alpha,\beta\in\mathbb{R}\setminus\{0\}$.
	A surjective map
	$\Phi:\mathbb{M}_n(\mathbb{C})\to\mathbb{M}_n(\mathbb{C})$
	satisfies
	\begin{equation}\label{eq:weighted1}
		\Phi(\Delta_{f,g}(\alpha A+\beta B))
		=\Delta_{f,g}(\alpha\Phi(A)+\beta\Phi(B))
	\end{equation}
	for all $A,B\in\mathbb{M}_n(\mathbb{C})$
	if and only if one of the following holds:
	\begin{enumerate}[label=(\roman*)]
		\item $f$ is nonconstant, and there exist
		a unitary matrix $U$ and a nonzero scalar
		$c\in\mathbb{C}$ with $|c|\in G_{f,g}$
		such that either
		\[
		\Phi(A)=cUAU^*
		\quad\text{or}\quad
		\Phi(A)=cU\overline{A}U^*
		\]
		for all $A\in\mathbb{M}_n(\mathbb{C})$.
		Moreover, $G_{f,g}$ is one of
		\[
		\{1\},\quad
		q^{\mathbb{Z}}\text{ for some }q>1,
		\quad\text{or}\quad
		(0,\infty).
		\]
		
		\item $f$ is constant, and there exist
		a surjective $\mathbb{Q}(\alpha,\beta)$-linear
		map $S$ on $\mathbb{M}_n(\mathbb{C})$
		and a matrix $D\in\mathbb{M}_n(\mathbb{C})$
		with $(\alpha+\beta-1)D=0$ such that
		\[
		\Phi(A)=S(A)+D
		\qquad(A\in\mathbb{M}_n(\mathbb{C})).
		\]
	\end{enumerate}
\end{theorem}
\begin{proof}
If $\Phi$ satisfies \eqref{eq:weighted1}
and $f$ is nonconstant, then
Theorem \ref{addthm} and
Proposition \ref{mapequ} imply that
$\Phi$ is a surjective additive map
satisfying
$\Phi\circ\Delta_{f,g}
=\Delta_{f,g}\circ\Phi$.
Finally, by Theorem \ref{kkeylem}
and Proposition \ref{basicprop3} (i),
we obtain (i).

If $\Phi$ satisfies \eqref{eq:weighted1}
and $f$ is constant, then by
Proposition \ref{basicprop2} (i)
and \eqref{eq:weighted1}, we have
\[
\Phi(\alpha A+\beta B)
=\alpha\Phi(A)+\beta\Phi(B)
\]
for all $A,B\in\mathbb{M}_n(\mathbb{C})$.
Taking $A=B=0$, we obtain
$(\alpha+\beta-1)\Phi(0)=0$.
Set $D=\Phi(0)$ and define
$S(X)=\Phi(X)-D$ for
$X\in\mathbb{M}_n(\mathbb{C})$.
Then $S$ is surjective and satisfies
\[
S(\alpha X+\beta Y)
=\alpha S(X)+\beta S(Y)
\]
for all $X,Y\in\mathbb{M}_n(\mathbb{C})$.
In particular,
\begin{equation}\label{ablinear}
S(\alpha X)=\alpha S(X),\quad S(\beta Y)=\beta S(Y)
\end{equation}
Since $\alpha,\beta\neq0$,
for any $X,Y\in\mathbb{M}_n(\mathbb{C})$,
we have
\[
\begin{aligned}
	S(X+Y)=S\left(
	\alpha\frac{X}{\alpha}
	+\beta\frac{Y}{\beta}
	\right)=\alpha S\left(\frac{X}{\alpha}\right)
	+\beta S\left(\frac{Y}{\beta}\right)=S(X)+S(Y).
\end{aligned}
\]
Thus $S$ is additive. Moreover, it follows from \eqref{ablinear} that $S$ is
$\mathbb{Q}(\alpha,\beta)$-linear.
Therefore, (ii) holds.

Conversely, it is straightforward to verify
that the maps in (i) and (ii) satisfy
\eqref{eq:weighted1}.

This completes the proof.
\end{proof}

We conclude the paper with the following
example, which shows that all three
possibilities for $G_{f,g}$ in (i) of
the preceding theorem do occur, and
that the maps in (ii) need not be
real-linear or multiplicative.

\begin{example}\label{ex}
	\begin{enumerate}[label=(\roman*)]
		\item 
		Let $f(x)=1+x$ and $g(x)=\frac{x}{1+x}$
		for $x\geq0$. Then $G_{f,g}=\{1\}$.
		Indeed, by Proposition \ref{basicprop3} (i),
		for any $r\in G_{f,g}$, there exists
		$c>0$ such that
		$1+rx=c(1+x)$ for $x>0$, which implies $r=c=1$.
		
		\item 
		Let $q>1$. Define
		\[
		f(x)=
		\begin{cases}
			\displaystyle
			q^{k/2}+
			\frac{x-q^k}{q^{k/2}(1+\sqrt q)},
			& q^k\leq x\leq q^{k+1},
			\quad k\in\mathbb Z,\\[6pt]
			0, & x=0,
		\end{cases}
		\quad \text{and}\quad
		g(x)=
		\begin{cases}
			\displaystyle\frac{x}{f(x)}, & x>0,\\[6pt]
			0, & x=0.
		\end{cases}
		\]
		Then $G_{f,g}=q^{\mathbb Z}$.
		Indeed, it is easy to see that
		$f(x)\leq\sqrt{x}$ for $x>0$,
		with equality if and only if
		$x\in q^{\mathbb Z}$.
		Clearly, $f(q^k x)=q^{k/2}f(x)$
		for $x>0$, $k\in\mathbb{Z}$,
		so $q^{\mathbb Z}\subseteq G_{f,g}$.
		Conversely, let $r\in G_{f,g}$.
		Then there exists $c>0$ such that
		$f(rx)=cf(x)$ for $x>0$.
		Since $f(1)=1$, we have
		$f(r)=c$ and $f(1/r)=1/c$.
		Consequently,
		$c=f(r)\leq\sqrt r$ and
		$\frac1c=f(1/r)\leq\frac1{\sqrt r}$.
		Thus $c=\sqrt r$, and hence
		$f(r)=\sqrt r$.
		It follows that $r\in q^{\mathbb Z}$.
		Therefore, $G_{f,g}=q^{\mathbb Z}$.
		
		\item Let $\lambda\in(0,1)$, $f(x)=x^\lambda$ and $g(x)=x^{1-\lambda}$ for $x\ge0$. Let $r>0$. Since $f(rx)=r^\lambda f(x)$, by Proposition \ref{basicprop3} (i),
		we have $G_{f,g}=(0,\infty)$.
		
		\item 
		Let $f(x)=1$ and $g(x)=x$
		for $x\geq0$.
		If $\alpha=\beta=\frac12$,
		the map $\Phi(A)=A+I$ for
		$A\in\mathbb{M}_n(\mathbb{C})$
		satisfies \eqref{eq:weighted1},
		but is not additive, real-linear,
		or multiplicative.
		
		More generally, for arbitrary
		$\alpha,\beta\in\mathbb{R}\setminus\{0\}$,
		since $\mathbb{Q}(\alpha,\beta)$
		is countable, choose
		$\xi\in\mathbb{R}\setminus
		\mathbb{Q}(\alpha,\beta)$
		and a Hamel basis $\mathcal{B}$
		of $\mathbb{R}$ over
		$\mathbb{Q}(\alpha,\beta)$
		containing $1$ and $\xi$.
		Define a $\mathbb{Q}(\alpha,\beta)$-linear
		map $h:\mathbb{R}\to\mathbb{R}$ by
		$$h(\xi)=2\xi,\quad h(b)=b\quad\text{for}\quad b\in\mathcal{B}\setminus\{\xi\}.$$
		Then $h(1)=1$, but $h$ is not
		real-linear, since
		$h(\xi)=2\xi\neq\xi h(1)$.
		Define
		$$\Phi(A)=A+
		h(\frac{a_{11}+\overline{a_{11}}}{2})E_{12}$$
		for $A=(a_{ij})\in\mathbb{M}_n(\mathbb{C})$.
		Then $\Phi$ is bijective and satisfies
		\eqref{eq:weighted1}.
		However, $\Phi$ is not real-linear, since
		$\Phi(\xi I)=\xi I+2\xi E_{12}
		\neq\xi\Phi(I)$.
		Moreover, $\Phi$ is not multiplicative,
		because
		$\Phi(I^2)=I+E_{12}
		\neq I+2E_{12}=\Phi(I)^2$.
		
	\end{enumerate}
\end{example}

	\bibliographystyle{amsplain}
	
\end{document}